\pdfoutput=1
\documentclass[12pt]{amsart}
\usepackage{amsmath,amsfonts,amssymb,amsxtra,latexsym,amscd,enumerate,amsthm}
\usepackage{mathtools}
\mathtoolsset{showonlyrefs}
\usepackage{fullpage}
\usepackage{calc}
\usepackage{accents}
\usepackage{cancel}
\usepackage{mathrsfs}
\usepackage{slashed}

\usepackage[utf8]{inputenc}
\usepackage[T1]{fontenc}
\usepackage[english]{babel}
\usepackage[hmarginratio=1:1,body={6.0in,9.0in}]{geometry}
\usepackage{xspace}
\usepackage[pdftex,dvipsnames]{xcolor}
\usepackage[pdftex,hyperindex,pagebackref]{hyperref}
\newcommand{\RbritoSectioning}{section}
\numberwithin{equation}{\RbritoSectioning}
\providecommand{\url}[1]{\texttt{#1}} %

\newcommand\restr[2]{{              %
  \left.\kern-\nulldelimiterspace   %
  #1                                %
  \vphantom{\big|}                  %
  \right|_{#2}                      %
  }}

\DeclareMathOperator{\supp}{supp}
\newcommand{\dif}{\ensuremath{\,\text{{\normalfont d}}}}

\newcommand{\R}{\ensuremath{\mathbb{R}}}
\newcommand{\C}{\ensuremath{\mathbb{C}}}

\newcommand{\Dcal}{\ensuremath{\mathcal{D}}\xspace}

\newcommand{\Hcal}{\ensuremath{\mathcal{H}}\xspace}

\newcommand{\Mcal}{\ensuremath{\mathcal{M}}\xspace}

\newcommand{\Scal}{\ensuremath{\mathcal{S}}\xspace}

    \theoremstyle{plain}
    \newtheorem{thm}[equation]{Theorem}
    \newtheorem{prop}[equation]{Proposition}
    \newtheorem{lemma}[equation]{Lemma}
    \newtheorem{coro}[equation]{Corollary}

    \theoremstyle{definition}
    \newtheorem{defi}[equation]{Definition}
    \newtheorem{remark}[equation]{Remark}

    \newtheorem*{notation}{Notation}

    \newtheorem*{exer*}{Exercise}

\newcounter{numexercs}
{\hrule
  \medskip
  \noindent {\sc Exercícios}
  \begin{list}
    {\textbf{\arabic{numexercs}.}}%
    {
      \usecounter{numexercs}%
      \setlength{\rightmargin}{\leftmargin}
    }
  }
  {\end{list}
  \hrule
  \medskip
}

\newcounter{romanitemnum}
{\begin{list}
  {\emph{\roman{romanitemnum}.}}
  {
    \usecounter{romanitemnum}%
    \setlength{\rightmargin}{\leftmargin}
  }
}
{\end{list}}

\newcounter{alphaitemnum}
{\begin{list}
  {\textbf{\alph{alphaitemnum}.}}
  {
    \usecounter{alphaitemnum}%
    \setlength{\labelwidth}{0.5em}%
    \setlength{\itemindent}{0em}
  }
}
{\end{list}}

\definecolor{cite_color}{HTML}{00629B}
\definecolor{link_color}{HTML}{182B49}
\hypersetup{
    colorlinks=true,
    linkcolor=link_color,
    filecolor=cite_color,
    urlcolor=cite_color,
    citecolor=cite_color,
    linkbordercolor=cite_color,
    urlbordercolor=cite_color,
    citebordercolor=cite_color,
}
\usepackage{dsfont} %

\usepackage{tikz}
\usetikzlibrary{arrows.meta,patterns}
\usepackage{marvosym}

\usepackage{etoolbox}
\robustify\eqref
\makeatletter
\patchcmd{\@tocline}{#6}{%
  \ifcase #1\or\or\hskip 1em\or\hskip 2em\else\hskip 3em\fi
  #6%
}{}{\PackageError{dkg-wp}{Unable to set table-of-contents indentation}{}}
\patchcmd{\thebibliography}{\labelsep .5em}{\labelsep 1em}{}%
  {\PackageError{dkg-wp}{Unable to set bibliography label spacing}{}}
\makeatother

\theoremstyle{definition}
\newtheorem*{convention}{Convention}

\NewDocumentEnvironment{nalign}{o}{
    \IfValueTF{#1}{\begin{equation*}\tag{#1}}{\begin{equation}}
    \begin{aligned}
}{
    \end{aligned}
    \IfValueTF{#1}{\end{equation*}}{\end{equation}}
    \ignorespacesafterend
}

\newcommand{\1}{\mathds{1}}
\renewcommand{\S}{\mathbb{S}}
\renewcommand{\epsilon}{\varepsilon}
\renewcommand{\Re}{\,\mathrm{Re}}
\renewcommand{\Im}{\,\mathrm{Im}}
\DeclareMathOperator{\vol}{vol}

\newcommand{\arctanh}{\text{arctanh}\,}

\newcommand{\vartilde}[1]{\widetilde{\raisebox{0pt}[0.85\height]{$\widetilde{#1}$}}}

\newcommand{\phase}{\pm i\Lambda}

\newcommand{\diag}{\mathrm{diag}}

\newcommand{\Cin}{C^{\text{wp}}}

\newcommand{\err}{\mathrm{Err}}

\newcommand{\wavehom}{\phi^{\rm{hom}}}
\newcommand{\waveinhom}{\phi^{\rm{inhom}}}
\newcommand{\phimain}{\phi_{\rm{main}}}
\newcommand{\phiap}{\phi_{\rm{app}}}
\newcommand{\phif}{\phi_{\rm{free}}}
\newcommand{\phit}{\phi_{\rm{tmlk}}}
\newcommand{\psiap}{\psi_{\rm{app}}}
\newcommand{\psia}{\psi_{\rm{asymp}}}
\newcommand{\Gin}{G^{\rm{in}}}
\newcommand{\Gout}{G^{\rm{out}}}
\newcommand{\Din}{D^{\rm{int}}}

\newcommand{\Fren}{F_{\text{rad}}}
\newcommand{\Sm}{S_{\rm{main}}}
\newcommand{\T}{\mathscr T}

\newcommand{\bad}{\Psi_{\mathrm{bad}}^{\pm}}

\renewcommand{\wp}{\Psi}
\renewcommand{\H}{\mathbb{H}}

\title{Modified scattering and asymptotic completeness for a massive
Dirac--wave system in dimension three}
\author{Vitor Borges}
\address[V. Borges]{Department of Mathematics, University of California, San Diego, La Jolla 92093, USA}
\email{vborges@ucsd.edu}

\begin{document}

\begingroup
\raggedbottom
\makeatletter
\patchcmd{\@maketitle}{\topskip42\p@}{\topskip24\p@}{}%
  {\PackageError{dkg-wp}{Unable to compact title spacing}{}}
\patchcmd{\@setauthors}{\@topsep30\p@}{\@topsep18\p@}{}%
  {\PackageError{dkg-wp}{Unable to compact author spacing}{}}
\patchcmd{\@maketitle}{\dimen@34\p@}{\dimen@18\p@}{}%
  {\PackageError{dkg-wp}{Unable to compact space below title}{}}
\patchcmd{\@starttoc}{\addvspace{32\p@\@plus14\p@}}{}{}%
  {\PackageError{dkg-wp}{Unable to remove trailing contents space}{}}
\makeatother

\begin{abstract}
 We establish global well-posedness and modified scattering for the
 three-dimensional Dirac--Klein--Gordon system with a massive spinor and a
 massless scalar field, for small, sufficiently regular and localized initial
 data. This answers a question posed by Dong, Li, and Yuan in
 \cite[Case IV]{DongLiYuan2023}. We also identify a Majorana-type condition on
 the initial data that yields linear scattering, showing that both linear and
 modified scattering solutions coexist in this model.

 Conversely, we establish asymptotic completeness for small asymptotic profiles
 satisfying suitable regularity and decay assumptions. The forward argument
 combines the hyperboloidal foliation method with testing by wave packets,
 whereas the backward construction relies on a contraction argument in suitably
 chosen norms. 
\end{abstract}
\maketitle

\begingroup
\setlength{\baselineskip}{13pt}
\setlength{\linespacing}{13pt}
\hypersetup{linkcolor=black,linkbordercolor=link_color,pdfborder={0 0 1}}
\tableofcontents
\endgroup
\newpage
\endgroup

\section{Introduction}\label{sec:intro}

The Dirac--Klein--Gordon system \eqref{gdkg} is a classical relativistic model
for the Yukawa interaction between a Dirac field of mass $M\geq0$ and a
spin-zero bosonic field of mass $m\geq0$ \cite{BjDr}.
Variants of this system also arise in simplified descriptions of Higgs--fermion
and meson--nucleon interactions \cite{AitchisonHey2012}. 

In three spatial dimensions, the system takes the form
\begin{nalign}[$\rm{gDKG}$]\label{gdkg}
  &(-i\gamma^{\mu}\partial_{\mu}+M)\psi = \phi\psi\\
  &(\square+m^2)\phi = \langle\psi,\gamma^0\psi\rangle,
\end{nalign}
where $M,m\geq0$ are the Dirac and scalar masses,
$\psi:\R_{t,x}^{1+3}\to\C^4$ is a spinor,
$\phi:\R_{t,x}^{1+3}\to\R$ is a real scalar field, and
$\square=\partial_t^2-\Delta_x$.

In this paper, we study the long-time dynamics in the massive--massless case
$M>0$, $m=0$. Normalizing the Dirac mass to one, we consider
\begin{nalign}[DKG]\label{dkg}
  &(-i\gamma^{\mu}\partial_{\mu}+1)\psi = \phi\psi\\
  &\square\phi = \langle\psi,\gamma^0\psi\rangle
\end{nalign}
with prescribed initial data at $t=0$,
\begin{align}\label{eq:i-cond}
  \psi(0,x)=\psi_0(x),\quad
  \phi(0,x)=\phi_0(x),\quad
  \partial_t\phi(0,x)=\phi_1(x).
\end{align}
Although the scalar equation in \eqref{dkg} is a wave equation,
we retain the name Dirac--Klein--Gordon to emphasize its relation
to the general system \eqref{gdkg}.

\vspace{.8em}
\noindent
Our main results address three fundamental questions:
\begin{enumerate}[(i)]
  \item The global existence of solutions to \eqref{dkg}, for small initial data
    satisfying some mild regularity and spatial decay assumptions,

  \item The asymptotic behavior of the solutions at infinity, where we
    establish modified scattering,

  \item Asymptotic completeness: given the (modified) scattering data at
    infinity satisfying certain high-regularity and decay conditions, we
    construct solutions to \eqref{dkg} whose behavior at
    infinity is given by the prescribed scattering data. 

\end{enumerate}
These correspond to Theorems \ref{main-thm-gwp}, \ref{main-thm-scattering}, and \ref{main-thm-completeness}, respectively.
\vspace{.5em}

The modified scattering in (ii) involves a logarithmic phase
correction in the Dirac component, reflecting the long-range
influence of the scalar field. Related long-range effects have
also been studied in the physics literature in a quantum Yukawa model with a
massive fermion and a massless scalar field \cite{BiswasSemenoff2022}.

\subsection{Previous works on \eqref{gdkg}}

Here we present our results within the broader context of the
well-posedness and scattering literature for the
Dirac--Klein--Gordon system \eqref{gdkg} in the case of general masses $ M,m\geq
0 $ and in different spatial dimensions.
Rather than giving a complete historical account, we focus on the most recent
results and those that have been of particular importance to the present work.

\smallskip
\noindent
\textit{Three spatial dimensions.}
The literature on \eqref{gdkg} is quite vast.
Early global results go back to Chadam--Glassey \cite{ChGl}, who, in the
case $ M,m > 0 $, established global well-posedness and scattering for very
regular initial data satisfying certain algebraic (Majorana-type) conditions.

The best general local well-posedness results were established by
D'Ancona, Foschi and Selberg \cite{DaFoSe} for initial data
\begin{align}
  \psi_0&\in H^s(\R^3),\quad
  (\phi_0,\phi_1)\in H^{s+\frac12}(\R^3)
  \times H^{s-\frac12}(\R^3),\quad s>0,
\end{align}
arbitrarily close to the scaling-critical space,\footnote{The scaling-critical
  space $ L_{ x }^{ 2 } \times H_{ x }^{ 1 /2 } \times H_{ x  }^{ - 1 /2 }  $
  corresponds to $ s = 0 $.} 
and by Tesfahun \cite{Tes}, who extended \cite{DaFoSe} away from the
scaling-critical line\footnote{Tesfahun \cite[Theorem 1]{Tes} obtains local well-posedness
  for initial data $ (\psi_0,\phi_0,\phi_1)\in H_{ x }^{ s }\times H_{ x }^{ r
  } \times H_{ x }^{ r-1} $ provided $s>0$ and
  \begin{align}
    \max\left\{\frac{1}{2}+\frac{s}{3},\frac{1}{3}+\frac{2s}{3},s\right\}
    &<r<\min\left\{\frac{1}{2}+2s,1+s\right\}.
  \end{align}
}.
These results treat the massive contributions $ M\psi $ and $ m^2\phi $ as
perturbations of a system of wave equations and 
hold for any $ M,m\geq 0 $. 
Local well-posedness in the scaling-critical space $ L_{ x }^{ 2 } \times H_{ x
}^{ 1 /2 } \times H_{ x }^{ -1 /2 } $ remains open.

Global theories, on the other hand, are sensitive to the mass parameters $ M,m
$ and mostly concern \textit{small} initial data. 
In the \textit{nonresonant} case $ 2M > m > 0 $, Bejenaru and Herr
\cite{BHdkg} established global well-posedness and scattering for small data in
$ H_{ x }^{ \epsilon } \times H^{ 1 /2 + \epsilon }_{ x } \times H_{ x }^{ - 1
/2 + \epsilon}$. 
Around the same time, Wang \cite{Wa} established scattering for initial data
at the critical Besov regularity with one additional angular derivative in the
cases $ M = m = 0 $ and $ 2M > m > 0 $.
The extra angular regularity assumption acts as a null condition, see
\cite{St}, and improves the behavior of near-resonant (parallel) interactions.
Later, Candy and Herr \cite{CandyHerr2018Transference} established scattering in the massive resonant
case $ m\geq 2M > 0 $ in the critical Sobolev space with some small amount of
angular regularity.

For high-regularity, spatially localized initial data, Zhang \cite{Zh}
established small-data global well-posedness and scattering when \(m=1\), with
estimates uniform in the Dirac mass \(M\in[0,1]\). In particular, this includes
the mixed wave–Klein–Gordon configuration \((M,m)=(0,1)\).

The remaining mixed configuration \((M,m)=(1,0)\) is the subject of the present
paper. To the best of our knowledge, small-data global well-posedness for
\eqref{dkg} had not been established prior to the present work. 
Closely related results concern the pseudoscalar system with the same mass
configuration
\begin{nalign}\label{eq:pseudo}
  &(-i\gamma^{ \mu }\partial_{ \mu }+1)\psi = i\phi\gamma^{ 5 }\psi\\
  &(\square\phi + m^2) = i\langle \psi, \gamma^{ 0 }\gamma^{ 5 }\psi\rangle,
\end{nalign} 
where $ \gamma^{ 5 } := i\gamma^{ 0 }\gamma^{ 1 }\gamma^{ 2 }\gamma^{ 3 } $.
For this modified system, Bachelot \cite{Ba} showed scattering for small,
highly regular, and compactly supported initial data. Later, in \cite[Case III]{DongLiYuan2023}, Dong, Li, and Yuan
improved on Bachelot's result by establishing scattering, and novel uniform
energy bounds, for localized and small initial data.

In three dimensions \eqref{dkg} is \textit{charge critical} 
and the problem of global $ C^{ \infty } $ regularity for \textit{large}
initial data is still largely open, and is only understood for special classes
of solutions \cite{Ba-large, CandyHerr2018Majorana}, with conditional scattering
results obtained by Candy--Herr \cite{CandyHerr2018Conditional}.
Most notably, in the case $ (M,m) = (0,1) $, Cai--Dong--Li--Zhao \cite{CaDoLiZh} have very recently
established global existence and linear scattering for initial data with large scalar but small
spinor weighted norms.

\smallskip \smallskip
\noindent
\textit{Two spatial dimensions.}
In two dimensions \eqref{gdkg} is \textit{charge subcritical} and
Grünrock--Pecher \cite{GrPe} established global well-posedness for initial data of arbitrary
size in the space $ L_{ x }^{ 2 } \times H_{ x }^{ 1 /2  } \times H_{ x }^{ -1
/2 } $ using the conservation of mass $ \int \lvert \psi (t)\rvert ^2\dif x =
\textit{constant} $ for the Dirac component. 
Local well-posedness is known for initial data in $ H_{
x}^{ s } \times H_{ x }^{ r } \times H_{ x }^{ r-1 }  $ in the range
$s>-\frac{1}{5}$,
\begin{align}
  \max\left\{\frac{1}{4}-\frac{s}{2},\frac{1}{4}+\frac{s}{2},s\right\}
  &<r<\min\left\{\frac{3}{4}+2s,\frac{3}{4}+\frac{3s}{2},1+s\right\},
\end{align}
see \cite[Theorem 1]{DaFoSe-2D}.
There is a gap between the best known local well-posedness regularity
and the scaling-critical space $ H_{ x }^{ -\frac{1}{2}  } \times L _{ x }^{ 2
} \times H_{ x }^{ -1 }  $; however, as noted in \cite{DaFoSe-2D}, this  
is expected in dimension two. 
The sharp well-posedness threshold is currently open.
Unconditional uniqueness is known in the range $ s\geq 0 $, $ r = s +
\frac{1}{2}  $, by the work of Selberg and Tesfahun \cite{SeTe}.

In the mixed case $ (M,m) = (0,1) $, small-data scattering was first obtained
by Dong--Wyatt \cite{DoWy} for compactly supported initial data using the
hyperboloidal foliation method.
Under the same mass assumption, Dong, Li, Ma, and Yuan
\cite{DoLiMaYu}, and Zhang \cite{Zh} removed the compact support
assumption from \cite{DoWy} and established small-data scattering in weighted spaces.
In \cite{Zh-2d-large}, Zhang showed scattering for a class of (weighted)
initial data with \textit{large} scalar mass and small spinorial mass.
Scattering for \textit{nonlocalized} small data in $ H_{ x }^{ \frac{1}{2} +
\epsilon } \times H_{ x }^{ 1+\epsilon } \times H_{ x }^{ \epsilon }  $ was
obtained by Bejenaru and the author \cite{BejB} in the nonresonant case $ 2M >
m > 0 $ using Fourier-space methods.

\smallskip \smallskip
\noindent
\textit{Scattering from infinity.}
The literature on the problem of scattering from infinity is more limited.
Hayashi--Ikeda--Naumkin \cite{HayashiIkedaNaumkin2011} established
the existence of a wave operator for \eqref{gdkg} in three spatial
dimensions when $M,m>0$, for sufficiently small final data in weighted
Sobolev spaces.
In two dimensions, Ikeda \cite{Ikeda2013} obtained a corresponding result
under the nonresonance condition $m\neq 2M$, with $M,m>0$, using
a normal form transformation and Strichartz estimates.
For a related system of nonlinear Klein--Gordon equations similar to 
the pseudoscalar system \eqref{eq:pseudo}, Wibowo \cite{Wibowo2009} constructed
wave and scattering operators in three dimensions, again for $ M,m > 0 $
and small data in weighted Sobolev spaces.

\subsection{Main results}

We establish \textit{the first small-data global well-posedness}
result for \eqref{dkg} in three space dimensions.
We also establish energy and dispersive bounds for solutions arising
from \textit{regular and localized} initial data.
This gives a positive answer to a question posed in \cite[Case IV]{DongLiYuan2023}.

\begin{thm}\label{main-thm-gwp}
  Assume that the initial data $ (\psi_0,\phi_0,\phi_1) $ satisfy 
  \begin{align}
    \sum_{ 3j+k \leq 9 }
    \lVert \lvert x \rvert^{ j } \partial_{ x }^{ \leq j+k } \psi_0 \rVert _{ L_{ x }^{ 2 }  }
      + \lVert \lvert x \rvert^{ j } \partial_{ x } ^{ j+k } \phi_0 \rVert_{ \dot H^{
      \frac{1}{2} }_{ x }  } 
    + \lVert \lvert x \rvert^{ j } \partial_{ x } ^{ j+k } \phi_1  \rVert _{ \dot
      H_{ x } ^{ -\frac{1}{2}  } }
    \leq \epsilon 
  \end{align} 
  as well as the low-frequency bound
  \begin{align}\label{eq:low-freq}
    \lVert \phi_0 \rVert_{ \dot H^{ \frac{1}{2} -\nu }_{ x }  }
    + \lVert \phi_1 \rVert_{ \dot H_{ x } ^{ -\frac{1}{2} -\nu } }
    \leq \epsilon,
    \qquad \nu > 0.
  \end{align} 
  Here $j,k$ are nonnegative integers, and $\partial_x^m$ denotes the family
  of all spatial derivatives of order $m$.
  If $ \epsilon $ is small enough, the Cauchy problem \eqref{dkg} admits
  global-in-time solutions satisfying the vector field bounds
  \begin{align}
    \sum_{ 3\lvert J \rvert + \lvert K \rvert \leq 9 }
      \lVert  \widehat \Omega^{ J } \partial^{ K } \psi (t)\rVert _{ L_{ x }^{ 2 }  }
      + \lVert \Omega^{ J } \partial ^{ K } \phi (t)\rVert_{ \dot H^{ \frac{1}{2}
      }_{ x }  } 
      + \lVert \Omega^{ J }  \partial ^{ K } \partial_{ t } \phi (t) \rVert _{
      \dot H_{ x } ^{ -\frac{1}{2} } }
    \lesssim \epsilon \langle t\rangle^{ C\epsilon }
  \end{align} 
  as well as the pointwise bounds
  \begin{align}
    \lvert \psi(t,x) \rvert
    &\lesssim \frac{ \epsilon }
      { \langle t+ \lvert x \rvert \rangle^{ \frac{3}{2} - C\epsilon }
        \langle t-\lvert x \rvert\rangle^{ C\epsilon }},
    \qquad
    \lvert \phi(t,x) \rvert
    \lesssim
    \frac{\epsilon}{\langle t+\lvert x \rvert  \rangle}. 
  \end{align}
  Furthermore, there exists some $ \delta > 0 $ such that, if $ \lvert x
  \rvert  - t \geq \langle t \rangle  ^{          \frac{1}{3} } $, 
  then $ \psi $ satisfies the improved estimate 
  \begin{align}
    \lvert \psi(t,x) \rvert
    &\lesssim \frac{ \epsilon }
      { \langle t+ \lvert x \rvert \rangle^{ \frac{3}{2} } \langle t\rangle^{
      \delta }  }.
  \end{align} 
\end{thm}
  
Our second result concerns asymptotics. 
Unlike in the other mass configurations in \eqref{gdkg}, small solutions to
\eqref{dkg} do not scatter linearly and instead undergo \textit{modified
scattering}.
  
\begin{thm}[Modified scattering]\label{main-thm-scattering}
  There exists $ \delta > 0 $ such that, for all solutions to \eqref{dkg} in 
  Theorem \ref{main-thm-gwp}, there exist asymptotic profiles
  \begin{align}
    (\rho^{ \pm }_{ \infty },\phi_{ \infty })
    &\in C^{ \frac{1}{2}  }\left( 
      \overline{B(0,1)}
    \right), 
  \end{align} 
  vanishing at the boundary, so that inside the light cone we have the
  asymptotic expansions
  \begin{align}\label{eq:phi-asympt-mod}
    \phi(t,x) = (t^2-x^2)^{ -\frac{1}{2}  } \phi_{ \infty }(x /t) 
    + O \left( 
      \epsilon\langle t \rangle  ^{ -1 }\langle t-r \rangle  ^{ -\delta } 
    \right) 
  \end{align} 
  and
  \begin{align}\label{eq:psi-asympt-mod}
    \psi(t,x)
    &= (t^2-x^2)^{ -\frac{3}{4}  } 
    \sum_{ \pm }e^{ \pm i \sqrt{ t^2-x^2 } \mp \frac{i}{2} \phi_{ \infty }(x /t)
    \log (t^2-x^2)  }\rho_{ \infty }^{ \pm }(x /t) 
    + O(\epsilon\langle t \rangle  ^{ -\frac{3}{2}  }\langle t-r \rangle^{
    -\delta }).
  \end{align} 
  Moreover, the profile $ \phi_{ \infty } $ is uniquely determined by $ \rho_{
  \infty }^{ \pm } $ via the elliptic equation
  \begin{align}\label{eq:cancel}
    (-1-\Delta_{ \H })\phi_{ \infty } = \lvert \rho_{ \infty }^{ - } \rvert_{ H
    }^2  
    - \lvert \rho_{ \infty }^{ + } \rvert _{ H }^2 
  \end{align} 
  where $ \Delta_{ \H } $ is the hyperbolic Laplacian in the Klein-Beltrami disk model.
\end{thm}

\begin{remark}
  The emergence of modified scattering is especially striking  when compared 
  with the two-dimensional results \cite{DoWy, DoLiMaYu, BejB} for 
  \eqref{gdkg}. 
  In the mass configurations treated there, solutions disperse more slowly and
  yet scatter linearly. 
  Thus, the long-range interaction in the case $ (M,m) = (1,0) $ reflects a
  feature of the coupling that is not apparent from the decay rates alone;
  see subsection \ref{sec:major-difficulties} for a more detailed explanation.
\end{remark}

\begin{remark}\label{rmk:wave-scattering-intro}
  The discussion on scattering for the wave term $ \phi $ is more subtle.
  On the one hand, the asymptotic equation \eqref{eq:phi-asympt-mod} displays sharp $ t^{ -1 }
  $-decay along timelike rays, which is weaker than the expected free-wave
  decay $ t^{ -1 }\langle t-|x| \rangle  ^{ -\delta } $.
  However, as noted in \cite{IoPa}, this pointwise decay is compatible with
  linear scattering, provided that we measure scattering in a
  sufficiently weak norm $ \dot H^{ \frac{1}{2} + \epsilon } $, $ \epsilon > 0 $.
  We expand this discussion in Remark \ref{rmk:wave-scattering}.
\end{remark}

Theorems \ref{main-thm-gwp} and \ref{main-thm-scattering} above have direct analogues for the Dirac--Maxwell
(DM) system, see \cite[Theorems 1.1, 1.3]{HIS}, and the reader may wonder what
differentiates the dynamics of the two models. 
Here the difference appears in the equation for the asymptotic wave profile.
Unlike its (DM) counterpart, the source \eqref{eq:cancel} has no fixed sign,
and the two spinorial contributions may cancel.
In particular, if $ \lvert \rho_{ \infty }^{ + } (v)\rvert_{ H }  = \lvert \rho_{
\infty}^{ - } (v)\rvert _{ H }   $ for all admissible velocities $ \lvert  v
\rvert < 1 $, then $ \phi_{ \infty } \equiv 0 $ and the logarithmic phase
correction disappears. In this case, \eqref{eq:psi-asympt-mod} describes
\textit{linear scattering}.

The condition $ \phi_{ \infty }\equiv 0 $ is expressed in terms of the
asymptotic profiles, which are themselves determined by the nonlinear
evolution. 
The next proposition gives an \textit{explicit} sufficient condition on the
Dirac initial condition that leads to linear scattering.

\begin{prop}[Linear scattering for the Majorana class]\label{prop:majorana-scattering}
  Let $ z\in \C $, $ \lvert z \rvert = 1 $.
  If $ (\psi_0,\phi_0,\phi_1) $ are initial data satisfying the conditions of
  Theorem \ref{main-thm-gwp}, and the spinorial
  component $ \psi_0 $ satisfies the Majorana-type condition
  \begin{align}
    \psi_{ 0 } + z\gamma^{ 2 }\psi_0^{ * } = 0,
  \end{align}
  where $ \psi_{ 0 }^{ * } $ denotes the complex conjugate of $ \psi_0 $,
  then: 
  \begin{enumerate}[(i)]
    \item The nonlinearity $ \langle \psi, \gamma^{ 0 }\psi \rangle = 0 $ for
      all $ t>0 $, and

    \item the Dirac asymptotic profiles $ \rho_{ \infty }^{ \pm } $ satisfy the
      algebraic condition
      \begin{align}\label{eq:algebraic-intro}
        \rho_{ \infty }^{ \pm  }(v) = -z\gamma^{ 2 }\left( 
          \rho_{ \infty }^{ \mp } (v)
        \right)^{ * },
        \quad v\in B(0,1).
      \end{align} 
  \end{enumerate}
  In particular, $ \phi $ is a free wave, $ \phi_{ \infty }\equiv 0 $, 
  and $ \psi $ linearly scatters to a free solution
  \begin{align}
    \psi(t,x)
    &= (t^2-x^2)^{ -\frac{3}{4}  } 
    \sum_{ \pm }e^{ \pm i \sqrt{ t^2-x^2 }} \rho_{ \infty }^{ \pm }(x /t) 
      + O(\epsilon\langle t \rangle  ^{ -\frac{3}{2}  }\langle t-r \rangle^{
      -\delta })
  \end{align}  
  inside the light cone $ \left\{ t \geq  \lvert x \rvert + 1\right\} $.
\end{prop}

We would also like to understand how much of the scattering set is
captured by the Majorana condition. Let $ \Scal $ denote the collection
of initial data $ (\psi_0,\phi_0,\phi_1) $ in the weighted class of
Theorem \ref{main-thm-gwp} for which the corresponding global solutions
$ (\psi,\phi) $ scatter linearly. Proposition
\ref{prop:majorana-scattering} shows that the Majorana class $ \Mcal $
satisfies $ \Mcal \subset \Scal $. Is $ \Scal = \Mcal $?

To answer this question, we study the more general problem of
constructing solutions with prescribed asymptotic states.
For now, we state the result somewhat informally.
The precise assumptions on the asymptotic data are specified in Proposition
\ref{prop:main-back}.

\begin{thm}[Asymptotic completeness, informal statement]\label{main-thm-completeness}
  Let $\rho_\infty^\pm(v)\in V_v^\pm$ be prescribed
  Dirac profiles with sufficiently high regularity and
  decay as $|v|\to1$, and sufficiently small in suitable
  weighted norms.

  Then there exists a global solution $(\psi,\phi)$ to
  \eqref{dkg} whose initial data satisfy the hypotheses
  of Theorem \ref{main-thm-gwp}, and whose Dirac component
  has the modified scattering expansion
  \eqref{eq:psi-asympt-mod} with the prescribed profiles
  $\rho_\infty^\pm$. The phase profile $\phi_\infty$ is
  determined by $\rho_\infty^\pm$ through \eqref{eq:cancel}.
\end{thm}

In particular, we may choose admissible compactly supported final
Dirac data $ \rho_{\infty}^{\pm} $ satisfying
$ \lvert \rho_{\infty}^{+}(v)\rvert_H
= \lvert \rho_{\infty}^{-}(v)\rvert_H $
for all $ v\in B(0,1) $, but not \eqref{eq:algebraic-intro}.
The corresponding solution scatters linearly but lies outside the
Majorana class, giving a negative answer to the question above.

\begin{remark}
  Theorem \ref{main-thm-completeness}, as stated, does not
  establish uniqueness: prescribing only the Dirac profiles
  $\rho_\infty^\pm$ leaves an independent free wave radiation
  term unspecified. 
  We omit this datum here to keep
  the statement simple. Proposition \ref{prop:main-back}
  gives the precise assumptions on the full asymptotic
  data and the corresponding asymptotic condition for
  the free wave component.
\end{remark}

\subsection{Major difficulties and challenges}\label{sec:major-difficulties}

The factorization
\begin{align}
  (\square + 1)\psi
  &=(i\gamma^\mu\partial_\mu + 1)(-i\gamma^\nu\partial_\nu + 1)\psi
\end{align}
allows us to view \eqref{dkg} as a Klein--Gordon-type spinor coupled to a
wave potential. The two fields have different decay rates: $\psi$ is
expected to decay like $t^{-3/2}$ along timelike rays, whereas $\phi$
can have a leading term of order $t^{-1}$. This borderline decay is
already an obstruction to a direct energy argument. Schematically, for
a commuted energy $E(t)$, one obtains 
\begin{align}
  E(t)
  &\lesssim E(1)+\epsilon\int_1^t s^{-1}E(s)\,\dif s
  \lesssim t^{C\epsilon}E(1).
\end{align}
The resulting growth is small, but it prevents us from closing the
sharp decay bootstrap using energy estimates alone.

The Klein--Gordon formulation further exposes an issue regarding loss of
derivatives.
Squaring the Dirac equation in \eqref{dkg} gives 
\begin{align}
  (\square+1)\psi
  &=i\gamma^\mu(\partial_\mu\phi)\psi+(2\phi-\phi^2)\psi.
\end{align}
The term $(\partial\phi)\psi$ places an additional derivative on the
scalar field and breaks the derivative balance between $ \psi $ and $ \phi $. 
This is especially relevant if we compare with the general framework of
\cite{S12} for systems of Klein--Gordon equations, which crucially requires
power-type nonlinearities.

Compared with the other mass configurations, the choice $(M,m)=(1,0)$
combines two difficulties. First, although the scalar source
$\langle\psi,\gamma^0\psi\rangle$ has a spinorial null structure, the positive
Dirac mass prevents cancellations on parallel frequency interactions, see
\cite[Equation 2.3 and Lemma 2.1]{BHdkg}. This is particularly relevant at
high frequencies, where the Klein--Gordon dispersion surface approaches the
wave cone and parallel interactions can be nearly resonant. 
This contrasts with the more favorable null structure in the mass
configurations $(M,m)=(0,1)$ and $(M,m) = (0,0)$ \cite{DongLiYuan2023,Wa}, as well as with the
nonresonant massive regime $2M>m>0$, which prevents these interactions
altogether; see \cite[Lemma 2.2]{BHdkg}. Second, the
coupling in the Dirac equation involves a wave potential $\phi$ with only
$t^{-1}$ decay, rather than the $t^{-3/2}$ decay of the Klein--Gordon
potential in the configuration $(M,m)=(0,1)$. The forcing $\phi\psi$ is thus
borderline for time integration in the energy estimates. The configuration
considered here consequently combines incomplete null cancellation in the
scalar equation with weaker decay in the Dirac coupling. We summarize this
comparison in Figure \ref{fig:mass-comparison}.

\begin{figure}[ht]
\centering
\begin{minipage}[t]{0.48\linewidth}
  \centering
  \begin{align}
    \begin{gathered}
    M=0,\quad m=1 \\[.6em]
    \left\{\begin{aligned}
      &-i\gamma^\mu\partial_\mu\psi
      =\underbrace{\phi\psi}_{\makebox[0pt]{$\scriptstyle\text{expected decay: }
      \lvert \phi \rvert \sim t^{-3/2}$}} \\[10pt]
      &(\square+1)\phi
      =\underbrace{\langle\psi,\gamma^0\psi\rangle}_{\makebox[0pt]{$\scriptstyle\substack{\text{full
      null condition,}\\ \text{cancels parallel resonances} }$}}
    \end{aligned}\right.
    \end{gathered}
  \end{align}
\end{minipage}\hfill
\begin{minipage}[t]{0.48\linewidth}
  \centering
  \begin{align}
    \begin{gathered}
    M=1,\quad m=0 \\[.6em]
    \left\{\begin{aligned}
      &(-i\gamma^\mu\partial_\mu+1)\psi
      =\underbrace{\phi\psi}_{\makebox[0pt]{$\scriptstyle\text{expected decay: }
      \lvert \phi \rvert \sim t^{-1}$}} \\[10pt]
      &\square\phi
      =\underbrace{\langle\psi,\gamma^0\psi\rangle}_{\makebox[0pt]{$\scriptstyle\substack{\text{partial null condition}\\
      \vphantom{\text{cancels parallel resonances}}}$}}
    \end{aligned}\right.
    \end{gathered}
  \end{align}
\end{minipage}
\caption{Null structure and scalar decay in the two mixed mass configurations.}
\label{fig:mass-comparison}
\end{figure}

In the context of the  vector field method one encounters a separate issue:
the scaling field $S=t\partial_t+x\cdot\nabla_x$, which commutes with the free
wave equation, does not commute with the Klein--Gordon propagator. 
As Klainerman noticed in \cite{Kl}, one may avoid this by working with
hyperboloids $H_\Lambda$, where
$\Lambda=(t^2-\lvert x\rvert^2)^{1/2}$. 
This approach was later reintroduced and expanded by LeFloch and Ma \cite{LeF-Ma}
and is referred to as the ``hyperboloidal foliation method''.
These techniques are especially suitable when $ \psi $ is
compactly supported inside the (translated) light cone $ t = \lvert x \rvert +
1$. By finite speed of propagation, they can be adapted to data compactly
supported on $ \left\{ t = 0 \right\} $; see \cite{LeF-Ma,DoWy}, for instance.
For data with spatial decay alone, one must also control the regions near and
outside the cone. 

As we explain in the next subsection, our argument also relies on
the method of \textit{testing by wave packets} of  Ifrim--Tataru
\cite{IfrimTataru3nls,IfrimTataruSurvey}. 
Under this framework, the Klein--Gordon dispersion makes this construction very
sensitive to frequency: the group velocities
$v_\pm(\xi)=\pm\xi/\sqrt{1+\lvert\xi\rvert^2}$ are contained inside a ball of
radius one. In particular, the Klein--Gordon dispersion combines some of the
worst features of both the Schrödinger and wave dispersion relations: 
low-frequency packets move slowly at a ``Schrödinger rate'' and obstruct decay,
while high-frequency packets approach the light cone, where dispersive
estimates degenerate.
In particular, when testing along constant-velocity rays $x=tv$, 
one only enters the ``hyperboloid-admissible region'' $ \{ t \geq
\lvert x \rvert + 1\}  $ once $t\geq(1-\lvert v\rvert)^{-1}$, which tends to
infinity as $\lvert v\rvert\to1$. 
Passing to a limit along each fixed ray is therefore not enough to control the
profiles uniformly up to the cone.
This is especially problematic in the backward problem, where we need to
construct suitable notions of approximate solutions that are defined globally
in space and, in particular, across the light cone.

Finally, the backward problem requires us to combine the timelike asymptotics,
which come from the Klein--Gordon behavior, with an independent free-wave
radiation field that is supported along null directions.
Both pieces must be incorporated into a single
approximate solution and matched across the cone boundary
$\lvert v\rvert=1$. Controlling this transition, together with the
derivatives and residual errors of the approximation, imposes
substantial regularity and decay requirements on the admissible
final data.

\subsection{Outline of the proof}

We now discuss how we overcame the difficulties above and delineate our proof
strategy.
This is a heuristic and fairly informal discussion, and so we deliberately
trade technical precision for a cleaner presentation.

Our basic strategy for the forward problem (Theorems \ref{main-thm-gwp} and \ref{main-thm-scattering}) is
based on the hyperboloidal foliation method \cite{Kl,LeF-Ma} and the method of
testing by wave packets of Ifrim--Tataru
\cite{IfrimTataru3nls,IfrimTataruSurvey} together with its adaptation to
the Dirac--Maxwell system by Herr--Ifrim--Spitz \cite{HIS}. 

\smallskip
\noindent
\textit{Global well-posedness.}
This is based on a bootstrap argument.
The massive and massless waves have different decay rates, which we bring
together by assuming
\begin{align}
  \lVert\psi(t)\rVert_{L_x^6}
  +\lVert\phi(t)\rVert_{L_x^\infty}
  &\leq C\epsilon\langle t\rangle^{-1}
\end{align}
as our bootstrap hypothesis.
Under this assumption, the natural scale-invariant energy $ E[\psi,\phi] $ for
\eqref{dkg} is not uniformly bounded, but instead has a small growth factor in
time
\begin{align}
  E[\psi,\phi] (t) 
  &\lesssim 
  \epsilon \langle t \rangle  ^{ C\epsilon } 
  E[\psi,\phi] (0).
\end{align} 
Higher-order vector field estimates produce the same $ \langle t \rangle  ^{
C\epsilon }  $ factor. 
 Klainerman--Sobolev estimates suffice near
and outside the light cone, but retain this growth in the timelike
interior.

To recover the missing decay, we use testing by wave packets.
Along each ray of constant velocity $x=vt$, $\lvert v\rvert<1$, 
we construct evolving wave packets $ \wp_{ v }^{ \pm }  $ that isolate the two
spinorial modes with oscillations $e^{\pm i\sqrt{t^2-\lvert x\rvert^2}}$.
Testing the solution against $ \wp_{ v
}^{ \pm }  $ produces spinorial amplitudes $ \rho^{ \pm }(t,v) $ satisfying an
approximate ODE:
\begin{align}\label{eq:intro-ode}
  i\partial_t\rho^\pm(t,v)
  &=
  \pm\sqrt{1-\lvert v\rvert^2}\,\phi(t,vt)\rho^\pm(t,v)
  +\text{time-integrable error}.
\end{align}
The coefficient $ \sqrt{ 1-\lvert v \rvert  ^2 } \phi $ is real, so the
leading interaction only affects the phase.
Estimating the error therefore gives bounded amplitudes without the energy
growth factor. Reconstructing the spinor from these profiles recovers the
$L_x^6$ bound. 
Kirchhoff's formula, together with the improved spinor bounds, then gives the
wave estimate and closes the bootstrap.

\smallskip
\noindent
\textit{Identifying the logarithmic phase.}
The ODE \eqref{eq:intro-ode} also explains the logarithmic modification to
scattering: the coefficient $ \phi(t,vt) = O(t^{ -1 }) $ 
can accumulate a logarithmically growing phase after integrating in
time.
After applying a phase correction and taking a limit $ t\to\infty $, we obtain
time-independent amplitudes $ \rho_{ \infty }^{ \pm }(v) $.
To identify this phase correction explicitly, we first need asymptotics for the
wave $ \phi $.

The nonlinearity $ \langle \psi, \gamma^{ 0 }\psi \rangle $ is effectively a
wave source that scales with $(t^2-\lvert x\rvert^2)^{-3/2}$. 
This determines the wave asymptotics to be of the form
\begin{align}
  \phi(t,x) 
  = \frac{1}{\sqrt{ t^2-x^2 } }   \phi_{ \infty }(x /t) + \text{ integrable
  error}.
\end{align} 
Substituting this into the wave equation gives
\begin{align}\label{eq:scattering-intro}
  (-1-\Delta_\H)\phi_\infty(v)
  &=
  \lvert\rho_\infty^-(v)\rvert_H^2
  -\lvert\rho_\infty^+(v)\rvert_H^2,
\end{align}
which fully determines the wave profile $ \phi_{ \infty } $ via an elliptic
equation.
Thus, the leading wave along timelike rays is determined by the spinor
profiles. 
Plugging the wave asymptotic description back into the profile ODE gives the
coefficient $\pm\phi_\infty(v)/t$, and hence the modified scattering expansion
\begin{align}\label{eq:modified-intro}
  \psi(t,vt)
  &\sim
  (t^2-\lvert vt\rvert^2)^{-3/4}\sum_\pm
  e^{\pm i\sqrt{t^2-\lvert vt\rvert^2}\mp i\phi_\infty(v)\log\sqrt{t^2-\lvert
  vt\rvert^2}}
  \rho_\infty^\pm(v).
\end{align}
Along a fixed ray, $\log\sqrt{t^2-\lvert vt\rvert^2}$ differs from $\log t$
by a constant, which we absorb into the limiting profile.

The signed source in \eqref{eq:scattering-intro} also explains why some
solutions scatter linearly. 
The Majorana condition $\psi_0+z\gamma^2\psi_0^\ast=0$, with $\lvert
z\rvert=1$, is preserved by the evolution and, asymptotically, exchanges  
 the two spinorial modes $ \rho_{ \infty }^{ \pm } $, producing the relation
\begin{align}
  \rho_\infty^+(v)
  &=-z\gamma^2\bigl(\rho_\infty^-(v)\bigr)^\ast
\end{align}
between the limiting profiles. 
Consequently, $ \lvert \rho_\infty^{ + } (v) \rvert_{ H } = \lvert
\rho_\infty^{ - }(v) \rvert_{ H }    $ and $\phi_\infty\equiv 0$, causing the
logarithmic phase
to disappear and the spinor to scatter to a free Dirac solution.

\smallskip
\noindent
\textit{Constructing solutions from their asymptotics.}
We now reverse the question: 

\par\medskip
{
\leftskip=1.2em %
\noindent
\textit{Given sufficiently regular, decaying
profiles $\rho_\infty^\pm$, can we find a solution with the expansion above?}
\par %
}

\par\medskip
\noindent
We build these asymptotics into an approximate solution and write
\begin{align}
  \psi&=\psiap+\Psi,
  &\phi&=\phiap+\Phi.
\end{align}

The goal is to solve the resulting equations for $(\Psi,\Phi)$ with zero
final state at $t=\infty$. Our approach draws inspiration from the works
of Ifrim--Tataru \cite{IfrimTataru3nls} on the one-dimensional cubic NLS
and Chen \cite{Chen2024} on the massive Maxwell--Klein--Gordon system.
Adapting these ideas to \eqref{dkg} requires a nonresonant correction
to the asymptotic Dirac ansatz. To recover the forward weighted class,
we also supplement the energy bounds with low-frequency estimates
for the wave component.

For the spinor, the forward expansion \eqref{eq:modified-intro} is the natural
first candidate for $ \psiap $.
However, substituting it into the Dirac equation leaves an error of size
$t^{-1}$ in $L_x^2$, which we correct for by adding a nonresonant term.
Imposing sufficient decay of the profiles as
$\lvert v\rvert\to1$ also lets us extend \eqref{eq:modified-intro} across the
cone and work on complete constant-time slices.
This approach differs from Chen \cite{Chen2024}, which uses a mixed
hyperbolic-flat foliation.

The wave requires another piece of asymptotic data.
Although $\rho_\infty^\pm$ determine the leading timelike profile $\phi_\infty$
through \eqref{eq:scattering-intro}, they do not determine the radiation
escaping to null infinity. We therefore use
\begin{align}
  \phiap&=\phit+\phif,
\end{align}
where $\phit$ is the wave generated by the leading spinor density, with
timelike behavior $(t^2-\lvert x\rvert^2)^{-1/2}\phi_\infty(x/t)$, and
$\phif$ is a free wave
carrying a prescribed additional radiation profile $\Fren$, which we address in
detail in Proposition \ref{prop:main-back}.
Our final data consist of the triplet $(\rho_\infty^+,\rho_\infty^-,\Fren)$.
This decomposition into timelike and null radiation is already present in
\cite{Chen2024}.

With the choice of $ (\psiap,\phiap)$ in place, we replace $t=\infty$ by a large finite
time $T$ and solve the difference equations for $ (\Psi,\Phi) $ with
\begin{align}
  (\Psi_T,\Phi_T,\partial_t\Phi_T)\big|_{t=T}&=0.
\end{align}
A contraction argument gives corrections decaying like $t^{-\delta}$,
uniformly in $T$. To let $T\to\infty$, we compare two such corrections at
a late time and propagate their difference backward to a fixed slice
$t=T_0$. Gr\"onwall costs only a small power of time, which is outweighed
by the decay when $\epsilon$ is small. The corrections therefore converge
at $T_0$; adding back the backgrounds and continuing to $t=0$ produces the
desired solution.

\subsection{Related works}

Here we discuss related work on other dispersive models and their relation to
this paper.
The work of Ionescu--Pausader \cite{IoPa} is one of the most influential works on coupled
wave-Klein--Gordon systems.
Their approach relies on a mix of vector-field and harmonic-analysis
techniques and establishes modified scattering for a large class of
semilinear/quasilinear systems. 
Our approach to the forward problem differs from theirs in a couple of key
aspects.
First, we work directly at the level of the Dirac equation, instead of the
Klein--Gordon equation, which requires addressing the spinorial structure and
understanding the propagation along the subspaces $ V^{ \pm }  $, which are to
be defined later. 
Second, while our arguments require the same number of vector fields as in
\cite{IoPa}, we require fewer standard derivatives.
This derivative improvement comes from the Klainerman-Sobolev theory developed
in \cite{HIS}.

For the massive Maxwell--Dirac system (DM), Cho--Lee \cite{ChoLee-MD} obtained
global existence and modified scattering by combining a  spacetime resonance analysis \cite{GermainMasmoudiShatah2009} with
the vector-field methods of \cite{IoPa}. 
Their work establishes an alternative approach to modified scattering of (DM)
and was obtained independently and around the same time as \cite{HIS}.
This suggests that a combination of vector-field estimates and spacetime
resonance analysis may also provide an alternative approach to \eqref{dkg}.
Such an adaptation would need to account for the spinorial interaction symbols.

Related results on asymptotic completeness include
the work by Chen--Lindblad
\cite{ChenLindblad2023} for wave--Klein--Gordon systems and by He \cite{He2021}
for the massless Maxwell--Klein--Gordon system in Lorenz gauge.
In \cite{Ouyang2023}, Ouyang constructed modified wave operators from small
weighted Sobolev final data for a quasilinear wave-Klein--Gordon model.

Dong--Guo--Li \cite{DongGuoLi-KGZ} have recently observed that, as in
Theorem \ref{main-thm-scattering}, the two-dimensional Klein--Gordon--Zakharov system
\begin{align}
  &(\square+1)E = -nE\\
  &\square n = \Delta |E|^2
\end{align}
admits both scattering and modified scattering solutions.
In fact, the scattering set of this model can be fully characterized in
terms of the initial data:
for small, smooth, localized data, $n$ scatters linearly in energy,
whereas $E$ scatters linearly when $\int_{\R^2}n_1\dif x=0$, and
modified scatters otherwise. 
Here $n_1$ denotes the initial wave velocity.
This criterion follows from the derivative structure of the wave source; the
leading behavior of $ E $ is explicitly determined by $ \int n_1 \dif x $.

An initial-data characterization of scattering for \eqref{dkg} was also one of the
motivations for the present work. Unlike \cite{DongGuoLi-KGZ}, the coefficient
of the logarithmic phase in \eqref{eq:modified-intro} for \eqref{dkg} is purely
nonlinear, as expressed by \eqref{eq:scattering-intro}.
The Majorana condition (Proposition \ref{prop:majorana-scattering})  gives a sufficient condition for
scattering, and the backward theorem (Theorem \ref{main-thm-completeness}) shows that the Majorana condition does
not provide a full characterization of the scattering set. 
We leave such a characterization of the scattering set in terms of the initial
data as an interesting open question.

\addtocontents{toc}{\protect\setcounter{tocdepth}{1}}
\subsection*{Organization of the paper}

This paper is organized as follows.
In section \ref{sec:notation} we introduce most of the relevant notation for
the forward problem. 
In section \ref{sec:prelim} we define our bootstrap assumption and briefly address some
basic results for linear Dirac equations from the setup in \cite{HIS}. 
The arguments in sections \ref{sec:wp}, \ref{sec:wavebounds} and \ref{sec:scattering} comprise
the proof of the bootstrap argument and are divided into the wave packet
analysis for the Dirac equation
(section \ref{sec:wp}), 
analysis of the wave equation
(section \ref{sec:wavebounds}),
and the derivation of the full asymptotic expansions
\eqref{eq:phi-asympt-mod} and \eqref{eq:psi-asympt-mod}
in section \ref{sec:scattering}.
Finally, in section \ref{sec:completeness} we prove asymptotic completeness.
We collect various facts about analysis on
hyperboloids and Sobolev-type estimates in the appendices.

\addtocontents{toc}{\protect\setcounter{tocdepth}{1}}
\subsection*{Acknowledgments}
The author is grateful to Ioan Bejenaru and Cécile Huneau for helpful and insightful
discussions, and to Ioan Bejenaru for his feedback on earlier drafts of this paper.

\addtocontents{toc}{\protect\setcounter{tocdepth}{2}}

\section{Notation}\label{sec:notation}
\addtocontents{toc}{\protect\setcounter{tocdepth}{1}}
We work on $\mathbb R^{1+3}$ with coordinates $x^0=t$ and
$x=(x^1,x^2,x^3)$.  Greek indices range over $0,1,2,3$, Latin indices
range over $1,2,3$, and repeated upper and lower indices are summed.  We
raise and lower indices using
\[
  (g_{\mu\nu})=(g^{\mu\nu})=\diag(-1,1,1,1).
\]
In particular,
\[
  \square
  :=\partial_t^2-\Delta_x
  =-g^{\mu\nu}\partial_\mu\partial_\nu,
  \qquad r=|x|,
  \qquad q=t-r,
  \qquad \langle s\rangle=(1+s^2)^{\frac12}.
\]
We also set
\[
  \langle t-r\rangle_-
  :=
  \left\{
    \begin{aligned}
      1,
      &\quad\text{ if } \lvert t \rvert  \geq  r    \\
      \langle t-r \rangle,
      &\quad\text{ if } \lvert t \rvert  <  r .
    \end{aligned}
  \right.
\]
Our Hermitian product on $\mathbb C^4$ is linear in the first entry:
\[
  \langle \xi,\eta\rangle:=\eta^\dagger\xi,
  \qquad |\xi|^2=\langle\xi,\xi\rangle,
\]
where $ \dagger $ denotes the Hermitian adjoint.

The Dirac matrices $ \left\{ \gamma^{ \mu } \right\}_{ \mu = 0 }^{ 3 } $
associated to the three-dimensional $ \square $ operator are $ 4\times 4 $
complex-valued matrices given by
\begin{align}
  \gamma^0
  &=\begin{pmatrix} I_2 & 0 \\ 0 & -I_2 \end{pmatrix},
  \qquad
  \gamma^j
  =\begin{pmatrix} 0 & \sigma^j \\ -\sigma^j & 0 \end{pmatrix},
  \qquad j=1,2,3,
\end{align}
where $ \left\{ \sigma^{ j }  \right\} _{ j=1,2,3 } $ denote the Pauli matrices 
\begin{align}
  \sigma^1
  &=\begin{pmatrix} 0 & 1 \\ 1 & 0 \end{pmatrix},
  \qquad
  \sigma^2
  =\begin{pmatrix} 0 & -i \\ i & 0 \end{pmatrix},
  \qquad
  \sigma^3
  =\begin{pmatrix} 1 & 0 \\ 0 & -1 \end{pmatrix}.
\end{align}
Here $I_4$ denotes the $4\times 4$ identity matrix.
In particular, this family satisfies
\begin{align}\label{eq:gamma-adjoint-relations}
  (\gamma^0)^2
  &=I_4,
  \qquad
  (\gamma^j)^2=-I_4,
  \qquad
  (\gamma^\mu)^\dagger\gamma^0=\gamma^0\gamma^\mu
  \qquad
  (\gamma^0)^\dagger=\gamma^0,
  \qquad
  (\gamma^j)^\dagger=-\gamma^j,
\end{align}
as well as the anti-commutation relation
\begin{align}\label{eq:gamma-relations}
  \gamma^{ \mu }\gamma^{ \nu } + \gamma^{ \nu }\gamma^{ \mu }
  &= -2g^{ \mu\nu }I_{ 4 }.
\end{align} 
If $ \eta\in \C^{ 4 } $ is a spinor, it is typical in the literature to write
$\bar\eta := \eta^{ \dagger }\gamma^{ 0 }$ so that 
\begin{align}
  \langle \xi, \gamma^{ 0 }\eta \rangle = \bar\eta\xi,
  \quad \xi,\eta\in \C^{ 4 }.
\end{align} 
To avoid possible confusion, we denote complex conjugation 
by $ z^{ \ast } = \Re z - i\Im z $, $ z\in\C $.

If $I$ is a multi-index, then $\partial^I$ denotes a product of
spacetime derivatives.  We use the convention
\[
  \lVert\partial^k u\rVert_Y^2
  :=\sum_{|I|=k}\lVert\partial^Iu\rVert_Y^2,
\]
and similarly for $ \partial^{ \leq k }u $. 
For a wave component, we write
\[
  \phi[s]:=(\phi(s),\partial_t\phi(s)).
\]

\subsection{Commuting vector fields}

For $0\leq\mu<\nu\leq3$, define the Lorentz generators and their
spinorial counterparts by
\begin{equation}\label{def:commuting-vector-fields}
  \Omega_{\mu\nu}
  :=x_\nu\partial_\mu-x_\mu\partial_\nu,
  \qquad
  \widehat\Omega_{\mu\nu}
  :=\Omega_{\mu\nu}+\frac12\gamma_\mu\gamma_\nu,
\end{equation}
which commute with the linear wave and Dirac propagators, respectively.
Moreover, the vector fields $ \Omega $ and $ \widehat \Omega $ interact with the
nonlinearities in \eqref{dkg} via the following product rule, which follows
from a straightforward computation.
\begin{lemma}
  The vector fields $ \{ \Omega_{ \mu\nu }, \widehat \Omega_{ \mu\nu } \}_{ 0\leq
  \mu < \nu \leq 3 } $ commute with the Dirac--Klein--Gordon
  equations \eqref{dkg} and satisfy the product rule
  \begin{align}
    \widehat\Omega_{\mu\nu}(\phi\psi)
    &=(\Omega_{\mu\nu}\phi)\psi
      +\phi\widehat\Omega_{\mu\nu}\psi, \\
    \Omega_{\mu\nu}\langle\psi,\gamma^0\psi\rangle
    &=\langle\widehat\Omega_{\mu\nu}\psi,\gamma^0\psi\rangle
      +\langle\psi,\gamma^0\widehat\Omega_{\mu\nu}\psi\rangle.
  \end{align} 
\end{lemma}

For a spinor and a scalar, respectively, we use the weighted families
\begin{align}\label{def:spinor-gamma-family}
  \Gamma^{\leq k}\psi
  &:=\left\{
    \widehat\Omega^J\partial^I\psi:
    |I|+3|J|\leq k
  \right\},
  \quad
  \Gamma^{\leq k}\phi
  :=\left\{
    \Omega^J\partial^I\phi:
    |I|+3|J|\leq k
  \right\}.
\end{align}
Whenever $\Gamma^{\leq k}$ is applied to the pair $(\psi,\phi)$, the
first convention is used for $\psi$ and the second for $\phi$.
We also let $ \T $ denote normalized derivatives tangent to the light cone,
that is,
\begin{align}
  \T \in 
  \left\{ \partial_{ r }  + \partial_{ t } , \nabla_{\omega} \right\},
\end{align} 
where $\nabla_{\omega}$ denotes normalized angular derivatives.

\subsection{Hyperboloids and spinorial geometry}

Inside the forward light cone, set
\begin{equation}\label{def:hyperbolic-coordinates}
  \Lambda=(t^2-r^2)^{\frac12},
  \qquad
  v=\frac{x}{t}\in B(0,1).
\end{equation}
For $\Lambda>0$, let $ H_{ \Lambda } $ denote the hyperboloid
\begin{align}
  H_\Lambda
  :=\{(t,x):t>0,\ t^2-r^2=\Lambda^2\},
  \qquad
  H := H_{ 1 } 
\end{align}
and, for $0<\Lambda_0<\Lambda_1$, we denote the region between $ H_{ \Lambda_0
} $ and $ H_{ \Lambda_1 } $ by
\begin{align}
  H_{[\Lambda_0,\Lambda_1]}
  :=\{(t,x):t>r,\ \Lambda_0\leq(t^2-r^2)^{\frac12}
    \leq\Lambda_1\}.
\end{align}
We will often identify hyperboloids $ H_{ \Lambda } $ with the Klein ball via
\begin{align}
  \mathbb H:=B(0,1),
  \qquad
  \iota_\Lambda(v)
  :=\frac{\Lambda}{\sqrt{1-|v|^2}}(1,v)\in H_\Lambda.
\end{align}
The hyperbolic distance and volume measure on $ \H $ are denoted by
$d_{H}$ and $\dif V_{H}$; in velocity coordinates,
\begin{equation}\label{eq:hyperbolic-volume}
  d_{ H } (v,w) = 
  \operatorname{arcosh}
  \left( 
    \frac{ 1-v\cdot w }{ \sqrt{ (1-v^2)(1-w^2) }  }  
  \right),
  \quad
  \dif V_{H}(v)=(1-|v|^2)^{-2}\dif v
\end{equation}
for $ v,w\in \H $.
If $\dif V_{H_\Lambda}$ is the volume measure induced by the Minkowski
metric on $H_\Lambda$, then
\begin{equation}\label{eq:hyperboloid-measure-scaling}
  \int_{H_\Lambda}f(v)\,\dif V_{H_\Lambda}(v)
  =\Lambda^3\int_{\H}
  f (\iota_\Lambda^{-1}(w)) \dif V_{H}(w).
\end{equation}
Henceforth, the pull-back above should be understood implicitly, as we will
simply write
\begin{align}
  \int_{ H_{ \Lambda } } f \dif V_{ H_{ \Lambda } }  = \Lambda^{ 3 }\int_{ \H
  }f \dif V_{ H } .
\end{align} 
The invariant measure $ \dif V_{ H }  $ and the Euclidean measure $ \dif
S$ on $ H $ are related via 
\begin{align}
  \dif S 
  &= \sqrt{ t^2 + x^2 } (t^2-x^2) \dif V_{ H } .
\end{align} 
We will often work with truncated hyperboloids, so we define the following
conventions. 
Let $ T > 1 $ be fixed, and let $ 1\leq \Lambda \leq 2T $.
Define
\begin{align}
  H_{ \Lambda, < T } = \left\{ (t,x) \in H_{ \Lambda } : 1\leq t \leq 4T
  \right\} 
\end{align} 
and its pullback $ \H_{ \Lambda, < T } = \iota_{ \Lambda }^{-1}(H_{ \Lambda, < T }) \subset \H
$.

For each $ (t,x) $, define
\begin{equation}
  \gamma^H
  :=\frac{x_\mu\gamma^\mu}{\Lambda}
  =\frac{-\gamma^0+v_j\gamma^j}{\sqrt{1-|v|^2}},
\end{equation}
and the hyperboloidal inner products
\begin{equation}\label{def:H-inner-product}
  \langle\xi,\eta\rangle_H
  :=-\langle\gamma^0\gamma^H\xi,\eta\rangle,
  \qquad
  |\xi|_H^2 := \langle\xi,\xi\rangle_H,
\end{equation}
where $ \xi,\eta\in \C^{ 4 } $.
The notation $ \gamma^{ H }  $ indicates that the matrix $ \gamma^{ H } =
\gamma^{ H } (v) $ only depends on the hyperbolic variable $ v\in \H $, and not
on $ \Lambda $. 
Similarly, for each $ v\in \H $, one should think of the inner product 
$ \langle \cdot \,,\, \cdot \rangle_{ H } =  \langle \cdot \,,\, \cdot
\rangle_{ H(v) } $ as a $ v $-dependent metric on $ \C^{ 4 } $.

We also define the $ L^{ p }  $-spaces, $ 1\leq p\leq\infty $,
\begin{align}
  \lVert f\rVert_{L^p(\mathbb H)}
  &:=\left(\int_{\mathbb H}|f(v)|^p\,\dif V_{H}(v)\right)^{1/p},
\end{align}
for scalar functions, and
\begin{align}
  \lVert\xi\rVert_{L^p(\mathbb H)}
  :=\left(\int_{\mathbb H}|\xi|_H^p
    \,\dif V_{H}\right)^{1/p}
\end{align}
for spinorial functions, with the usual change when $p=\infty$.
The same convention is used on $H_\Lambda$ with
$\dif V_{H_\Lambda}$.

Since $(\gamma^H)^2=I_4$, we introduce
\begin{equation}\label{def:velocity-projectors}
  P_v^\pm:=\frac12(I_4\mp\gamma^H),
  \qquad
  V_v^\pm:=\ker P_v^\pm.
\end{equation}
The spaces $ V_{ v }^{ \pm }  $ correspond to the $ (\pm 1) $-eigenspaces of $
\gamma^{ H }  $, which are orthogonal under the $ H $-inner product 
\[
  \mathbb C^4=V_v^+ \oplus_{ H }  V_v^-,
\]
and $P_v^\pm$ is the $H$-orthogonal projection onto $V_v^\mp$.
We will repeatedly use
\begin{equation}\label{eq:H-Euclidean-duality}
  P_v^+ - P_v^-=-\gamma^H,
  \qquad
  \langle\eta,\gamma^0\xi\rangle
  =\langle\eta,(P_v^+-P_v^-)\xi\rangle_H.
\end{equation}
The projectors $ P_{ v }^{ \pm }  $ play a key role in the construction of the
wave packets for the Dirac equation; see subsection \ref{sec:wp-construction}.

For the sake of completeness, we also define the angular projectors
\begin{equation}\label{def:angular-projectors}
  P_\pm^\theta:=\frac12(I_4\pm\gamma^0\gamma^\theta),
  \qquad
  \gamma^\theta:=\theta_j\gamma^j,
  \qquad
  \theta:=\frac{x}{r}.
\end{equation}
These operators are very useful for studying the pointwise behavior of $ \psi $ near
the light cone.
However, since we suppress the proof of Proposition
\ref{prop:klainerman-sobolev}, their main role in our analysis is simply
to provide an explicit relationship between the norms $\left\lvert \,\cdot\,
\right\rvert _{ H} $ and $ \left\lvert \,\cdot\, \right\rvert  $:
\begin{equation}\label{eq:H-angular-decomposition}
  |\xi|_H^2
  =\left(\frac{t-r}{t+r}\right)^{\frac12}|P_+^\theta\xi|^2
   +\left(\frac{t+r}{t-r}\right)^{\frac12}|P_-^\theta\xi|^2.
\end{equation}
In particular,
\begin{equation}\label{eq:H-Euclidean-comparison}
  \left(\frac{t-r}{t+r}\right)^{\frac14}|\xi|
  \leq |\xi|_H
  \leq
  \left(\frac{t+r}{t-r}\right)^{\frac14}|\xi|.
\end{equation}

\subsection{Conserved quantities for the Dirac equation}
Finally, if
\[
  (-i\gamma^\mu\partial_\mu+1)\psi=F,
\]
then the Dirac current
\[
  J^\alpha[\psi]
  :=\langle\psi,\gamma^0\gamma^\alpha\psi\rangle
\]
satisfies
\begin{equation}\label{eq:dirac-current}
  \partial_\alpha J^\alpha[\psi]
  =-2\Im\langle F,\gamma^0\psi\rangle.
\end{equation}
For the free flow, and more generally when $F=\phi\psi$ with $\phi$ real,
the current is divergence free.  If $n$ is the future unit normal to a
space-like hypersurface $\Sigma$, the corresponding positive flux is
\[
  E_\Sigma(\psi):=-\int_\Sigma n_\alpha J^\alpha[\psi]\dif V_\Sigma.
\]
On $H_\Lambda$ this becomes
\begin{equation}\label{eq:hyperboloidal-dirac-energy}
  E_{H_\Lambda}(\psi)
  =\int_{H_\Lambda}|\psi|_H^2\,\dif V_{H_\Lambda}.
\end{equation}

\addtocontents{toc}{\protect\setcounter{tocdepth}{2}}

\section{Preliminaries}\label{sec:prelim}
\subsection{Bootstrap setup}\label{subsec:bootstrap-setup}

We perform the bootstrap argument not on constant-time slices, but on the
spacetime regions $ C_{ T }  $ defined below.
For dyadic $T\geq1$, define
\begin{align}
  C_T
  &:=\left\{(t,x):T\leq t\leq4T,\quad
    t^2-r^2\leq4T^2\right\},
  \label{def:CT}\\
  C_{<T}
  &:=\left\{(t,x):1\leq t\leq4T,\quad
    t^2-r^2\leq4T^2\right\}.
  \label{def:ClessT}
\end{align}
For $ t\in [1,4T] $, we also write
\begin{align}
  C_{ <T } (t) = \left\{ x : (t,x) \in C_{ <T }  \right\} 
\end{align} 
for a constant-time slice inside $ C_{ <T }  $.

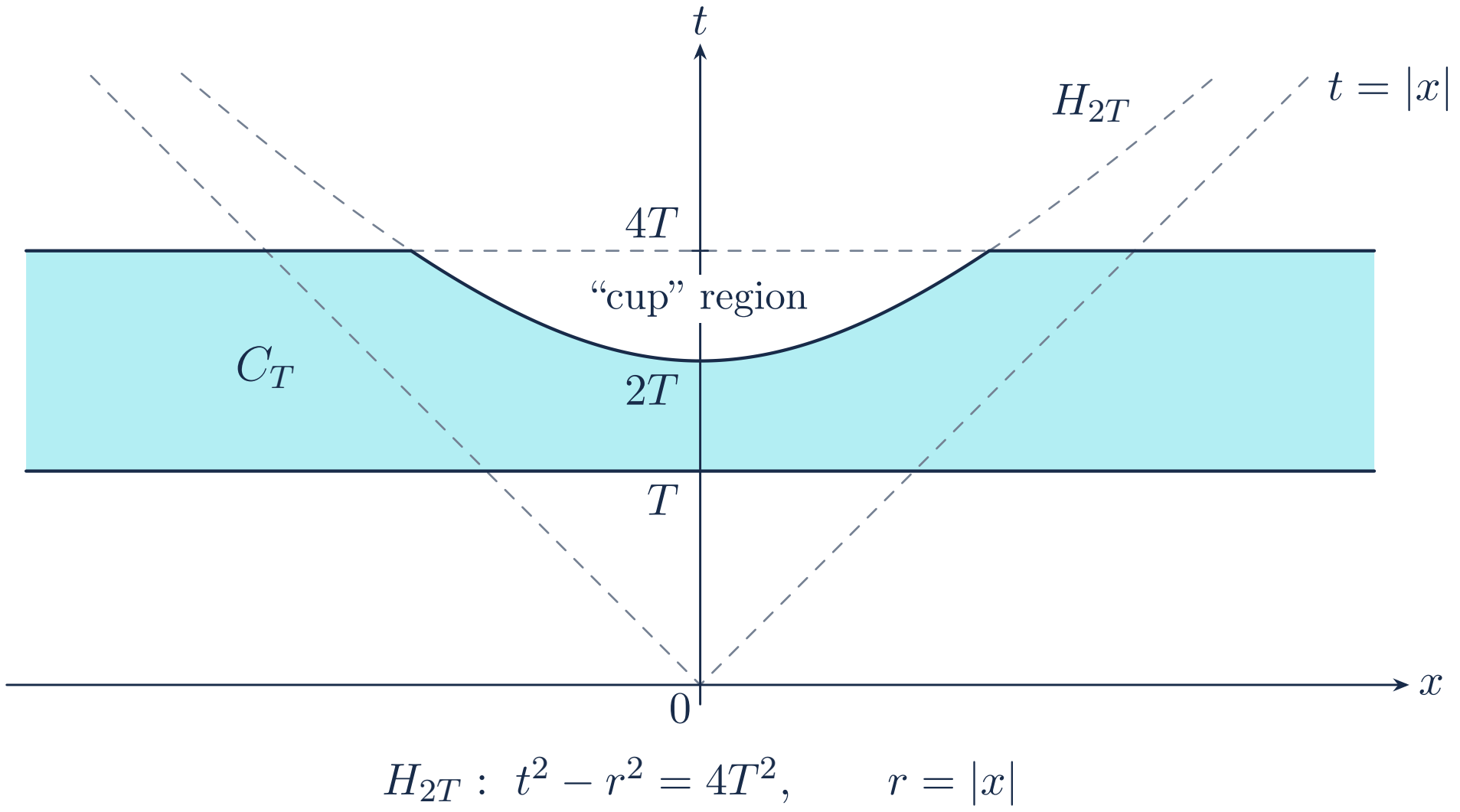
\begin{figure}[h]
  \centering
  \begingroup
\fontsize{10.95}{13.6}\selectfont
\newcommand{\CTTimeT}{1.65}
\newcommand{\CTTimeTwoT}{2.50}
\newcommand{\CTTimeFourT}{3.35}
\definecolor{UCSDNavy}{HTML}{182B49}
\definecolor{UCSDTurquoise}{HTML}{00C6D7}
\providecommand{\CTRegionColor}{UCSDTurquoise!30!white}
\colorlet{ctboundary}{UCSDNavy}
\colorlet{ctfill}{\CTRegionColor}
\colorlet{axisgray}{UCSDNavy}
\colorlet{guidegray}{UCSDNavy!60!white}

\begin{tikzpicture}[
  x=1.12cm,y=1.12cm,
  line cap=round,line join=round,
  axis/.style={draw=axisgray,line width=.55pt,-{Stealth[length=4pt,width=3.2pt]}},
  guide/.style={draw=guidegray,line width=.5pt,dash pattern=on 3pt off 3pt},
  boundary/.style={draw=ctboundary,line width=.8pt},
  every node/.style={text=UCSDNavy,inner sep=2pt}
]
  \pgfmathsetmacro{\crossing}{sqrt(\CTTimeFourT^2-\CTTimeTwoT^2)}

  \fill[ctfill]
    (-5.2,\CTTimeT) -- (5.2,\CTTimeT)
    -- (5.2,\CTTimeFourT) -- (\crossing,\CTTimeFourT)
    -- plot[domain=\crossing:-\crossing,samples=141]
       (\x,{sqrt(\x*\x+\CTTimeTwoT^2)})
    -- (-5.2,\CTTimeFourT) -- cycle;

  \draw[guide] (-4.7,4.7) -- (0,0) -- (4.7,4.7);
  \draw[guide] plot[domain=-4:-\crossing,samples=45]
    (\x,{sqrt(\x*\x+\CTTimeTwoT^2)});
  \draw[guide] plot[domain=\crossing:4,samples=45]
    (\x,{sqrt(\x*\x+\CTTimeTwoT^2)});
  \draw[guide] (-\crossing,\CTTimeFourT) -- (\crossing,\CTTimeFourT);

  \draw[boundary] (-5.2,\CTTimeT) -- (5.2,\CTTimeT);
  \draw[boundary] (-5.2,\CTTimeFourT) -- (-\crossing,\CTTimeFourT);
  \draw[boundary] (\crossing,\CTTimeFourT) -- (5.2,\CTTimeFourT);
  \draw[boundary] plot[domain=-\crossing:\crossing,samples=141]
    (\x,{sqrt(\x*\x+\CTTimeTwoT^2)});

  \draw[axis] (-5.35,0) -- (5.47,0) node[right] {$x$};
  \draw[axis] (0,-.15) -- (0,4.95) node[above] {$t$};
  \node[below left] at (0,0) {$0$};
  \foreach \height/\label/\offset in {
    \CTTimeT/T/-7pt,\CTTimeTwoT/2T/-7pt,\CTTimeFourT/4T/7pt}{
    \draw[axisgray,line width=.5pt] (-.06,\height) -- (.06,\height);
    \node[anchor=east,xshift=-3pt,yshift=\offset] at (0,\height) {$\label$};
  }

  \node[text=ctboundary,font=\fontsize{12}{14}\selectfont] at (-3.35,2.45) {$C_T$};

  \node[font=\fontsize{10}{12}\selectfont,fill=white,inner xsep=3pt,inner ysep=1.5pt]
    at (0,2.98) {``cup'' region};
  \node[anchor=west] at (4.77,4.60) {$t=|x|$};
  \node[anchor=south east] at (3.40,4.28) {$H_{2T}$};

  \node[anchor=north] at (0,-.49)
    {$H_{2T}:\ t^2-r^2=4T^2,\qquad r=|x|$};
\end{tikzpicture}
\endgroup

  \caption{The spacetime region $C_T$.}
\end{figure}

Fix a large constant $C_0$.  
On $C_{<T}$ we make the bootstrap assumption
\begin{equation}\label{eq:bootstrap-assumption}
  \lVert\psi(t)\rVert_{L^6(C_{<T}(t))}
  +\lVert\phi(t)\rVert_{L^\infty(C_{<T}(t))}
  \leq C_0\epsilon \langle t \rangle  ^{ -1 },
  \qquad 1\leq t\leq4T.
\end{equation}
Global well-posedness will then follow from establishing
\begin{align}
  \lVert\psi(t)\rVert_{L^6(C_{<T}(t))}
  +\lVert\phi(t)\rVert_{L^\infty(C_{<T}(t))}
  \leq C'\epsilon \langle t \rangle  ^{ -1 },
  \qquad 1\leq t\leq4T,
\end{align}
with an improved constant $ C'\in (0,C_0) $.
This is the content of Corollary \ref{coro:close-bootstrap-psi} and Proposition
\ref{prop:close-bootstrap-phi}, which are the
main results of sections \ref{sec:wp} and \ref{sec:wavebounds}, respectively.

\subsection{Energy and pointwise estimates}\label{subsec:imported-estimates}

Our proof uses the commuted energy and  
 Klainerman--Sobolev estimates developed in
\cite[Proposition~4.1, Theorem~5.1]{HIS}.  
Since the present paper is centered around the wave packet analysis of
\eqref{dkg}, we keep this section and the justification for the results below
deliberately short;  we refer the reader to \cite[sections 4,5]{HIS} for more
details.

For a spinor $\Psi$ and a real scalar $\Phi$, define the energy norm
\begin{equation}\label{def:energy-space}
  \lVert(\Psi,\Phi)(s)\rVert_{\Hcal^0}^2
  :=\lVert\Psi(s)\rVert_{L_x^2}^2
    +\lVert\Phi(s)\rVert_{\dot H_x^{\frac12}}^2
    +\lVert\partial_t\Phi(s)\rVert_{\dot H_x^{-\frac12}}^2.
\end{equation}
The choice of $L^2\times\dot H^{\frac12}\times\dot H^{-\frac12}$ corresponds
to the scale-invariant space for \eqref{dkg} in three space
dimensions.

We first recall the basic energy estimates for the Dirac and wave equations.
If
\begin{align}
  (-i\gamma^\mu\partial_\mu+1)\Psi
  &=F, \notag\\
  \square\Phi
  &=G,
  \label{eq:linear-equations}
\end{align}
on $C_T$, then
\begin{nalign}\label{eq:linear-energy-estimates}
  \lVert\Psi\rVert_{L_t^\infty L_x^2(C_T)}
  &\lesssim \lVert\Psi(T)\rVert_{L_x^2}
    +\lVert F\rVert_{L_t^1L_x^2(C_T)},\\
  \left\lVert(\Phi,\partial_t\Phi)\right\rVert_{
    L_t^\infty\dot H_x^{\frac12}(C_T)
    \times L_{ t }^{ \infty } \dot H_x^{-\frac12}(C_T)}
  &\lesssim
    \left\lVert(\Phi(T),\partial_t\Phi(T))\right\rVert_{
      \dot H_x^{\frac12}\times\dot H_x^{-\frac12}}
    +\lVert G\rVert_{L_t^1\dot H_x^{-\frac12}(C_T)}.
\end{nalign}
Since the norms above are nonlocal and unwieldy in the regions $ C_{ T }  $, 
we use the local isotropic norm
\begin{align}
  \lVert(\Psi,\Phi)\rVert_{X_T}
  &:=
  \lVert(\Psi,\Phi)(T)\rVert_{\Hcal^0}
  +T^{\frac12}\lVert(-i\gamma^\mu\partial_\mu+1)\Psi\rVert_{L_{t,x}^2(C_T)}
  +T^{\frac13}\lVert\square\Phi\rVert_{L_{t,x}^{\frac32}(C_T)}
\end{align}
which, after applying the Sobolev embeddings $\dot H^{\frac12}(\R^3)\hookrightarrow
L^3(\R^3)$ and $L^{\frac32}(\R^3)\hookrightarrow\dot H^{-\frac12}(\R^3)$, gives
\begin{align}
  \lVert(\Psi,\Phi)\rVert_{L_t^\infty\Hcal^0(C_T)}
  &\lesssim\lVert(\Psi,\Phi)\rVert_{X_T}.
\end{align} 

The following energy estimates follow directly from \cite[Proposition~4.1]{HIS}.
The argument there relies on the quadratic structure of the nonlinearities,
their compatibility with the commuting vector fields, and
Gagliardo--Nirenberg-type interpolation inequalities for vector fields.
It does not use Maxwell-specific gauge constraints or null structure.
Since the required product and commutation properties also hold for
\eqref{dkg}, the analysis carries over without major modifications.

\begin{prop}[Commuted DKG energy estimate]\label{prop:commuted-energy-black-box}
  Let $(\psi,\phi)$ solve \eqref{dkg} on $[0,4T]\times\R^3$ and satisfy the
  bootstrap assumption \eqref{eq:bootstrap-assumption} in $C_{<T}$.
  If $C_0\epsilon$ is sufficiently small, then every dyadic
  $1\leq T_1\leq T$ satisfies
  \begin{align}\label{eq:energy-bound}
    \lVert\Gamma^{\leq9}(\psi,\phi)\rVert_{X_{T_1}}
    &\lesssim
    \langle T_1\rangle^{c\epsilon}
    \lVert\Gamma^{\leq9}(\psi,\phi)(0)\rVert_{\Hcal^0},
  \end{align}
  where $ c\approx C_0 $.
\end{prop}

Since a significant portion of our analysis is done on hyperboloids, it is
convenient to rewrite \eqref{eq:energy-bound} in terms of the hyperboloidal
energy $ E_{ H_{ \Lambda } }(\widehat \Omega^{ \leq 3 }\psi) $.
However, the growth factor $ \langle T \rangle  ^{ c\epsilon }  $ 
prevents global bounds on $ E_{ H_{ \Lambda } } (\widehat\Omega^{ \leq 3 }\psi)
$.
Instead, our regularity analysis in section \ref{sec:scattering} relies on the
truncated-hyperboloid replacement below;
the proof follows directly from \eqref{eq:energy-bound} and the divergence
theorem.

\begin{coro}[Hyperboloidal flux bound]\label{coro:truncated-energy}
  Under the hypotheses of Proposition
  \ref{prop:commuted-energy-black-box}, assume
  $\lVert\Gamma^{\leq9}(\psi,\phi)(0)\rVert_{\Hcal^0}\leq\epsilon$.
  For $0<\Lambda\leq2T$, set
  \begin{align}
    H_{\Lambda,<T}:=H_\Lambda\cap C_{<T}.
  \end{align}
  The corresponding truncated hyperboloidal flux satisfies
  \begin{align}\label{eq:commuted-hyperboloidal-energy}
    E_{H_{\Lambda,<T}}(\widehat\Omega^{\leq3}\psi)
    &\lesssim\epsilon^2 \langle T\rangle^{c\epsilon},
  \end{align}
  where $ c \approx C_0 $.
\end{coro}

The pointwise estimates below are established via trace and interpolation
theorems in suitable subregions of $ C_{ <T }  $.
We emphasize that Proposition \ref{prop:klainerman-sobolev} is
\textit{lossless}, in the sense that it simply propagates the a priori energy
bounds from Proposition \ref{prop:commuted-energy-black-box}.

\begin{prop}[Klainerman--Sobolev bounds]\label{prop:klainerman-sobolev}
  Let $T\geq2$ be dyadic, and let $(\Psi,\Phi)$ solve
  \eqref{eq:linear-equations} on $C_T\cup C_{T/2}$.  
  Assume
  \begin{align}
    \lVert\Gamma^{\leq9}(\Psi,\Phi)\rVert_{X_T}\leq M  .
  \end{align}
  Then, there exists $ \delta > 0 $ independent of $ T $ such that, in $ C_{ T
  } , $
  \begin{align}
    \lvert\Psi(t,x)\rvert
    &\lesssim 
      M\langle t+r\rangle^{-\frac32}
      \langle t-r\rangle_-^{-\delta},
    \label{eq:dirac-global-pointwise}\\
    \lvert\partial\Phi(t,x)\rvert
    &\lesssim 
      M t^{-1 }\langle t-r\rangle^{-\frac12}, \notag\\
    \lvert\T\Phi(t,x)\rvert
    &\lesssim M t^{-\frac32 }. \notag
  \end{align}
  In the interior of the light cone, we have the improved bound
  \begin{align}
    \lvert\T\Phi(t,x)\rvert
    &\lesssim 
      M\langle t\rangle^{-\frac32 }\langle
      t-r\rangle^{-\frac12}
  \end{align}
  along normalized null directions.
  Moreover, the angular components from \eqref{def:angular-projectors} satisfy
  \begin{align}\label{eq:refined-near-cone}
    \lvert P_+^\theta\Psi(t,x)\rvert
    &\lesssim M t^{ -\frac{3}{2} }
    \qquad
    \lvert P_-^\theta\Psi(t,x)\rvert
    \lesssim 
      M(t^2-r^2)^{\frac12}t^{-\frac52 }
  \end{align}
  on $ C_{ T } \cap \left\{ t-r\geq 1 \right\}  $.
\end{prop}

While \eqref{eq:refined-near-cone} is not explicitly stated in \cite[Theorem
5.1]{HIS}, it follows from its interior-cone proof, see \cite[section 5]{HIS}.
For our purposes, we restate \eqref{eq:refined-near-cone} in terms of the $ H
$-norm.
Combining
\eqref{eq:energy-bound}  and \eqref{eq:H-angular-decomposition} gives us
\begin{align}
  \lvert\psi(t,x)\rvert_H
  &\lesssim\epsilon\Lambda^{\frac12}t^{-2+C\epsilon},
  \qquad t-r\geq1.
  \label{eq:refined-H-bounds}
\end{align}
This bound shall be used in section \ref{sec:wp} to close the bootstrap assumption in
$ \psi $.

\section{Wave packet analysis}\label{sec:wp}

Our main goal in this section is to prove the Proposition and Corollary below.
\begin{prop}\label{prop:close-bootstrap-psi}
  Fix $ C_0 $ sufficiently large. There exists $ \epsilon_0=\epsilon_0(C_0)>0 $,
  independent of $ T $, such that for every dyadic $ T\geq1 $ and
  $ 0<\epsilon\leq\epsilon_0 $, if the initial data satisfy the hypotheses of
  Theorem \ref{main-thm-gwp} and the solution satisfies the bootstrap assumption
  \eqref{eq:bootstrap-assumption}, then $ \psi $ satisfies the pointwise bounds
  \begin{align}
    \lvert \psi (t,x) \rvert 
    &\lesssim \epsilon \langle t+r \rangle  ^{ -\frac{3}{2} +C\epsilon } \langle
    t-r\rangle ^{ -C\epsilon } 
  \end{align} 
  on $ \{r-t\lesssim t^{ \frac{1}{3}  }\}\cap C_{ <T } $, and
  \begin{align}
    \lvert \psi(t,x) \rvert 
    &\lesssim  \epsilon \langle t+r \rangle  ^{ -\frac{3}{2} -\delta }
  \end{align} 
  on $ \{ r-t \gtrsim t^{ \frac{1}{3}  }\}\cap C_{ <T }   $, for some $ \delta
  > 0 $ independent of $ T $.
\end{prop}
\noindent
Integrating the pointwise bounds above in $ L^{ 6 } $ closes the
bootstrap argument \eqref{eq:bootstrap-assumption} in $ \psi $.
\begin{coro}\label{coro:close-bootstrap-psi}
  For $ C_0 $ sufficiently large, there exists $ \epsilon_0=\epsilon_0(C_0)>0 $,
  independent of $ T $, such that for every
  dyadic $ T\geq1 $ and $ 0<\epsilon\leq\epsilon_0 $, under the hypotheses of
  Proposition \ref{prop:close-bootstrap-psi}, we obtain
  \begin{align}
    \lVert \psi \rVert_{ L_{ x }^{ 6 } (C_{ <T } (t)) } \leq
    \frac{C_0}{2} \epsilon\langle t \rangle
    ^{ -1 },\qquad 1\leq t\leq4T.
  \end{align} 
\end{coro}

Since the Dirac wave packets play a key role in our analysis, below we recall the
construction of spinorial wave packets adapted to the Dirac equation from
\cite[section 6]{HIS}.
The familiar reader may skip subsection \ref{sec:wp-construction} and refer to
Lemma \ref{lemma:goodwp} for the key properties of the construction.

This section is organized as follows.
In subsection \ref{sec:wp-construction}, we provide a complete construction of
wave packets for the Dirac equation.
In subsection \ref{sec:wp-testing}, we define the Dirac profile via wave packet
testing and establish key estimates for the asymptotic approximation and
asymptotic ODE.
In subsection \ref{sec:refined-psi}, we use the wave packet analysis to improve
the pointwise bounds on $ \psi $ and, finally, in subsection
\ref{sec:closing-psi}, we close the bootstrap argument for the Dirac component of
\eqref{dkg}.

\begin{convention}
  For the remainder of this section, unless otherwise stated, we will work in the interior of the light cone $ \left\{
  t-\lvert x \rvert  \geq 1\right\}  $, where the hyperboloidal time $ \Lambda $
  satisfies
  \begin{align}
    t^{\frac{1}{2}}\leq \Lambda\leq t.
  \end{align} 
  The implicit constants are uniform in this region, with 
  $ v $-dependent weights displayed explicitly.
  The region $ \left\{ t - \lvert x \rvert < 1 \right\}  $ is treated
  separately using the Klainerman--Sobolev bounds in Proposition
  \ref{prop:klainerman-sobolev}, see subsection \ref{sec:closing-psi}.
\end{convention}

\subsection{Wave packets for the Dirac equation}\label{sec:wp-construction}

This section is organized as follows.
\begin{itemize}
  \item First, we construct our spinorial wave packets and show they are approximate
    solutions to the free Dirac equation in an appropriate sense.
  \item Second, we define the spinorial profile $ \rho $ by testing $ \psi $ against a wave
    packet and establish basic estimates on $ \rho $.
  \item Third, we show that $ \rho $ is indeed an appropriate asymptotic profile for $ \psi
    $ by deriving good approximation bounds, and derive the asymptotic ODE for $ \rho $
    along constant-velocity rays.
\end{itemize}

To motivate the definition of the spinorial wave packets, we recall some useful
facts about wave packets for scalar equations. 
Consider the half-Klein-Gordon equations 
$ (i\partial_t \pm \langle D \rangle)u^{ \pm } = 0$
with initial condition $ u_0^\pm$.
For localized initial data, using stationary phase to expand the
solutions $ u^{ \pm } $ along the ray $ x = vt $, $ \lvert v \rvert < 1 $,
yields the asymptotic formula
\begin{align}\label{eq:scalar-asymptotics}
  u^\pm(t,vt)
  &= c_\pm t^{-\frac32}
  e^{\pm it\sqrt{1-|v|^2}}\rho^{ \pm }(v)
  + O_v(t^{-\frac52})
\end{align} 
for asymptotic profiles
\begin{align}
  \rho^{ \pm }(v) = (1-|v|^2)^{-\frac54}\widehat{u}_{ 0 }^{ \pm }\left(\mp
  \frac v{\sqrt{1-|v|^2}}\right)
\end{align}
and constants $ c^{ \pm } $.
This suggests that the associated wave packets take the form
\begin{align}\label{eq:scalar-wp}
  \wp_{ v,\text{scalar} }^{ \pm }(t,x)
  &= \chi
    \left( 
      \frac{ x-vt }{ t^{ \frac{1}{2}  } }   
    \right) 
    e^{ \pm i\sqrt{ t^2-x^2 }  }
\end{align} 
where $ \chi $ is a smooth compactly supported bump function with unit integral.
Testing $ u^{ \pm } $  against $ \wp_{ v } ^{ \pm } $ localizes $ u^{ \pm } $
around the ray $ x=vt $ at the sharp Fourier uncertainty scale, which determines
the rescaling factor $ t^{ -\frac{1}{2} } $.
This is known as the \textit{method of testing by wave packets}, see
\cite{IfrimTataruSurvey,IfrimTataru3nls} for
more details.

Adapting \eqref{eq:scalar-wp} to the Dirac propagator requires a few modifications
to the scalar setup:
First, since the Dirac propagator acts on spinors, we must include  
a spinorial coefficient.
Second, since the Dirac propagator is first order in time but contains
two packets with phases $ e^{ \pm i\sqrt{ t^2-x^2 }  } $,
we cannot expect the two profiles $ \rho^{ \pm } $ to be fully independent.
Thus, we impose the compatibility conditions $ \rho^{ \pm }\in V^{ \pm } $.
This choice is not arbitrary; it arises naturally in the proof of Lemma
\ref{lemma:goodwp}.
\begin{defi}\label{def:wp}
  Given an admissible velocity $ v\in B(0,1) $, we define the \textit{Dirac wave
  packets} $ \wp_{ v }^{ \pm }  $ localized at velocity $ v $ with principal
  direction $ \xi^{ \pm }\in V $ by
  \begin{align}\tag{WP}\label{wp}
    \wp_{ v }^{ \pm } (t,x)
    &:= \chi(y) e^{ \pm i\Lambda } \sigma^{ \pm }(x /t),
    \qquad
    y := \Lambda 
    \left( 
      d_{ H } ( x /t,v) 
    \right)^2 ,
  \end{align} 
  where $ \sigma^{ \pm }_{x /t} \in V_{ x /t } ^{ \pm } $ and $ \sigma_{ v }^{
  \pm} = \xi^{ \pm } $. 
   We recall that $ \Lambda = \sqrt{ t^2-x^2 } $. 
\end{defi}

\noindent
Thus, \eqref{wp} acts as a hyperbolic and spinorial version of
\eqref{eq:scalar-wp}.

One should think of  $\sigma^{ \pm }  $ as follows.  
Consider the trivial vector bundle $ E := \C^{ 4 }\times B(0,1) $ over $ B(0,1) $,
with fibers $ E_{ v }  $ canonically isomorphic to $ \C^{
4} $. Here, $ v\in B(0,1) $.
For each $ v $, the $ v $-dependent matrix $ \gamma^{ H }  $ induces a
decomposition
\begin{align}
  \C^{ 4 } =  V_{ v }^{ + } \oplus  V_{ v }^{ - }
\end{align} 
of the fiber $ E_{ v } $ into two $ \gamma^{ H }  $-invariant subspaces $ V_{ v
}^{ \pm }  $  of dimension $ 2 $.
Such a decomposition is smooth and locally trivializable, and so the disjoint
unions $ V^{ \pm } := \bigsqcup_{ v\in B(0,1) } V_{ v }^{ \pm } $ produce two
well-defined vector bundles over $ B(0,1) $.

Given a fixed base point $ v_0 $ and a fixed direction $ \sigma^{ \pm
}\in V_{ v_0 }^{ \pm }  $, $ \sigma^{ \pm } $ is a local section of the bundle
$ V^{ \pm } $ with $ \sigma^{ \pm }_{ v_0 } = \xi^{ \pm } $.
In practice, we want $ \sigma^{ \pm } $ to satisfy the following properties:
\begin{enumerate}[(i)]
  \item $ \sigma^{ \pm }_{ v_0 }  = \xi^{ \pm } $, 
  \item $ \sigma_{ w } ^{ \pm }\in V_{ w }^{ \pm }  $, and
  \item \label{eq:lorentz-invariance-sigma} 
    $ \widehat{\Omega}_{ w } \sigma^{ \pm }(w) = 0,$ 
\end{enumerate} 
for all $ w $ in a neighborhood of $ v_0 $ in $ B(0,1)$, where the last
property can be satisfied via the Frobenius theorem.
For practical purposes, our construction only depends on properties (i) --
(iii) and they will be assumed henceforth.

\begin{notation}
  We denote the value of a section $\sigma^{\pm}$ at $v\in B(0,1)$
  by either $\sigma^{\pm}_v$ or $\sigma^{\pm}(v)$.
  When we wish to indicate the choice of base point explicitly, we write
  $\sigma^{\pm}_v(w)$ for the value of $ \sigma^{ \pm  } $ at $w\in B(0,1)$ of
  the section chosen with base point $v$ and direction $ \sigma_{ v } ^{ \pm
  }\in V_{ v }^{ \pm}  $ fixed.
\end{notation}

\begin{remark}
  In fact, by taking the maximal such extension and using the uniqueness of the
  Frobenius extension, we obtain a unique \textit{global} $ \widehat\Omega
  $-invariant section $ \sigma^{ \pm } $ satisfying (i) -- (iii). 
  This section is fully determined by the initial condition $ \sigma_{ v_0 }^{
  \pm} =  \xi^{ \pm } $.
\end{remark}

\begin{remark}
  In reality, condition \eqref{eq:lorentz-invariance-sigma} cannot be satisfied
  by all six $ \widehat\Omega $ vector fields simultaneously, since the matrix
  corrections $ \gamma_{ \alpha }\gamma_{ \beta } $ are not directly compatible
  with Frobenius' theorem. 
  Indeed, for $ v = 0 $, one can directly observe that the rotation vector
  fields satisfy  $ \Omega_{ ij } |_{ v = 0 } = 0 $, so
  \eqref{eq:lorentz-invariance-sigma} can only
  be satisfied when $ \sigma^{ \pm }(0) = 0 $.
  This costs us the ability to extend
  an arbitrarily chosen $ \sigma^{ \pm }(v)\in V_{ v } ^{ \pm } $.
  However, the incompatibility can be fixed by working with a slightly different
  family of vector fields, e.g.,
  \begin{align}
    \widehat X_{ 1 } = \widehat\Omega_{ 01 } + \widehat\Omega_{ 31 },
    \qquad
    \widehat X_{ 2 } = \widehat\Omega_{ 02 } + \widehat\Omega_{ 32 },
    \qquad
    \widehat X_{ 3 } = \widehat\Omega_{ 03 }.
  \end{align} 
  The key compatibility property is that 
  \begin{align}
    [ \widehat X_1, \widehat X_2 ] = 0,
    \qquad
    [ \widehat X_1, \widehat X_3 ] = \widehat X_1,
    \qquad
    [ \widehat X_2, \widehat X_3 ] = \widehat X_2,
  \end{align} 
  which avoids the cancellations coming from the rotation vector fields.
  Moreover, the vector fields $ \widehat\Omega $ can be written in terms of $
  \widehat X$ up to zero-order corrections, and so all the relevant estimates for $
  \widehat \Omega $ can be immediately obtained from controlling the analogous
  estimates for $ \widehat X $. 
  Therefore, for simplicity, we harmlessly write the analysis below in terms of
  the more usual vector fields $ \Omega $ and assume $\widehat\Omega\sigma^{
  \pm }  = 0$.

\end{remark}

The definition of the wave packet is consistent with
the linear theory provided that $ \wp_{ v } ^{ \pm } $ is an approximate solution
to
\begin{align}
  \left(-i\gamma^{\mu}\partial_{\mu}+1\right)\wp_{v}^{\pm} \approx
  O(t^{-1})\wp_{v}^{\pm},
\end{align} 
see \cite{IfrimTataruSurvey};
below we adjust Definition \ref{def:wp} accordingly and establish the main
properties of the wave packet construction.

\begin{lemma}\label{lemma:goodwp}
  Let $ v\in B(0,1) $, $ \sigma^{ \pm } \in V^{ \pm } $, and $ \wp_{ v } ^{ \pm }
  $ be as in Definition \ref{def:wp}.
  Then, there exists a correction term $R^{\pm}$ with the same localization as
  $\wp_{v}^{\pm}$,
  \begin{align}
    \lvert R^{\pm}\rvert_{H} \lesssim \Lambda^{-1 /2}\lvert
    \sigma^{\pm}\rvert_{H},
  \end{align}
  such that the corrected wave packet
  \[
    \widetilde\wp_{v}^{\pm} = e^{\pm i\Lambda}\left(\chi(y)\sigma^{\pm} + R^{\pm}\right)
  \] 
  satisfies
  \[
    \left(-i\gamma^{\mu}\partial_{\mu} + 1\right)\widetilde\wp_{v}^{\pm} = \Lambda^{-1} e^{\pm
    i\Lambda} r_{v}^{\pm},
  \] 
  where, again, $r_{v}^{\pm}$ has similar size and localization to $\widetilde\wp_{v}^{\pm}$.
  Moreover, we can represent $r_{v}^{\pm}$ as
  \begin{equation}\label{eq:vfstructure}
    r_{v}^{\pm} = \Lambda^{-\frac{1}{2}}\left(\widehat\Omega\widetilde r_{v}^{\pm}
      + \vartilde r_{v}^{\pm}\right),
  \end{equation} 
  where $\widetilde r_{v}^{\pm}$ and $\vartilde r_{v}^{\pm}$ have the same size and localization
  as $\Psi_{v}^{\pm}$.
\end{lemma}

As preparation, we record various useful projector identities which can be checked
directly.

\begin{lemma}
  Let $P_{v}^{\pm}$ and $\Lambda$ be defined as in \eqref{def:velocity-projectors} and
  \eqref{def:hyperbolic-coordinates}. 
  Then, for each $ \mu
  =0,\ldots,3$,
  \begin{align}
    &\partial_{\alpha}\Lambda = \frac{-x_{\alpha}}{\Lambda},\\
    \label{rel2}&\left[P^{\pm}, \partial_{\mu}\right] 
    = \pm \frac{1}{2}\left(\frac{\gamma_{\mu}}{\Lambda} +
      \frac{x_{\alpha}\gamma^{\alpha}}{\Lambda^3}x_{\mu}\right),\\
      \label{rel3}&\gamma^{\alpha}\left[P^{\pm},\partial_{\alpha}\right] 
    = \mp \frac{3}{2\Lambda} \\
      \label{rel4}&P^{\pm}\gamma^{\alpha} 
    = \gamma^{\alpha}P^{\mp} \pm \frac{x^{\alpha}}{\Lambda}\\
      \label{rel5}&\gamma^{\alpha}(\partial_{\alpha}P^{\pm}) 
    =\pm\frac{3}{2\Lambda}. 
  \end{align}
  Furthermore,  
  $\gamma^{\alpha}\gamma_{\alpha} 
  = g_{\alpha\beta}\gamma^{\alpha}\gamma^{\beta} 
  = -g_{\alpha\beta}g^{\alpha\beta} 
  = -4$.
\end{lemma}

\begin{proof}[Proof of Lemma \ref{lemma:goodwp}]
We construct $ R^{ \pm } $ iteratively by applying the Dirac operator to $\wp_{v}$,
identifying the slowest decaying terms, adding appropriate corrections to the definition of 
$\wp_{v}^{\pm}$, and repeating the process until the remaining terms have
sufficient decay in $ \Lambda $.
By Lorentz invariance, it suffices to consider the case $ v = 0 $.

By direct computation,
\begin{align}
  (-i\gamma^{\mu}\partial_{\mu}+ 1)\wp_{v}^{\pm} 
  &= \chi(y)e^{\phase}\left(1\mp
    \frac{x_{\mu}\gamma^{\mu}}{\Lambda}\right)\left[\sigma^{\pm}\left(x /t\right)\right] 
    -i\chi(y)e^{\phase}\gamma^{\mu}\partial_{\mu}\left[\sigma^{\pm}\left(x
    /t\right)\right]\\
  &\hspace{2.5em} + i\gamma^{\mu}\chi'(y)
    d_{H}^{2}\frac{x_{\mu}}{\Lambda}e^{\phase}\sigma^{\pm}\left(x /t\right)
  - i\gamma^{\mu}\chi'(y)(\partial_{\mu}d _{H}^{2})\Lambda e^{\phase}\sigma^{\pm}
    \left(x /t\right)\\
  &=: R_{1}^{\pm}+R_2^{\pm} + R_3^{\pm}+R_{4}^{\pm}.
\end{align}
Since $ \sigma_{ x/ t }  ^{ \pm }\in V_{ x /t } ^{ \pm } = \ker P_{ x /t } ^{ \pm } $,
\begin{align}\label{eq:wp-compatibility}
  R_{1}^{\pm} = 2\chi(y)e^{\phase}P_{ x /t } ^{\pm}\sigma^{\pm}(x /t) = 0.
\end{align} 
In fact, this cancellation is the reason we require $ \sigma_{ x /t } ^{ \pm
}\in V_{ x /t } ^{ \pm } $.

For the remaining terms, we treat their $V^{ + } $ and $V^{ - }$ components
separately. 
For this first iteration, the $V^{\pm}$ components of $R^{\pm}_{j}$ have the correct
decay rate, whereas some of their $V^{\mp}$ components do not and need to be
incorporated in the definition of $\widetilde\wp_{v}^{\pm}$, as we check below.
We streamline the notation by writing $ P^{ \pm } $ instead of $ P_{ x /t }^{ \pm }  $.

By \eqref{rel4} and \eqref{rel3},
\begin{align}
  P^{\mp}R_{2}^{\pm} 
  &= -i\chi(y)e^{\phase} P^{\mp}\gamma^{\mu}\partial_{\mu}\left[\sigma^{\pm}(x /t)\right] \\
  &= -i\chi e^{\phase}\left(\gamma^{\mu}P^{\pm}\mp
   \frac{x^{\mu}}{\Lambda}\right)\partial_{\mu}\left[\sigma^{\pm}(x /t)\right]\\
  &= -i\chi e^{\pm
    i\Lambda}\left(\gamma^{\mu}P^{\pm}\partial_{\mu}\left[\sigma^{\pm}(x /t)\right]
    \mp \frac{x^{\mu}}{\Lambda}\partial_{\mu}\left[\sigma^{\pm }(x /t)\right]\right) \\
  &= -i\chi e^{\pm i\Lambda} \left(\mp \frac{3}{2\Lambda}\sigma^{\pm}(x /t)\mp
  \frac{x^{\mu}}{\Lambda}\partial_{\mu}\left[\sigma^{\pm}(x /t)\right]\right)\\
  &=\pm \frac{3}{2\Lambda} i\chi e^{\pm i\Lambda} \sigma^{\pm}(x /t),
\end{align}
where, since $ \sigma^{ \pm }(x /t) $ is constant along rays of constant velocity,
$ x^{ \mu }\partial_{ \mu }\sigma^{ \pm }(x /t) = 0 $.
This fact will be used repeatedly when analyzing the remaining terms.

Similarly, we use the fact that $P^{\pm}\sigma^{\pm} = 0$ and  \eqref{rel4} to compute
$P^{\mp}R_{3}^{\pm}$ and $P^{\mp}R_{4}^{\pm}$:
\begin{align}
  P^{\mp}R_{3}^{\pm} 
  &= iP^{\mp}\gamma^{\mu}\chi'(y)d _{H}^2 \frac{x_{\mu}}{\Lambda} e^{\pm
    i\Lambda}\sigma^{\pm}
    \nonumber\\
  &= i \left(\gamma^{\mu}P^{\pm}\mp \frac{x^{\mu}}{\Lambda}\right)\chi'(y)d _{H}^2
    \frac{x_{\mu}}{\Lambda} e^{\pm i\Lambda}\sigma^{\pm}(x /t)
    \nonumber\\
  &= \mp i \frac{x_{\mu}x^{\mu}}{\Lambda^2}\chi'(y)d _{H}^{2} e^{\pm
    i\Lambda}\sigma^{\pm}(x /t)
    \nonumber\\
  &= \pm i\chi'(y) d _{H}^{2}e^{\pm i\Lambda}\sigma^{\pm}(x /t)
    \nonumber\\
  \label{eq:hitschi}&= \pm i\Lambda^{-1}\chi'(y) \underbrace{d _{H}^{2}\Lambda}_{=O(1)}e^{\pm i\Lambda}\sigma^{\pm}(x /t)
\end{align}
so we also conclude that $P^{\mp}R_{3}^{\pm} = R_{3}^{\pm}$, and
\begin{align}
P^{\mp}R_{4}^{\pm} 
&= -i P^{\mp}\gamma^{\mu}\chi'(y) \partial_{\mu}d _{H}^{2}\Lambda
  e^{\pm i\Lambda}\sigma^{\pm}(x /t)\\
&= -i \left(\gamma^{\mu}P^{\pm}\mp \frac{x^{\mu}}{\Lambda}\right)\chi'(y)
  \partial_{\mu}d _{H}^2 \Lambda e^{\pm i\Lambda}\sigma^{\pm}(x /t)\\
&= \pm i\chi'(y) x^{\mu}\partial_{\mu}d _{H}^2 e^{\pm i\Lambda}\sigma^{\pm}(x /t)\\
&= 0.
\end{align}

Therefore,
\begin{nalign}\label{eq:firstiteration}
  (-i\gamma^{\mu}\partial_{\mu}+1)\wp_{v}^{\pm} 
  &=  P^{ \mp }
    \left( 
      R_{ 2 }^{ \pm } + R_{ 3 }^{ \pm } + R_{ 4 }^{ \pm } 
    \right)
   + P^{ \pm }
    \left( 
      R_{ 2 }^{ \pm } + R_{ 3 }^{ \pm } + R_{ 4 }^{ \pm } 
    \right) 
  \\
  &=\underbrace{\pm\frac{i}{2\Lambda}e^{\pm i\Lambda}\left(3\chi(y) + 2\chi'(y)
  d_{H}^{2}\Lambda\right)\sigma^{\pm}(x /t)}_{V^{\pm}\text{ component}} 
  + \underbrace{P_{x /t}^{\pm}\left(R_{2}^{\pm}+R_{4}^{\pm}\right)}_{V^{\mp}\text{
  component}}
\end{nalign} 
where 
\begin{align}
  P_{x /t}^{\pm}\left(
    R_{2}^{\pm} + R_{4}^{\pm}\right)
  &= -ie^{\pm i\Lambda}
    \left(\vphantom{\frac 3\Lambda}\chi(y)\gamma^{\mu}\partial_{\mu}[\sigma^{\pm}(x
    /t)]\right.\\
  &\hspace{6em}
  \left.\pm \frac{3}{2\Lambda}\chi(y)\sigma^{\pm}(x /t)
   + \chi'(y)\gamma^{\mu}\partial_{\mu}d_{H}^{2}\Lambda \sigma^{\pm}(x /t)\right)\\
  &=:\Psi_{\text{bad}}^{\pm}.
\end{align}
Since the $ V^{ \pm } $-components above have sufficient decay $
O(\Lambda^{ -1 }) $, we focus on $ \Psi_{ \text{bad} }^{ \pm } $.
While $P_{{x}/{t}}^{\pm}R_{2}^{\pm}$ gains  a factor of $\Lambda^{-1}$, the term
$P_{{x}/{t}}^{\pm}R_{4}^{\pm}$ only improves by a factor of $\Lambda^{-\frac{1}{2}}$.
Indeed, by Lemma \ref{lemma:metric-expansion},
$d_{H}^{2}(w,0) = \lvert w\rvert ^2 + O(\lvert w\rvert ^{4})$ and
\begin{equation}\label{eq:hitsd}
  \Lambda\left\lvert \nabla_{x}d_{H}^{2}\left(\frac{x}{t},0\right)\right\rvert
  \approx \Lambda \frac{\lvert x\rvert }{t}\cdot t^{-1}
  \approx \Lambda^{\frac{1}{2}}t^{-1}
  \lesssim \Lambda^{-\frac{1}{2}},
\end{equation} 
up to errors of order $ O(\Lambda^{ -\frac{3}{2}  }) $.
We redefine
\begin{align}
  \widetilde\wp_{v}^{\pm} := \wp_{v}^{\pm} - \frac{1}{2}P^{\pm}\left(R_{2}^{\pm} + R_{4}^{\pm}\right)
\end{align} 
and repeat the process, computing $ (-i\gamma^{ \mu }\partial_{ \mu 
}+1)\widetilde\wp_{ v }^{ \pm }  $.

First, we write
\[
  \bad = e^{\pm i\Lambda}P^{\pm}R,
  \quad
  R:= e^{\mp i\Lambda}\left(R_{2}^{\pm}+ R_{4}^{\pm}\right)
\] 
and compute
\begin{align}
  (-i\gamma^{\mu}\partial_{\mu}+1)\left(\frac{1}{2}\bad\right)
  &= \frac{1}{2}e^{\pm i\Lambda}\left(\mp
    \frac{x_{\mu}\gamma^{\mu}}{\Lambda}P^{\pm}R -
    i\gamma^{\mu}\partial_{\mu}\left(P^{\pm}R\right) + P^{\pm}R\right) \nonumber\\
  &= \frac{1}{2}e^{\pm i\Lambda}\left(\left(1\mp
    \frac{x_{\mu}\gamma^{\mu}}{\Lambda}\right)P^{\pm}R
     -i\left(\pm \frac{3}{2\Lambda}\right)R
     -i\gamma^{\mu}P^{\pm}\partial_{\mu}R \right)
   \nonumber\\
  &= \frac{1}{2}e^{\pm i\Lambda}\left(2 P^{\pm}R \mp \frac{3i}{2\Lambda}R
   - i\gamma^{\mu}P^{\pm}\partial_{\mu}R\right) \nonumber\\
  &= \bad \mp \frac{3i}{4\Lambda}(R_{2}^{\pm} + R_{4}^{\pm})
  -\frac{i}{2}\gamma^{\mu}e^{\pm i\Lambda}P^{\pm}\partial_{\mu}R.
  \label{eq:wp-correction}
\end{align}

The first term in \eqref{eq:wp-correction} cancels out the slowly decaying terms in
\eqref{eq:firstiteration},  whereas the second term has improved $O(\Lambda^{-\frac{3}{2}})$
decay.
For the last term, we consider the action of $\gamma^{\mu}P^{\pm}\partial_{\mu}$
separately on each component of $R$:
\[
  \gamma^{\nu}P^{\pm}\partial_{\nu}R
  = \gamma^{\nu}P^{\pm}\partial_{\nu}\left(e^{\mp i\Lambda}R_2^{\pm} 
    + e^{\mp i\Lambda}R_{4}^{\pm}\right).
\] 
  
\noindent
\underline{The $\gamma^{\nu}P^{\pm}\partial_{\nu}(e^{\mp i\Lambda}R_{2}^{\pm})$ component:} 
  When the derivative $\partial_{\nu}$ hits $\chi(y)$, we get terms resembling
  \eqref{eq:hitschi} and \eqref{eq:hitsd}, which improve by at least
  $O(\Lambda^{-\frac{1}{2}})$. Similarly, when the derivative hits either
  $\sigma^{\pm}$ or $\Lambda^{-1}$, the decay rate improves by $O(\Lambda^{-1})$.
  Since $R_{2}^{\pm} = O(\Lambda^{-1})$, we obtain $\gamma^{\nu}P^{\pm}\partial_{\nu}(e^{\mp
  i\Lambda}R_{2}^{\pm}) = O(\Lambda^{-\frac{3}{2}})$.

  \vspace{1em}
  \noindent
  \underline{The $\gamma^{\nu}P^{\pm}\partial_{\nu}(e^{\mp i\Lambda}R_{4}^{\pm})$ component:} 
  The same
  considerations as in the previous case still yield decay $O(\Lambda^{-3 /2})$ for most
  terms; however, since $P^{\pm}R_{4}^{\pm}=R_{4}^{\pm}
  =O(\Lambda^{-\frac{1}{2}})$, the term below
  only decays as $\Lambda^{-1}$:
  indeed,
  \begin{align}
    &-\frac{i}{2}\gamma^{\nu}e^{\pm i\Lambda}P^{\pm}
      \partial_{\nu}\left(e^{\mp i\Lambda}R_{4}^{\pm}\right)
     = -\frac{1}{2}e^{\pm i\Lambda}\gamma^{\nu}P^{\pm}\gamma^{\mu}
      \partial_{\nu}\left[\chi'(y)\Lambda
      \left(\partial_{\mu}d_{H}^{2}\right)\sigma^{\pm}(x/t)\right].
  \end{align}
  Since $\partial_{\nu}y=\Lambda\partial_{\nu}d_H^2-x_{\nu}\frac{ d_{ H }^{ 2 }
  }{ \Lambda } $,
  the terms not already of size $O(\Lambda^{-3/2})$ can be combined into 
  \begin{align}
    R_{5}^{\pm}
    &:= -\frac{1}{2}e^{\pm i\Lambda}
      \gamma^{\nu}P^{\pm}\gamma^{\mu}\left[\chi''(y)\Lambda^2
      \left(\partial_{\nu}d_{H}^{2}\right)\left(\partial_{\mu}d_{H}^{2}\right)\sigma^{\pm}(x
      /t)
        + \chi'(y)\Lambda \left(\partial_{\nu\mu}^{2}d_{H}^{2}\right)\sigma^{\pm}(x /t)
    \right].
  \end{align} 
Since the $V^{\mp}$-component can be corrected for again, we project into $V^{\pm}$. 

To simplify the expression for $R_{5}^{\pm}$, we first note that
\begin{align}
  P^{\mp}\gamma^{\nu}P^{\pm}\gamma^{\mu} 
  = \gamma^{\nu}\gamma^{\mu}P^{\mp} 
  \pm \frac{\gamma^{\nu}x^{\mu}}{\Lambda} 
  \mp \frac{\gamma^{\mu}x^{\nu}}{\Lambda}P^{\mp}
  - \frac{x^{\mu}x^{\nu}}{\Lambda^{2}}.
\end{align} 
Using $x^{\mu}\partial_{\mu}d_H^2=0$ and \eqref{eq:gamma-relations}, we obtain
\begin{align}
  P^{\mp}R_{5}^{\pm}
  &= \frac{1}{2}e^{\pm i\Lambda}g^{\mu\nu}\left[\chi''(y)\Lambda^2
      \left(\partial_{\nu}d_{H}^{2}\right)\left(\partial_{\mu}d_{H}^{2}\right)\sigma^{\pm}(x
      /t) 
    + \chi'(y)\Lambda \left(\partial_{\nu\mu}^{2}d_{H}^{2}\right)\sigma^{\pm}(x
      /t)\right].
\end{align} 

Using the  identities
\begin{align}
  g^{\mu\nu}\partial_{\mu\nu}^2d_{H}^{2}
  = -\square d_H^2
  = 6\Lambda^{-2} + O(\Lambda^{-3})
  \quad\text{and}\quad 
  g^{\mu\nu}\partial_{\mu}d_{H}^{2}\partial_{\nu}d_{H}^{2} =
  4d_{H}^2\Lambda^{-2},
\end{align} 
which are proven in Lemma \ref{lemma:metric-expansion}, we get
\begin{align}
  P^{\mp}R_{5}^{\pm} = \frac{1}{2}\Lambda^{-1}e^{\pm i\Lambda}
  \left[4\chi''(y) d_{H}^{2}\Lambda 
    + 6\chi'(y) 
    \right]\sigma^{\pm}(x /t)+O(\Lambda^{-2}).
\end{align} 
Refining the definition of the wave packet one last time to include the nonresonant
contribution of $R_{5}^{\pm}$, we set
\begin{align}
  \widetilde\wp_{v}^{\pm} :=
  \wp_{v}^{\pm}-\frac{1}{2}P^{\pm}\left(R_{2}^{\pm}+R_{4}^{\pm}\right)
  +\frac{1}{2}P^{\pm}R_{5}^{\pm}.
\end{align} 
Then,
\begin{align}
  (-i\gamma^{\mu}\partial_{\mu}+1)\widetilde\wp_{v}^{\pm}
  = \pm\frac{i}{2\Lambda}e^{\pm i\Lambda}&\left[3\chi(y) + 2\chi'(y)
  d_{H}^{2}\Lambda
  \pm 4i\chi''(y) d_{H}^{2}\Lambda \pm 6i\chi'(y)\right]\sigma^{\pm }(x /t)\\
                                         &\hspace{2em} + O\left(\Lambda^{-\frac{3}{2}}\right).
\end{align}

It remains to prove \eqref{eq:vfstructure}, which we do via a direct computation.  
By Lorentz invariance, we again take $v=0$.  
Using $\widehat\Omega_{0j}\sigma^{\pm}=0$, we have
\begin{align}
  r_v^{\pm}
  &=\Lambda^{-\frac12}\widehat\Omega_{0j}
    \left[\Lambda^{\frac12}\frac{x^j}{\Lambda}
      \left(\pm\frac{i}{2}\chi(y)-\chi'(y)\right)\sigma^{\pm}(x/t)\right]
      \\
  &\qquad\pm\frac{i}{2}
    \left[3\left(1-\frac{t}{\Lambda}\right)\chi(y)
      +\chi'(y)\left(2\Lambda d_H^2-x^j\Omega_{0j}d_H^2\right)\right]
      \sigma^{\pm}(x/t)
      \\
  &\qquad-\left[3\left(1-\frac{t}{\Lambda}\right)\chi'(y)
      +\chi''(y)\left(2\Lambda d_H^2-x^j\Omega_{0j}d_H^2\right)\right]
      \sigma^{\pm}(x/t)
    +O\left(\Lambda^{-\frac12}\right).
\end{align}
We can also rewrite
\begin{align}
  \frac{t}{\Lambda}
  &=\cosh d_H,
  \quad
  x^j\Omega_{0j}d_H^2
  =2\Lambda d_H\sinh d_H.
\end{align}
Since $d_H=d_H(x/t,0)\lesssim\Lambda^{-\frac12}$ on the packet support, we get
\begin{align}
  1-\frac{t}{\Lambda}
  &=O(d_H^2)=O(\Lambda^{-1}), \\
  2\Lambda d_H^2-x^j\Omega_{0j}d_H^2
  &=2\Lambda d_H(d_H-\sinh d_H)=O(\Lambda d_H^4)=O(\Lambda^{-1}),
\end{align}
thus establishing \eqref{eq:vfstructure}.

\end{proof}

For convenience, we shall always denote the corrected wave packet $
\widetilde\wp_{ v,\sigma^{ \pm }}^{ \pm } $ by $ \wp_{ v, \sigma^{ \pm } }^{ \pm } $,
and, whenever the choice of $ \sigma^{ \pm } $ is evident, we drop the subscript
$ \sigma^{ \pm } $ and simply write $ \wp_{ v }^{ \pm} $ instead.

\subsection{Wave packet testing}\label{sec:wp-testing}
We are now ready to define the evolution profile $\rho^{\pm}$ via wave packet testing.
\begin{defi}
  Fix $t > 0$ and an admissible velocity $v\in B(0,1)$. 
  Let $\Lambda > 0$ be such that 
  \[
    (t,vt)\in H_{\Lambda} := \left\{t^2-\lvert x\rvert ^2 = \Lambda^2\right\}.
  \]
  We define the \emph{spinorial profile} $\rho^{\pm}\in V^{\pm}$ 
  in each direction $\sigma_{v}^{\pm}\in V_{v}^{\pm}$ via
  \begin{equation}\label{defi:profile}
      \left\langle \rho^{\pm}(t,v),\sigma_{v}^{\pm}\right\rangle _{H}
    := \int_{ H_{ \Lambda }  }\left\langle \psi,\wp_{v,\sigma_{ v }
      ^{\pm}}^{\pm}\right\rangle _{H}\dif V_{H_{ \Lambda } }
    = \Lambda^{3}\int_{ H }\left\langle \psi,\wp_{v,\sigma_{ v } ^{\pm}}^{\pm}\right\rangle
   _{H}\dif V_{H}
  \end{equation} 
  where $\wp_{v,\sigma^{\pm}}^{\pm}$ is a wave packet with principal direction
  $\sigma_{v}^{\pm}$.
  It is also useful to write $ \rho^{ \pm } $ in terms of the hyperboloidal
  coordinates $ (\Lambda,v) $, so we set
  \begin{align}
    \rho^{ \pm }(\Lambda,v) :=  \rho^{ \pm }(t(\Lambda,v),v),
    \qquad t(\Lambda,v) = \Lambda(1-\lvert v\rvert^2)^{-\frac{1}{2}}.
  \end{align} 
  
\end{defi}

Before proceeding with the analysis of $ \rho $, we describe some basic properties of $ \psi
$ in the $ H $-norm.
Starting from the hyperboloidal energy density
\begin{align}
  e_{H}(\psi) 
  = \frac{t}{\sqrt{t^2+x^2}}\lvert \psi\rvert^2 
    - \frac{\langle \psi,x_{j}\gamma^{0}\gamma^{j}\psi\rangle}{\sqrt{t^2 + x^2}}
  =\frac{t-r}{\sqrt{t^2+r^2}} |P_{ + }^{ \theta } \psi|^2 +\frac{t+r}{\sqrt{t^2+r^2}}
  |P_{ - }^{ \theta }\psi|^2,
\end{align} 
we write
\begin{align}
  |\psi|_H^2
  =\frac{t|\psi|^2-\langle\psi,x_j\gamma^0\gamma^j\psi\rangle}{\Lambda}
  =\frac{\sqrt{t^2+r^2}}{\Lambda}\,e_H(\psi).
\end{align} 
Combining both expressions, we arrive at
\begin{align}\label{eq:Hdecomp}
    \lvert \psi \rvert _{ H } ^2 
    = \left(\frac{t+r}{t-r}\right)^{\frac{1}{2}} \lvert P_-^\theta\psi\rvert^2
     + \left(\frac{t-r}{t+r}\right)^{\frac{1}{2}} \lvert P_+^\theta\psi\rvert^2.
\end{align} 
The proof of \eqref{eq:Hdecomp} also gives us an expression for the hyperboloidal energy
$E_{H}(\psi)$ in terms of the $H$-norm
\begin{equation}\label{eq:EnergyToL2Norm}
  E_{H}(\psi) = \Lambda^{3}\int_{H}\lvert \psi \rvert_{H}^2 \dif V_{H}.
\end{equation} 

\begin{prop}\label{prop:transferbounds}
  The profiles $ \rho^{ \pm }  $ satisfy
  \begin{align}
   &\lvert \rho^{\pm}\rvert_{H}\lesssim \epsilon t^{-2+C\epsilon}\Lambda^2 
     \label{eq:Hnormrho},\\
   &\lVert \widehat\Omega^{\leq 3}
     \rho^{\pm}(\Lambda)\rVert_{L^{2}(\H_{\Lambda,<T})}
   \lesssim \epsilon \langle T \rangle  ^{ C\epsilon }.
   \label{eq:apriorirhoL2}
 \end{align}  
\end{prop}

\begin{proof}
  Since $ \supp\wp_{ v }^{ \pm }  $ has volume $ O(\Lambda^{ -\frac{3}{2}  }) $ 
  for $ \Lambda $ sufficiently large, \eqref{eq:refined-H-bounds} yields
  \begin{align}
    \lvert \rho^{\pm}\rvert_{H} 
    &\lesssim \Lambda^{3} \int_{H}\lvert \psi\rvert_{H}\lvert \wp_{v}^{\pm}\rvert_{H}\dif
    V_{H} \\
     &\lesssim \Lambda^{3}\cdot \epsilon t^{-2+C\epsilon}\Lambda^{\frac{1}{2}}\cdot
     \Lambda^{-\frac{3}{2}}\\
     &\lesssim \epsilon t^{-2+C\epsilon}\Lambda^{2}.
  \end{align}
  
  By Corollary \ref{coro:truncated-energy}, \eqref{eq:apriorirhoL2} follows provided
  we show
  \begin{align}
    \lVert \widehat \Omega^{ \leq 3 }\rho^{ \pm }(\Lambda)
      \rVert_{ L^{ 2 }(\H_{\Lambda,<T}) }
    &\lesssim E_{ H_{ \Lambda, < 2T } }(\widehat \Omega^{ \leq 3 }\psi )^{\frac12}.
  \end{align} 
  We simplify the notation by dropping the $ T $-localization parameters and show
  \begin{align}
    \lVert \widehat \Omega^{ \leq 3 }\rho^{ \pm }(\Lambda)
      \rVert_{ L^{ 2 }(\H) }
    &\lesssim E_{ H_{ \Lambda} }(\widehat \Omega^{ \leq 3 }\psi )^{\frac12}.
  \end{align} 
  
  When no vector fields act on $ \rho $, one can directly check that
  \begin{nalign}\label{eq:L2rhobound}
    |\rho(\Lambda,v)|_H^2
    &\lesssim
    \Lambda^6\left(
      \int_H |\psi(w)|_H |\Psi_v(w)|_H \dif V_H(w)
    \right)^2
    \lesssim \Lambda^{\frac92}
      \int_{\supp\wp_{ v }} |\psi(w)|_H^2 \dif V_H(w).
  \end{nalign}
  Integrating in $ v $ and applying Fubini yields
  \begin{align}
    \lVert\rho(\Lambda)\rVert_{L^2(\H)}^2
    &\lesssim \Lambda^{\frac92}
      \int_H\int_{ \supp \wp_{ v }  }
      |\psi(w)|_H^2 \dif V_H(w) \dif V_H(v) \\
    &\lesssim \Lambda^3\int_H|\psi(w)|_H^2\dif V_H(w)\\
    &= E_{H_\Lambda}(\psi).
  \end{align}

  For the vector field estimates, we first note that the hyperbolic metric satisfies a
  simultaneous Lorentz covariance condition in velocity-space $ \H $:
  \begin{align}
    (\Omega_{ v } + \Omega_{ w } )d_{ H } (v,w)^2 = 0,
  \end{align} 
  which is the hyperbolic analogue of the Euclidean transport equation
  \begin{align}
    (\partial_{ x } + \partial_{ y } )(|x-y|^2) = 0.
  \end{align} 
  In particular, it can be directly checked that the wave packets satisfy
  \begin{align}\label{eq:simultaneous-lorentz}
    (\Omega_{ v } + \widehat \Omega_{ w } )\wp_{ v }^{ \pm } = 0.
  \end{align} 
  Then, by \eqref{eq:simultaneous-lorentz} and integration by parts, 
  \begin{align}
    \left\langle \widehat{\Omega}_v\rho(\Lambda,v),\sigma_v\right\rangle_H
    &=
      \Omega_v\left\langle \rho(\Lambda,v),\sigma_v\right\rangle_H \\
    &= \int_{H_\Lambda}\Omega_v  \left\langle \psi(w),\Psi_v(w)\right\rangle_H
      dV_{H_\Lambda}(w) \\
    &= \int_{H_\Lambda} \left\langle \psi(w), \Omega_v\Psi_v(w)
      \right\rangle_H dV_{H_\Lambda}(w) \\
    &= \int_{H_\Lambda} \left\langle \psi(w), -\widehat{\Omega}_w\Psi_v(w) \right\rangle_H
      dV_{H_\Lambda}(w) \\
    &= \int_{H_\Lambda} \left\langle \widehat{\Omega}_w\psi(w), \Psi_v(w) \right\rangle_H
      dV_{H_\Lambda}(w),
  \end{align}
  which has the same structure as \eqref{eq:L2rhobound} and can be similarly addressed.
  The higher-order estimates follow analogously.
\end{proof}

\begin{lemma}\label{asympdecomp}
  For each $t > 0$ and $v\in B(0,1)$, we have
  \begin{align}\label{eq:dirac-profile-decomposition}
    \psi(t,vt) = \sum_{\pm}\Lambda^{-\frac{3}{2}}e^{\pm i\Lambda}\rho^{\pm}(t,v)
    + O_{H}\left(\epsilon \Lambda ^{-\frac{7}{4}}t^{C\epsilon}\right),
  \end{align} 
  where $ 1\leq \Lambda \leq t $.
\end{lemma}

\begin{proof}
  Decompose $\psi = \psi^{+} + \psi^{-}$ in terms of the subspaces $V_{v}^{\pm}$, and fix an
  $ H $-orthonormal basis $ \left\{ \sigma_{ 1 }^{ \pm } ,\sigma_{ 2 }^{ \pm }
  \right\}  $ of $ V^{ \pm }  $.
  By orthogonality and the definition of the profile, it is enough to show
  \begin{align}
    \langle\psi(t,vt),\sigma_v^\pm\rangle_H
    &=   \Lambda^{ -\frac{3}{2}  } e^{ \pm i\Lambda }\langle \rho^{ \pm },
    \sigma_{ v }^{ \pm }  \rangle_{ H } 
    +O_H\left(\epsilon\Lambda^{-\frac74}t^{C\epsilon}\right),
  \end{align}
  for each unit basis vector $\sigma_v^\pm$.
  We treat the correction and principal terms separately.

  We start with the correction term.
  Since $\lvert R^\pm\rvert_H\lesssim\Lambda^{-\frac12}$ and $R^\pm$
  is supported on a set of volume $O(\Lambda^{-\frac32})$, we obtain from
  \eqref{eq:refined-H-bounds}
  \begin{align}\label{eq:nonresbound}
    \left\lvert\Lambda^{\frac32}\int_H
      \langle\psi,R^\pm\rangle_H\dif V_H\right\rvert
    &\lesssim\epsilon t^{-2+C\epsilon},
  \end{align}
  which is an acceptable error.

  For the principal component, keeping $\Lambda$ fixed, we set
  \begin{align}
    f(w) := \langle\psi(\Lambda,w),\sigma_{w}^{\pm}\rangle_{H},
    \quad
    w\in B(0,1).
  \end{align}
  In geodesic normal coordinates centered at $v$, the kernel
  $\Lambda^{\frac32}\chi(\Lambda\lvert z\rvert^2)$ is a normalized average on a
  ball $ B \subset \H$ of radius $\Lambda^{-\frac12}$.  
   Applying Corollary \ref{coro:rescaledaverage} to a cutoff of $f$ on an
   enlargement of $ B $ and using
   Lemma
  \ref{lem:hyperbolic-flat-average} gives
  \begin{align}
    \left\lvert f(v)-\Lambda^{\frac32}\int_H\chi(y)f\dif V_H\right\rvert
    &\lesssim\Lambda^{-\frac14}
      \lVert\Omega^{\leq2}f\rVert_{L^2(B)}
      +\Lambda^{-1}\lVert f\rVert_{L^\infty(B)}.
  \end{align}
  Then, 
  \eqref{eq:commuted-hyperboloidal-energy} yields
  \begin{align}
    \lVert\Omega^{\leq2}f\rVert_{L^2(B)}
    &\lesssim\Lambda^{-\frac32}
      E_{H_{\Lambda,<2T}}(\widehat\Omega^{\leq2}\psi)^{\frac12}
     \lesssim\epsilon\Lambda^{-\frac32}\langle T\rangle^{C\epsilon}.
  \end{align}
  Additionally, the volume-form error can be controlled by the a priori bounds
  on $ \psi $ \eqref{eq:refined-H-bounds} 
  \begin{align}
    \Lambda^{-1}\lVert f\rVert_{L^\infty(B)}
    &\lesssim\epsilon\Lambda^{-\frac12}t^{-2+C\epsilon}
     \lesssim\epsilon\Lambda^{-\frac74}t^{C\epsilon},
  \end{align}
  completing the proof.
\end{proof}

\begin{lemma}\label{lemma:asymptoticeq}
  For each $v\in B(0,1)$, the profiles $\rho^{\pm}$ satisfy the asymptotic equation
  \begin{equation}\label{eq:asymp}
    i\partial_{t}\rho^{\pm}(t,v) = \pm\sqrt{1-v^2}\phi(t,vt)\rho^{\pm}(t,v)
     + O_{H}(\epsilon t^{-1+C\epsilon}\Lambda^{-\frac{1}{2}}).
   \end{equation} 
\end{lemma}

\begin{proof}
  Fix a principal direction $\sigma_{v}^{\pm}\in V_{v}^{\pm}$, and let
  $\wp_{v}^{\pm}$ be its corresponding wave packet. Consider the current density
  \begin{align}
    J^{\alpha}:=\langle \gamma^{0}\gamma^{\alpha}\psi,\wp_{v}^{\pm}\rangle.
  \end{align}
  Since $(-i\gamma^{\mu}\partial_{\mu}+1)\psi=\phi\psi$, a direct computation gives
  \begin{align}
    \partial_{\alpha}J^{\alpha}
    &= i\phi\langle \gamma^{0}\psi,\wp_{v}^{\pm}\rangle
      -i\left\langle \psi,\gamma^{0}\bigl(-i\gamma^{\mu}\partial_{\mu}+1\bigr)
        \wp_{v}^{\pm}\right\rangle.
  \end{align}
  Fix $ \Lambda_1 > \Lambda_0 > 0 $. 
  Integrating the identity above over $H_{[\Lambda_{0},\Lambda_{1}]}$ and using
  the divergence theorem yields
  \begin{align}
    \left\langle
      \rho^{\pm}(\Lambda_{1},v)-\rho^{\pm}(\Lambda_{0},v),\sigma_{v}^{\pm}\right\rangle_{H}
    &= \int_{ \partial H_{ [\Lambda_0,\Lambda_1] }  } \langle \gamma^{ 0 } \gamma^{ \alpha }
      \psi,\wp_{ v }^{ \pm }  \rangle\dif S\\
    &= \int_{H_{[\Lambda_{0},\Lambda_{1}]}}
      i\phi\langle \gamma^{0}\psi,\wp_{v}^{\pm}\rangle
      - i\left\langle \psi,\gamma^{0} (-i\gamma^{\mu}\partial_{\mu}+1 )
        \wp_{v}^{\pm}\right\rangle\dif x\dif t.
  \end{align}
  Dividing by $ \Lambda_1-\Lambda_0 $, letting $ \Lambda_1\to\Lambda_0 $, and
  using $ \dif \Lambda /\dif t = \sqrt{ 1-v^2 }  $, we get
  \begin{align}\label{eq:rho-derivative-t}
    \partial_t\langle\rho^\pm(t,v),\sigma_v^\pm\rangle_H
    &=+i\sqrt{1-v^2}\,\Lambda^3\int_H
      \phi\langle\gamma^0\psi,\wp_v^\pm\rangle\dif V_H
      + \err
  \end{align}
  where
  \begin{align}
    \err
    &:=-i\sqrt{1-v^2}\,\Lambda^3\int_H
      \left\langle\psi,\gamma^0
      (-i\gamma^\mu\partial_\mu+1)\wp_v^\pm\right\rangle\dif V_H.
  \end{align}

  We first estimate the error term. 
  By Lemma \ref{lemma:goodwp},
  \begin{align}
    (-i\gamma^\mu\partial_\mu+1)\wp_v^\pm
    &=\Lambda^{-\frac32}e^{\pm i\Lambda}
      (\widehat\Omega\widetilde r_v^\pm+\vartilde r_v^\pm),
  \end{align}
  and so we integrate by parts, see \eqref{eq:ibp}, to rewrite
   \begin{align}
    \err
    &= -i\sqrt{ 1-v^2 } \Lambda^{ \frac{3}{2}  } e^{ \mp i\Lambda }
      \int_{ H } \left( -\left\langle \gamma^{ 0 } \widehat\Omega \psi, \widetilde r_{ v
        }^{ \pm } \right\rangle
      + \left\langle \gamma^{ 0 } \psi, {\vartilde r}_{ v }^{ \pm }
        \right\rangle\right)\dif V_{ H }.
  \end{align}
  Then,
  \begin{align}
    \lvert\err\rvert
    &\lesssim \Lambda t^{ -1 }\cdot \Lambda^{ \frac{3}{2}  }\Lambda^{ -\frac{3}{2}
      }\cdot \lVert \widehat\Omega^{ \leq 1 }\psi \rVert_{ L^{ \infty }(\H) } 
    \lesssim \epsilon t^{ -1+C\epsilon }\Lambda^{ -\frac{1}{2}  },
  \end{align}
  where the pointwise bounds on $ \psi $ and $ \widehat\Omega\psi $ follow,
  respectively, from \eqref{eq:refined-H-bounds} and the local Sobolev embedding
  \begin{align}
    \lvert \widehat\Omega\psi (v)\rvert_{ H } 
    &\lesssim \lVert \widehat\Omega ^{ \leq 3 }\psi \rVert_{ L^{ 2 }(B) }
    \lesssim \Lambda^{ -\frac{3}{2}  } E_{ H_{\Lambda, <2t} }(\widehat \Omega^{ \leq
      3}\psi)^{ \frac{1}{2}  }
    \lesssim \epsilon \Lambda^{ -\frac{3}{2}  }\langle t \rangle^{ C\epsilon },
    \quad \Lambda \leq t ,
  \end{align} 
  applied to a ball $ B \subset \H$ of radius $ 1 $ centered at $ v\in \H$,
  together with Corollary \ref{coro:truncated-energy}.
  
  We next freeze $\phi$ at the packet center.  For $ (s,x) $ in the packet support,
  Proposition \ref{prop:klainerman-sobolev}, the Fundamental Theorem of Calculus,
  and $(t-r)^{-\frac12}\lesssim t^{\frac12}\Lambda^{-1}$ give
  \begin{align}\label{eq:frozen-scalar-difference}
    \lvert\phi(s,x)-\phi(t,vt)\rvert
    &\lesssim d_H(x/s,v)\lVert\Omega\phi\rVert_{L^\infty(H)}
    \lesssim\Lambda^{-\frac12}t
      \cdot\epsilon t^{-\frac32+C\epsilon}(t-r)^{-\frac12}
     \lesssim\epsilon\Lambda^{-\frac32}t^{C\epsilon}.
  \end{align}
  Hence, 
  \begin{align}
    &\sqrt{1-v^2}\,\Lambda^3
      \int_H\lvert\phi(s,x)-\phi(t,vt)\rvert
      \lvert\psi\rvert_H\lvert\wp_v^\pm\rvert_H\dif V_H \\
    &\hspace{12em}\lesssim
      \Lambda t^{ -1 }
      \cdot \Lambda^{ 3 }
      \cdot \epsilon\Lambda^{-\frac32}t^{C\epsilon}
      \cdot \epsilon\Lambda^{\frac12}t^{-2+C\epsilon}
      \Lambda^{-\frac32}\\
    &\hspace{12em}\lesssim
    \epsilon t^{-1+C\epsilon}\Lambda^{-\frac12}.
  \end{align}

  It remains to identify the leading integral.  
  By \eqref{eq:H-Euclidean-duality},
  \begin{align}
    \langle \gamma^{0}\psi,\wp_{v}^{\pm}\rangle
    &= \langle \psi,(P^{+}-P^{-})\wp_{v}^{\pm}\rangle_{H}.
  \end{align}
  The principal part of $ \wp_{ v }^{ \pm } $ lies in $V^\pm$, whereas the
  correction has $H$-size
  $O(\Lambda^{-\frac12})$.  Thus, the definition of the profile and
  \eqref{eq:refined-H-bounds} imply
  \begin{align}
    \Lambda^3\int_H\langle\gamma^0\psi,\wp_v^\pm\rangle\dif V_H
    &=\mp\langle\rho^\pm(t,v),\sigma_v^\pm\rangle_H
      +O_H(\epsilon\Lambda^{\frac32}t^{-2+C\epsilon}).
  \end{align}
  Substituting into \eqref{eq:rho-derivative-t} gives
  \begin{align}
    \partial_{t}\left\langle \rho^{\pm}(t,v),\sigma_{v}^{\pm}\right\rangle_{H}
    &= \mp i\sqrt{1-v^{2}}\phi(t,vt)
      \left\langle \rho^{\pm}(t,v),\sigma_{v}^{\pm}\right\rangle_{H}
      + O\bigl(\epsilon t^{-1+C\epsilon}\Lambda^{-1/2}\bigr).
  \end{align}
  Since this holds for every principal direction $\sigma_{v}^{\pm}\in V_{v}^{\pm}$, the
  vector-valued form \eqref{eq:asymp} follows.

\end{proof}

\begin{remark}
  Unlike \cite[Lemma 6.5]{HIS}, our asymptotic equation does not contain a nonresonant
  component with borderline $ t^{ -1 }  $-decay.
  In a sense, this seems to be the manifestation of the partial null condition 
  of \eqref{dkg}.
\end{remark}

\subsection{Refined pointwise bounds for \texorpdfstring{$ \psi $}{psi} and
proof of Proposition \ref{prop:close-bootstrap-psi}}\label{sec:refined-psi}

Now we are ready to prove the pointwise bounds in Proposition
\ref{prop:close-bootstrap-psi}.
Since we analyze $ \psi $ along each ray $ x = vt $, we use the relation $\Lambda
= t\sqrt{ 1-v^2 }$ to rewrite the error terms in Lemmas \ref{asympdecomp} and
\ref{lemma:asymptoticeq} in terms of $ t $ and $ v $.

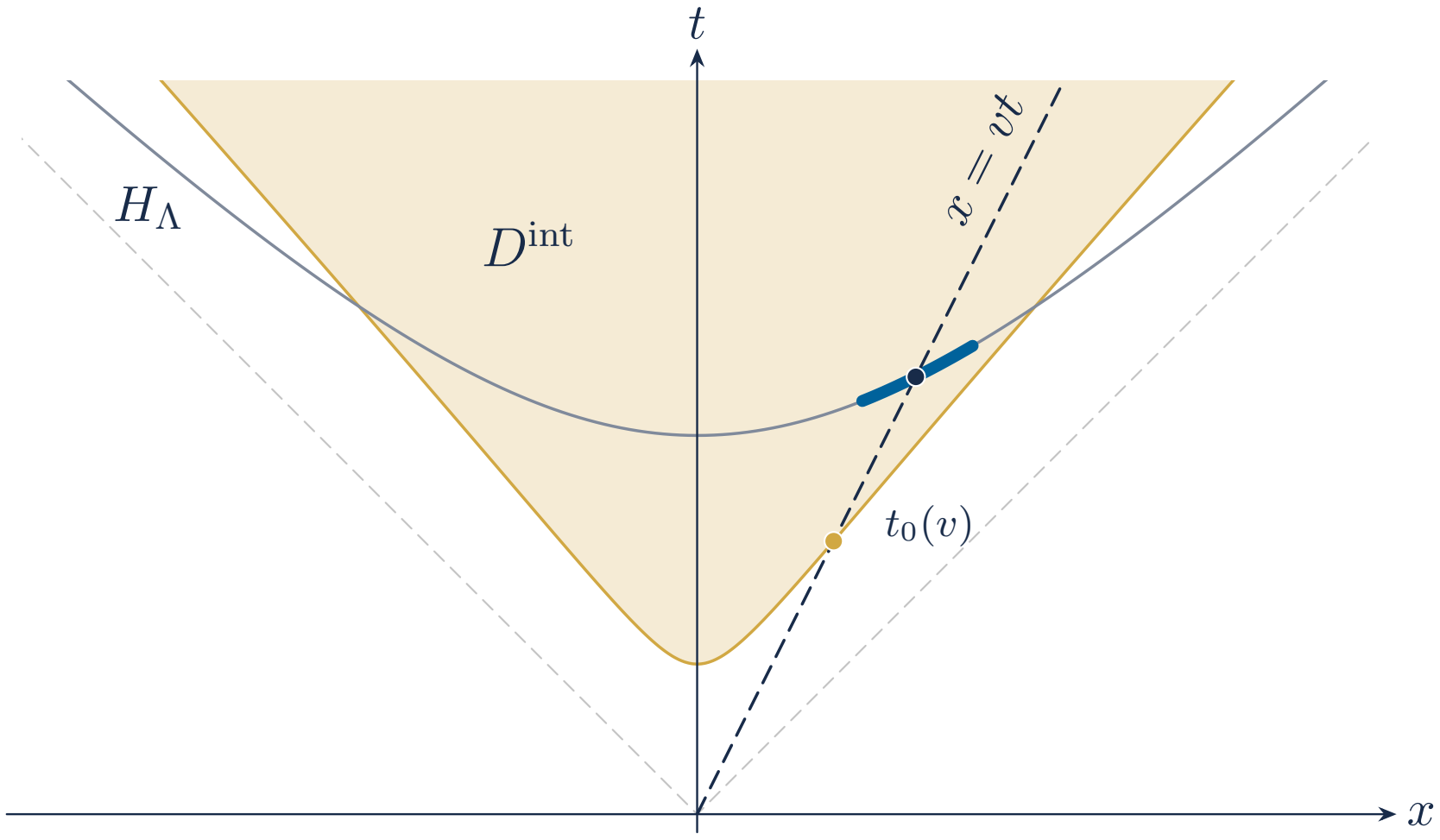
\begin{figure}[h]
  \centering
  \newcommand{\WPTestingVelocity}{0.5}
  \newcommand{\WPTestingSupportScale}{0.2}
  \begingroup
\normalfont\fontencoding{OT1}\fontsize{10.95}{13.6}\selectfont
\providecommand{\WPTestingVelocity}{0.6375}
\providecommand{\WPTestingSupportScale}{0.165}
\providecommand{\WPTestingHyperboloidHeight}{2.5233}
\providecommand{\WPTestingShapeExponent}{2.08}
\providecommand{\WPTestingPlotHeight}{4.885}
\providecommand{\WPTestingPlotHalfWidth}{4.5}
\definecolor{WPTNavy}{HTML}{182B49}
\definecolor{WPTBlue}{HTML}{00629B}
\definecolor{WPTGold}{HTML}{C69214}
\providecommand{\WPTestingRegionColor}{WPTGold!18!white}
\providecommand{\WPTestingHyperboloidColor}{WPTNavy!55!white}
\providecommand{\WPTestingSupportArcColor}{WPTBlue}
\providecommand{\WPTestingThresholdColor}{WPTGold!80!white}
\providecommand{\WPTestingLightConeColor}{black!23!white}
\providecommand{\WPTestingAxisColor}{WPTNavy}
\providecommand{\WPTestingRayColor}{WPTNavy}
\colorlet{WPTRegion}{\WPTestingRegionColor}
\colorlet{WPTHyperboloid}{\WPTestingHyperboloidColor}
\colorlet{WPTSupportArc}{\WPTestingSupportArcColor}
\colorlet{WPTThreshold}{\WPTestingThresholdColor}
\colorlet{WPTLightCone}{\WPTestingLightConeColor}
\colorlet{WPTAxis}{\WPTestingAxisColor}
\colorlet{WPTRay}{\WPTestingRayColor}

\begin{tikzpicture}[
  x=1.15cm,y=1.15cm,line cap=round,line join=round,
  every node/.style={text=WPTNavy,inner sep=2pt},
  axis/.style={draw=WPTAxis,line width=.45pt,-{Stealth[length=4pt,width=3.2pt]}},
  lightcone/.style={draw=WPTLightCone,line width=.4pt,dash pattern=on 3pt off 2pt},
  ray/.style={draw=WPTRay,line width=.65pt,dash pattern=on 4pt off 2.4pt},
  hyperboloid/.style={draw=WPTHyperboloid,line width=.65pt},
  threshold/.style={draw=WPTThreshold,line width=.65pt}
]
  \pgfmathsetmacro{\WPTTop}{\WPTestingPlotHeight}
  \pgfmathsetmacro{\WPTHalfWidth}{\WPTestingPlotHalfWidth}
  \pgfmathsetmacro{\WPTExponent}{\WPTestingShapeExponent}
  \pgfmathsetmacro{\WPTLambda}{\WPTestingHyperboloidHeight}
  \pgfmathsetmacro{\WPTVelocity}{\WPTestingVelocity}
  \pgfmathtruncatemacro{\WPTInvalid}{abs(\WPTVelocity)>=1}
  \ifnum\WPTInvalid=1
    \PackageError{wave-packet-testing}{WPTestingVelocity must satisfy abs(v)<1}
      {Choose a velocity strictly between -1 and 1.}
  \fi
  \pgfmathsetmacro{\WPTUMax}{sqrt(1-pow(\WPTTop+.5,-1/\WPTExponent))}
  \def\ThresholdPath{
    plot[variable=\u,domain=-\WPTUMax:\WPTUMax,samples=241]
      ({\u*pow(1-\u*\u,-\WPTExponent)},
       {pow(1-\u*\u,-\WPTExponent)})}
  \def\HyperboloidPath{
    plot[variable=\x,domain=-\WPTHalfWidth:\WPTHalfWidth,samples=161]
      (\x,{sqrt(\WPTLambda^2+(\x)^2)})}
  \pgfmathsetmacro{\WPTInitialTime}{pow(1-(\WPTVelocity)^2,-\WPTExponent)}
  \pgfmathsetmacro{\WPTTestingTime}{\WPTLambda/sqrt(1-(\WPTVelocity)^2)}
  \pgfmathsetmacro{\WPTTestingX}{\WPTVelocity*\WPTTestingTime}
  \pgfmathsetmacro{\WPTTheta}{.5*ln((1+\WPTVelocity)/(1-\WPTVelocity))}
  \pgfmathsetmacro{\WPTSupportRadius}{\WPTestingSupportScale/sqrt(\WPTLambda)}
  \pgfmathsetmacro{\WPTSupportStart}{\WPTTheta-\WPTSupportRadius}
  \pgfmathsetmacro{\WPTSupportEnd}{\WPTTheta+\WPTSupportRadius}
  \pgfmathtruncatemacro{\WPTLeftRay}{\WPTVelocity<0}
  \ifnum\WPTLeftRay=1
    \def\WPTSide{-1}\def\WPTInitialAnchor{east}
  \else
    \def\WPTSide{1}\def\WPTInitialAnchor{west}
  \fi
  \pgfmathsetmacro{\WPTLabelTransition}{min(1,abs(\WPTVelocity)/.15)}
  \pgfmathsetmacro{\WPTLabelDX}{\WPTSide*.275*\WPTLabelTransition}
  \pgfmathsetmacro{\WPTLabelDY}{
    max(0,.28-.7*(\WPTVelocity)^2)*min(1,abs(\WPTVelocity)/.25)
    -.22*(1-\WPTLabelTransition)}
  \pgfmathtruncatemacro{\WPTReversed}{\WPTInitialTime>=\WPTTestingTime}
  \ifnum\WPTReversed=1
    \PackageWarningNoLine{wave-packet-testing}{The testing point is not after t0(v).
      Increase WPTestingHyperboloidHeight or reduce abs(WPTestingVelocity)}
  \fi
  \pgfmathsetmacro{\WPTSupportExtent}{abs(\WPTTheta)+\WPTSupportRadius}
  \pgfmathtruncatemacro{\WPTOutside}{
    max(\WPTInitialTime,\WPTLambda*cosh(\WPTSupportExtent))>\WPTTop
    || max(abs(\WPTVelocity)*\WPTInitialTime,
            \WPTLambda*sinh(\WPTSupportExtent))>\WPTHalfWidth}
  \ifnum\WPTOutside=1
    \PackageWarningNoLine{wave-packet-testing}{A marker or the support exceeds the frame.
      Increase WPTestingPlotHeight and/or WPTestingPlotHalfWidth}
  \fi
  \coordinate (initial) at ({\WPTVelocity*\WPTInitialTime},\WPTInitialTime);
  \coordinate (testing) at (\WPTTestingX,\WPTTestingTime);
  \coordinate (rayend) at ({\WPTVelocity*\WPTTop},\WPTTop);
  \coordinate (initiallabel) at
    ({\WPTVelocity*\WPTInitialTime+\WPTLabelDX},{\WPTInitialTime+\WPTLabelDY});

  \begin{scope}
    \clip (-\WPTHalfWidth,0) rectangle (\WPTHalfWidth,\WPTTop);
    \fill[WPTRegion] \ThresholdPath -- cycle;
    \draw[lightcone] (-\WPTTop,\WPTTop) -- (0,0) -- (\WPTTop,\WPTTop);
    \draw[threshold] \ThresholdPath;
    \draw[hyperboloid] \HyperboloidPath;
    \draw[draw=WPTSupportArc,line width=2.6pt]
      plot[variable=\a,domain=\WPTSupportStart:\WPTSupportEnd,samples=35]
        ({\WPTLambda*sinh(\a)},{\WPTLambda*cosh(\a)});
    \draw[ray] (0,0) -- (rayend);
  \end{scope}

  \draw[axis] ({-\WPTHalfWidth-.10},0) -- ({\WPTHalfWidth+.16},0) node[right] {$x$};
  \draw[axis] (0,-.12) -- (0,{\WPTTop+.215}) node[above] {$t$};
  \node[anchor=east] at
    ({-.7444*\WPTHalfWidth},{sqrt(\WPTLambda^2+(.7444*\WPTHalfWidth)^2)-.153})
    {$H_\Lambda$};
  \node[font=\fontsize{12}{14}\selectfont]
  at ({-\WPTSide*.2489*\WPTHalfWidth},{.778*\WPTTop}) {$D^{ \rm{int} }$};
  \path (0,0) -- (rayend)
    node[pos=.87,sloped,above=3pt] {$x=vt$};

  \filldraw[fill=WPTThreshold,draw=white,line width=.4pt]
    (initial) circle[radius=2pt];
  \node[font=\fontsize{9}{11}\selectfont,anchor=\WPTInitialAnchor]
    at (initiallabel) {$t_0(v)$};
  \filldraw[fill=WPTNavy,draw=white,line width=.4pt]
    (testing) circle[radius=2pt];
\end{tikzpicture}
\endgroup

  \caption{The setup for testing by wave packets. The dark dashed line is the ray,
  the thick arc is the wave packet support on $H_\Lambda$, and the
  colored region is $ D^{ \rm{int} }:= \bigcup_{ T > 1 } D_{<T}  $, which is where we do the
  wave packet analysis.}
  \label{fig:Cin}
\end{figure}

A priori, the profiles $ \rho^{ \pm } $ satisfy \eqref{eq:Hnormrho},
\begin{align}\label{eq:rho-apriori-growth}
  \lvert \rho^{ \pm }(t,v) \rvert_{ H }  \lesssim \epsilon t^{ C\epsilon } (1-v^2),
\end{align} 
which is not enough to close the bootstrap assumption, see subsection
\ref{sec:closing-psi}.
However, for large $ t $, the dynamics of $ \rho^{ \pm } $ are determined by the
asymptotic ODE \eqref{eq:asymp}:
\begin{align}
  \partial_{ t } | \rho^{ \pm }  |_{ H }  ^2
  &= 2\Re \langle \partial_{ t } \rho^{ \pm } , \rho^{ \pm }  \rangle_{ H } \\
  &= \underbrace{2\Re \left\langle \mp i\sqrt{ 1-v^2 }\phi(t,vt)\rho^{\pm}, \rho^{
  \pm }  \right\rangle_{ H }}_{ = 0 }
    + O\left(|\rho^{\pm}|_{H}\,\epsilon t^{-1+C\epsilon}\Lambda^{-1/2}\right) \\
  &\lesssim \lvert \rho^{ \pm }\rvert_{ H }\cdot\epsilon t^{-\frac{3}{2}
  +C\epsilon}(1-v^2)^{-\frac{1}{4} }.
\end{align}
After dividing both sides by $ \lvert \rho^{ \pm } \rvert _{ H } $ and
rewriting $ \partial_{ t } \lvert \rho^{ \pm } \rvert = \frac{1}{2} \lvert \rho^{
\pm } \rvert ^{ -1 } \partial_{ t } \lvert \rho^{ \pm } \rvert^2   $, we get
\begin{align}
  \left\lvert\partial_t\lvert\rho^{\pm}\rvert_H\right\rvert
  &\lesssim \epsilon t^{-1+C\epsilon}\Lambda^{-\frac12}
   =\epsilon t^{-\frac32+C\epsilon}(1-v^2)^{-\frac14}.
\end{align}
Since the estimate above is globally integrable in $ t $, integrating from a fixed
$t_0$ to $t$ gives
\begin{align}\label{eq:improvedrho}
  \lvert \rho^{ \pm } (t) \rvert _{ H } 
  &\lesssim \lvert \rho^{ \pm } (t_0)  \rvert_{ H }
    +\epsilon t_0^{-\frac12+C\epsilon}(1-v^2)^{-\frac14}.
\end{align} 
In particular, \eqref{eq:improvedrho} is a significant improvement over
\eqref{eq:rho-apriori-growth}, provided that we pick $ t_0 $ appropriately.
Let
\begin{align}
  t_0(v):= (1-v^2)^{ -10 }.
\end{align} 
Figure \ref{fig:Cin} illustrates the initial time $t_0(v)$ along the ray $x=vt$,
together with the wave packet support on $H_\Lambda$.

Then, using \eqref{eq:Hnormrho} for $ t = t_0(v) $,
\begin{align}\label{eq:finalrhobound}
  \lvert \rho^{ \pm } (t) \rvert _{ H } 
  &\lesssim \epsilon t^{ C\epsilon } (1-v^2) + \epsilon t^{ -\frac{1}{2} 
  + C\epsilon } (1-v^2)^{ -\frac{1}{4}  } \Big|_{ t=t_0(v) }  
  \lesssim \epsilon (1-v^2)^{ 1-C\epsilon },
  \quad t \geq t_0(v).
\end{align} 

\begin{figure}[h]
  \centering
  \begingroup
\normalfont\fontencoding{OT1}\fontsize{10.95}{13.6}\selectfont
\definecolor{RegionNavy}{HTML}{182B49}
\definecolor{RegionBlue}{HTML}{00629B}
\definecolor{RegionGold}{HTML}{C69214}
\definecolor{RegionSand}{HTML}{F5F0E6}
\colorlet{RegionWavePacketDefault}{RegionGold!35!white}
\colorlet{RegionNearConeDefault}{RegionSand}
\colorlet{RegionTransferDefault}{RegionNavy!30!white}
\colorlet{RegionExteriorDefault}{RegionBlue!32!white}
\colorlet{RegionAxis}{RegionNavy}
\colorlet{RegionGuide}{RegionNavy!60!white}
\providecommand{\WavePacketRegionColor}{RegionWavePacketDefault}
\providecommand{\NearConeRegionColor}{RegionNearConeDefault}
\providecommand{\TransferRegionColor}{RegionTransferDefault}
\providecommand{\ExteriorRegionColor}{RegionExteriorDefault}
\colorlet{RegionWavePacketInk}{RegionNavy}
\colorlet{RegionTransferInk}{RegionNavy}
\color{RegionNavy}

\begin{tikzpicture}[
  x=1.08cm,y=1.08cm,
  line cap=round,line join=round,
  axis/.style={draw=RegionAxis,line width=.55pt,-{Stealth[length=4pt,width=3.2pt]}},
  boundary/.style={draw=RegionAxis,line width=.75pt},
  guide/.style={draw=RegionGuide,line width=.5pt,dash pattern=on 3pt off 3pt},
  every node/.style={inner sep=2pt},
  declare function={
    deep(\s)=.91*\s;
    outer(\s)=\s+.6*(pow(\s+.3,1/3)-pow(.3,1/3));
    hyper(\r)=sqrt(\r*\r+2.5*2.5);
  }
]
  \def\top{3.5}
  \def\bottom{0}
  \pgfmathsetmacro{\cupend}{sqrt(\top*\top-2.5*2.5)}
  \pgfmathsetmacro{\deepend}{deep(\top)}
  \pgfmathsetmacro{\ch}{(exp(.105)+exp(-.105))/2}
  \pgfmathsetmacro{\sh}{(exp(.105)-exp(-.105))/2}
  \pgfmathsetmacro{\collarcupx}{\cupend*\ch-\top*\sh}
  \pgfmathsetmacro{\collarcupt}{\top*\ch-\cupend*\sh}
  \pgfmathsetmacro{\collarendtime}{\top/(\ch+.91*\sh)}
  \pgfmathsetmacro{\collarendx}{deep(\collarendtime)}

  \foreach \side in {-1,1}{
    \fill[fill=\ExteriorRegionColor]
      (0,0) -- ({\side*5.05},0) -- ({\side*5.05},\top)
      -- ({\side*outer(\top)},\top)
      -- plot[variable=\s,domain=\top:\bottom,samples=111]
         ({\side*outer(\s)},\s)
      -- cycle;
  }

  \foreach \side in {-1,1}{
    \fill[fill=\NearConeRegionColor]
      (0,0)
      -- plot[variable=\s,domain=\bottom:\top,samples=111]
         ({\side*outer(\s)},\s)
      -- ({\side*deep(\top)},\top)
      -- plot[variable=\s,domain=\top:\bottom,samples=111]
         ({\side*deep(\s)},\s)
      -- cycle;
  }

  \fill[fill=\WavePacketRegionColor]
    (0,\bottom)
    -- plot[variable=\s,domain=\bottom:\top,samples=111] ({deep(\s)},\s)
    -- (\cupend,\top)
    -- plot[variable=\r,domain=\cupend:-\cupend,samples=141] (\r,{hyper(\r)})
    -- (-\deepend,\top)
    -- plot[variable=\s,domain=\top:\bottom,samples=111] ({-deep(\s)},\s)
    -- cycle;

  \foreach \side in {-1,1}{
    \path[fill=\TransferRegionColor]
      ({\side*\cupend},\top) -- ({\side*\deepend},\top)
      -- plot[variable=\s,domain=\top:\collarendtime,samples=31]
         ({\side*deep(\s)},\s)
      -- ({\side*\collarcupx},\collarcupt)
      -- plot[variable=\r,domain=\collarcupx:\cupend,samples=31]
         ({\side*\r},{hyper(\r)}) -- cycle;
    \path[pattern=north east lines,pattern color=RegionTransferInk!65]
      ({\side*\cupend},\top) -- ({\side*\deepend},\top)
      -- plot[variable=\s,domain=\top:\collarendtime,samples=31]
         ({\side*deep(\s)},\s)
      -- ({\side*\collarcupx},\collarcupt)
      -- plot[variable=\r,domain=\collarcupx:\cupend,samples=31]
         ({\side*\r},{hyper(\r)}) -- cycle;
    \draw[RegionTransferInk,line width=.55pt]
      ({\side*\collarcupx},\collarcupt) -- ({\side*\collarendx},\collarendtime);
  }

  \draw[guide] (-4.1,4.1) -- (0,0) -- (4.1,4.1);
  \draw[guide] (-\cupend,\top) -- (\cupend,\top);

  \draw[boundary] (-5.05,\top) -- (-\cupend,\top);
  \draw[boundary] (\cupend,\top) -- (5.05,\top);
  \draw[boundary] plot[variable=\r,domain=-\cupend:\cupend,samples=141]
    (\r,{hyper(\r)});

  \draw[axis] (-5.2,0) -- (5.35,0) node[right] {$x$};
  \draw[axis] (0,-.14) -- (0,4.4) node[above] {$t$};
  \node[below left] at (0,0) {$0$};
  \node[anchor=east,xshift=-3pt,yshift=6pt] at (0,\top) {$4T$};
  \draw[RegionAxis,line width=.5pt] (-.06,\top) -- (.06,\top);

  \node[text=RegionWavePacketInk] at (-.76,1.74) {$C^{\mathrm{wp}}_{<T}$};
  \node[anchor=west] at (4.14,4.05) {$t=|x|$};

  \fill[fill=\ExteriorRegionColor] (-4.90,-.68) rectangle (-4.62,-.48);
  \node[anchor=west] at (-4.53,-.58) {exterior};
  \fill[fill=\NearConeRegionColor] (-2.50,-.68) rectangle (-2.22,-.48);
  \node[anchor=west] at (-2.13,-.58) {near cone};
  \fill[fill=\WavePacketRegionColor] (.25,-.68) rectangle (.53,-.48);
  \node[anchor=west] at (.62,-.58) {$C^{\mathrm{wp}}_{<T}$};
  \fill[fill=\TransferRegionColor] (2.35,-.68) rectangle (2.63,-.48);
  \fill[pattern=north east lines,pattern color=RegionTransferInk!65]
    (2.35,-.68) rectangle (2.63,-.48);
  \node[anchor=west] at (2.72,-.58)
    {$D_{<T}\setminus C^{\mathrm{wp}}_{<T}$};
\end{tikzpicture}
\endgroup

  \caption{To close the bootstrap estimate we use the thin support of the
    near-cone region to compensate for weaker pointwise decay, use the wave packet
    analysis in $\Cin_{<T}$, and transfer the improved bounds to $D_{ <T }
    \setminus \Cin_{ <T }$ by FTC.
    In the exterior region, the improved pointwise bounds \eqref{eq:dirac-global-pointwise} suffice.}
  \label{fig:regions}
\end{figure}
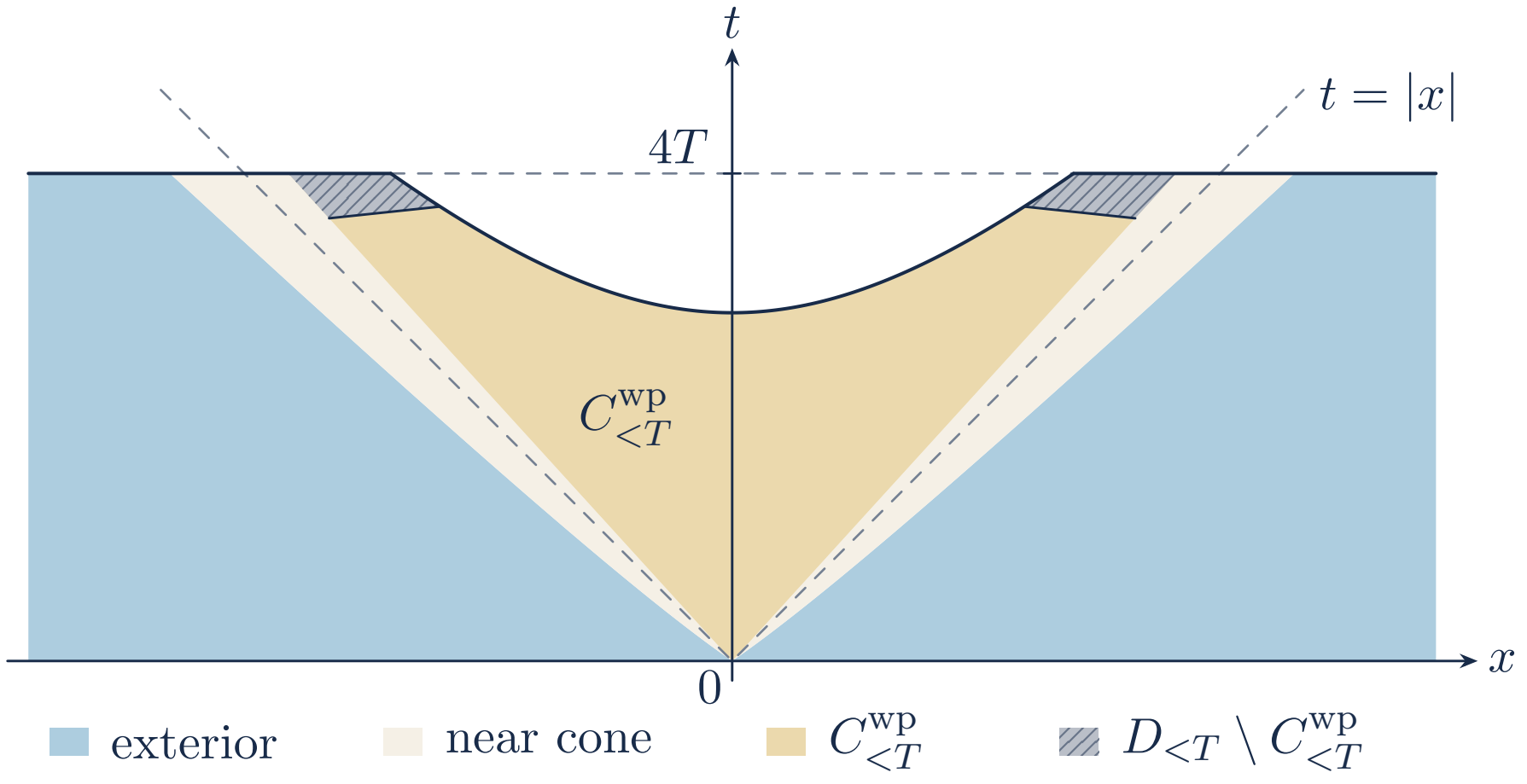

Combining \eqref{eq:Hdecomp}, \eqref{eq:finalrhobound} and Lemma \ref{asympdecomp},
we conclude that
\begin{align}
  \lvert P_+^\theta\psi  \rvert
  &\lesssim \left( \frac{ t+r }{ t-r }  \right) ^{ \frac{1}{4}  } \lvert \psi \rvert _{ H
    }\\
  &\lesssim \left( \frac{(t+r)^2 }{ \Lambda^{ 2 }  }  \right) ^{ \frac{1}{4}  }\Lambda^{
    -\frac{3}{2}  }  \lvert \rho \rvert _{ H } \\
  &\lesssim (t+r)^{ \frac{1}{2}  } \left( t\sqrt{ 1-v^2 }  \right) ^{ -2 } \cdot \epsilon 
  (1-v^2)^{ 1-C\epsilon } \\
  &\lesssim \epsilon t^{ -\frac{3}{2}  } (1-v^2)^{ -C\epsilon },
\end{align} 
and 
\begin{align}
  \lvert P_-^\theta\psi  \rvert
  &\lesssim \epsilon t^{ -\frac{3}{2}  } (1 -v^2)^{ \frac{1}{2} -C\epsilon }.
\end{align} 
Finally, rewriting $ (1-v^2) = (t-r)(t+r)t^{ -2 } $ along a ray, we get
\begin{align}\label{eq:fulldecaypsi}
  \lvert \psi \rvert 
  &\lesssim \epsilon t^{ -\frac{3}{2} + C\epsilon }(t-r)^{ -C\epsilon  }
  \quad\text{and}\quad 
  \lvert \psi \rvert _{ H } \lesssim \epsilon t^{ -\frac{3}{2}  } (1-v^2)^{ \frac{1}{4}
    -C\epsilon }.
\end{align} 
The pointwise decay rates in Proposition \ref{prop:close-bootstrap-psi}
would then follow by combining \eqref{eq:fulldecaypsi} inside the cone and the
Klainerman-Sobolev bound \eqref{eq:dirac-global-pointwise} outside the cone.
Before we can conclude, however, we must address a technical issue in the
asymptotic analysis above.
The definition of the profile $ \rho^{ \pm }(t,v) $ depends on testing $ \psi $ against the
wave packet $ \wp_{ v } ^{ \pm }  $ supported on a small hyperbolic ball
\begin{align}
  B_{ H_{ \Lambda }  } (v, c  \Lambda^{ -\frac{1}{2}  }) \subset  H_{ \Lambda },
  \quad \Lambda = t\sqrt{ 1-v^2 },
\end{align}
where $ c > 0 $ is some small constant.
Consequently, if the bootstrap assumption is only known in $C_{<T}$, the profile and
the asymptotic equation are directly available only at those points for which the above ball
is still contained in $C_{<T}$.
This gives the following overlapping cover of
$C_{<T}\cap\{t-r\geq1\}$, illustrated in Figure \ref{fig:regions}.
\begin{itemize}
  \item The \textit{deep-interior region} 
    \begin{align}
      D_{<T}
      &:=C_{<T}\cap\{t-r\geq1\}\cap\{t\geq t_0(v)\},
    \end{align}
    where the loss in \eqref{eq:rho-apriori-growth} prevents us from closing the
    bootstrap argument directly. 

  \item 
    The \textit{wave packet region}
      \begin{align}
        \Cin_{<T}
        :=\left\{(t,x)\in D_{ <T } : B_{H_{ \Lambda } }\left(x/t,2c
          \Lambda^{-\frac12}\right) \subset H_{\Lambda}\cap C_{<T}\right\},
      \end{align}
      where the packet support remains inside $C_{<T}$ with a margin,
      and one can use wave packet testing to improve the a priori profile bound
      \eqref{eq:rho-apriori-growth} to \eqref{eq:finalrhobound},

    \item The \textit{near-boundary region}
      \begin{align}
        \left\{ t-\lvert x \rvert  \geq 1 \right\} \setminus D_{ <T }
      \end{align}
      where \eqref{eq:dirac-global-pointwise} provides enough pointwise control
      independently of the wave packet analysis.
\end{itemize}
Indeed, on the thin near-boundary region above, $t\leq t_0(v)$ and hence
$\langle t-r\rangle\lesssim\langle t+r\rangle^{9/10}$.  Thus,
\eqref{eq:dirac-global-pointwise} 
gives
\begin{align}
  \lvert\psi(t,x)\rvert
  &\lesssim\epsilon\langle t+r\rangle^{-\frac32+C\epsilon}, 
\end{align}
thus establishing the first estimate in Proposition \ref{prop:close-bootstrap-psi}
in this region.

In particular, it only remains to establish the first estimate in Proposition
\ref{prop:close-bootstrap-psi} on
$ D_{ <T } \setminus \Cin_{ <T }  $.
Given $ (t,x) \in D_{ <T }  $, the only obstruction to $ (t,x) \in \Cin_{ < T} $ is that
$(t,x)$ lies along a hyperboloid $ H_{ \Lambda }  $ and is within hyperbolic distance
$O_{ H }(\Lambda^{-1/2})$ of the top boundary $t=4T$ of $ C_{ <T }  $.

Let $(t_1,x_1)\in D_{<T}\setminus \Cin_{<T}$ and set
\[
  v_1=\frac{x_1}{t_1},
  \qquad
  \Lambda =(t_1^2-x_1^2)^{\frac12}.
\]
Hence, there exists a point
$(t_2,x_2)\in \Cin_{<T}$ on the same hyperboloid $H_{\Lambda}$ such that
\begin{align}%
  d_{H}\left(\frac{x_1}{t_1},\frac{x_2}{t_2}\right)
  \lesssim \Lambda^{-\frac12}.
\end{align}

We compare the values of $\psi$ at $ (t_1,x_1) $ and $ (t_2,x_2) $ using the
Fundamental Theorem of Calculus along $ H_{ \Lambda } $. 
We use the Lorentz generators $\widehat\Omega$ as directional derivatives, since they
preserve the $H$-inner product.  
First, note that by a Sobolev embedding,
\eqref{eq:EnergyToL2Norm}, and Corollary \ref{coro:truncated-energy}, we have the
a priori vector field bound
\begin{align}\label{eq:OmegaHpointwise}
  \lvert \widehat\Omega\psi(t,x)\rvert_{H}
  &\lesssim \Lambda^{-\frac32}
    E_{H_{\Lambda,<T}}(\widehat\Omega^{\leq 3}\psi)^{\frac12}
   \lesssim \epsilon \Lambda^{-\frac32}T^{C\epsilon},
   \qquad (t,x)\in H_{\Lambda}\cap C_{<T}.
\end{align}

Therefore, if $\gamma$ is a unit-speed  geodesic in $ \H $ joining $x_2/t_2$ and 
$x_1/t_1$ with length $ \ell(\gamma) = O(\Lambda^{ -\frac{1}{2}  } ) $, then
\begin{align}
  \lvert \psi(t_1,x_1)\rvert_{H}
  &\leq \lvert \psi(t_2,x_2)\rvert_{H}
    + \int_{0}^{\ell(\gamma)}\lvert \widehat\Omega_{H}\psi(\gamma(s))\rvert_{H}\dif s \\
  &\lesssim \epsilon t_2^{-\frac32}(1-v_2^2)^{\frac14-C\epsilon}
    + \Lambda^{-\frac12}\epsilon \Lambda^{-\frac32}T^{C\epsilon} \\
  &\lesssim \epsilon t_1^{-\frac32}(1-v_1^2)^{\frac14-C\epsilon}
    + \epsilon \Lambda^{-2}T^{C\epsilon}.
\end{align}

In passing to the third line above, we have used that $ t_1, t_2\sim T $ and that $ (t_1,x_1) $
and $ (t_2,x_2) $ belong to the same hyperboloid $ H_{ \Lambda }  $, and so
\begin{align}
  t_1^{ -\frac{3}{2}  } (1-v_1^2)^{\frac{1}{4} }
  = \Lambda^{ \frac{1}{2}  } t_{ 1 } ^{ -2 } 
  \sim t_2^{ -\frac{3}{2}  } (1-v_2^2)^{\frac{1}{4} }.
\end{align} 
Moreover, in $ D_{ <T }  $, $ t_1\gtrsim (1-v^2)^{ -10 }  $, and so we directly obtain
\begin{align}
  \epsilon \Lambda^{ -2 } t_1^{ C\epsilon } 
  \lesssim \epsilon t_{ 1 }^{ -\frac{3}{2}  } (1-v_1^2)^{ \frac{1}{4} -C\epsilon },
\end{align} 
showing that $ \psi $ satisfies the pointwise bounds \eqref{eq:fulldecaypsi} in $ D_{ <T
}  $.
This completes the proof of Proposition \ref{prop:close-bootstrap-psi}.

\subsection{Closing the bootstrap in \texorpdfstring{$ \psi $}{psi} and the
proof of Corollary \ref{coro:close-bootstrap-psi}}\label{sec:closing-psi}

We can now close the $ L^{ 6 }  $ part of the bootstrap assumption
\eqref{eq:bootstrap-assumption} and prove Corollary
\ref{coro:close-bootstrap-psi}.
Fix $ T > 0 $ and a constant-time slice $ C_{<T}(t) \subset  C_{ <T }  $.
We split the $ L^6_x $ integral into the three overlapping regions given by Figure
\ref{fig:regions}.

First, in the region $ D_{ <T }   $, we use the improved pointwise bounds
\eqref{eq:fulldecaypsi} to estimate
\begin{align}
  \int_{D_{ <T }(t)}\lvert \psi(t,x)\rvert^{6}\dif x
  &\lesssim \epsilon^{6}t^{-9}
    \int_{D_{<T}(t)}(1-v^2)^{-6C\epsilon}\dif x \\
  &\lesssim \epsilon^{6}t^{-6}
    \int_{B(0,1)}(1-|v|^2)^{-6C\epsilon}\dif v \\
  &\lesssim \epsilon^{6}t^{-6}.
\end{align}
Here, the last integral is finite provided $6C\epsilon<1$.
This estimate is also the bottleneck preventing one from closing
\eqref{eq:bootstrap-assumption}
using the Klainerman-Sobolev analysis.

On the complement of $ D_{ <T }  $ inside the cone, we have
$t\lesssim t_0(v)=(1-v^2)^{-10}$ and hence $t-r\lesssim t^{9/10}$.  
In the thin neighborhood of the light cone $ \left\{ -t^{ \frac{9}{10}  } \lesssim r-t
\lesssim t^{ \frac{1}{3}  }  \right\}  $, it suffices to consider the Klainerman-Sobolev
bound \eqref{eq:dirac-global-pointwise}, together with \eqref{eq:energy-bound}:
\begin{align}
  \int_{\left\{-t^{9/10}\lesssim r-t\lesssim t^{1/3}\right\}}
    \lvert \psi(t,x)\rvert^{6}\dif x
  &\lesssim \epsilon^{6}t^{-9+C\epsilon}
    \,t^2\bigl(t^{9/10}+t^{1/3}\bigr) 
   \lesssim \epsilon^{6}t^{-6},
\end{align}
provided $C\epsilon<1/10$. 
Finally, outside the cone we use the exterior part of
\eqref{eq:dirac-global-pointwise}:
\begin{align}
  \int_{\left\{r\geq t+t^{ 1 /3 } \right\}}\lvert \psi(t,x)\rvert^{6}\dif x
  &\lesssim \epsilon^{6}t^{C\epsilon}
    \int_{ \left\{ r\geq t + t^{ 1 /3 }  \right\} }
    \frac{r^{2}}{\langle t+r\rangle^{9}\langle r-t\rangle^{6\delta}}\dif r 
  \lesssim \epsilon^{6}t^{-6},
\end{align}
provided, again, that $ C\epsilon $ is taken small enough.
Combining the three estimates gives
\begin{align}
  \lVert \psi(t)\rVert_{L^{6}_{x}}
  \lesssim \epsilon t^{-1},
\end{align}
with an implicit constant independent of the bootstrap constant.  Taking the bootstrap
constant $C_{0}$ large and then $\epsilon$ sufficiently small closes the $\psi$ component of
\eqref{eq:bootstrap-assumption}.

\section{Bounds on \texorpdfstring{$ \phi $}{phi}}\label{sec:wavebounds}

In this section, we prove the following proposition, which, when combined with
Proposition \ref{prop:close-bootstrap-psi} and Corollary
\ref{coro:close-bootstrap-psi}, closes the bootstrap argument \eqref{eq:bootstrap-assumption}.
\begin{prop}\label{prop:close-bootstrap-phi}
  For $ C_0 $ sufficiently large, there exist constants $ C'\in(0,C_0/2) $ and
  $ \epsilon_0=\epsilon_0(C_0)>0 $, independent of $ T $, such that for every
  dyadic $ T\geq1 $ and $ 0<\epsilon\leq\epsilon_0 $, if the initial data satisfy
  the hypotheses of Theorem \ref{main-thm-gwp} and the solution satisfies the
  bootstrap assumption \eqref{eq:bootstrap-assumption}, then
  \begin{align}
    \lVert \phi(t)\rVert_{L_x^\infty(C_{<T}(t))}
    &\leq C'\epsilon\langle t\rangle^{-1},\qquad 1\leq t\leq4T.
  \end{align}
  In particular,  $ \psi(t) $ and  $ \phi(t) $ are well-defined for all $ t > 0
  $ and satisfy the dispersive bound
  \begin{align}
    \lVert \psi(t) \rVert_{ L_{ x }^{ 6 }  }  + \lVert \phi (t)\rVert_{ L_{ x }^{
    \infty}  } \lesssim \epsilon\langle t \rangle  ^{ -1 }.
  \end{align} 
\end{prop}

We consider the homogeneous and inhomogeneous contributions of $\phi$ separately
\[
  \phi = \wavehom + \waveinhom.
\] 

In Section \ref{sec:prelim}, we used the linear Dirac estimates from \cite{HIS},
stated in Proposition \ref{prop:klainerman-sobolev}, which apply to both the
Dirac--Maxwell and \eqref{dkg} systems. We likewise assume the following estimate
for linear waves; a complete proof of \eqref{eq:free-wave-pointwise} is given in
\cite[Proposition 7.1]{HIS}.

\begin{prop}[Free wave estimate]\label{prop:wave-hom-bound}
  Assume that the initial data $ (\phi_0,\phi_1) $ are as in Theorem
  \ref{main-thm-gwp}.
  Then, the solution $ \wavehom $ to 
  \begin{align}
    &\square\wavehom = 0,\\
    &(\wavehom,\partial_{ t } \wavehom)\big|_{ t=0 } = (\phi_0,\phi_1)
  \end{align} 
  satisfies
  \begin{align}\label{eq:free-wave-pointwise}
    \lvert \wavehom (t,x) \rvert 
    &\lesssim \frac{ \epsilon }{ \langle t + r \rangle   \langle t-r \rangle ^{
        \nu /6 } },
  \end{align} 
  uniformly in $ (t,x) \in \R^{ 1+3 } $.
\end{prop}

We now focus on the inhomogeneous contribution
\begin{align}\label{eq:wave-inhom}
  \square \waveinhom = F,
  \qquad
  F:=\langle \psi,\gamma^{0}\psi\rangle,
  \qquad
  \waveinhom[0]=0.
\end{align}
We will use the improved pointwise bounds for $\psi$ from
Proposition~\ref{prop:close-bootstrap-psi}. Inside the cone these give
\begin{align}
  |\psi(t,x)|
  \lesssim
  \epsilon \langle t+r\rangle^{-\frac32+C\epsilon} \langle t-r\rangle^{-c\epsilon},
\end{align}
and therefore, 
\begin{align}
  |F(t,x)|
  \lesssim
  \epsilon^{2} \langle t+r\rangle^{-3+c\epsilon} \langle t-r\rangle^{-c\epsilon}.
\end{align}
Outside the cone, the improved exterior decay gives
\begin{align}
  |\psi(t,x)|
  \lesssim
  \epsilon
  \langle t+r\rangle^{-\frac32+c\epsilon}
  \langle r-t\rangle^{-\delta},
\end{align}
and hence
\begin{align}
  |F(t,x)|
  \lesssim
  \epsilon^{2} \langle t+r\rangle^{-3+c\epsilon} \langle r-t\rangle^{-\delta}.
\end{align}
As we shall see, by the positivity of the wave kernel, it suffices to estimate 
a wave with positive radial source
\begin{nalign}\label{eq:sourcebounds}
  \Gin(t,r)
  &:= \epsilon^{2} \langle t+r\rangle^{-3+c\epsilon} \langle
    t-r\rangle^{-c\epsilon},\\
  \Gout(t,r)
  &:=
    \epsilon^{2} \langle t+r\rangle^{-3+c\epsilon} \langle r-t\rangle^{-\delta},
\end{nalign}
localized respectively to the interior and exterior of the cone.

\begin{prop}\label{prop:wave-inhom-bound}
The inhomogeneous wave $ \waveinhom $ satisfies
  \begin{align}\label{eq:inhom-bound}
    \lvert \waveinhom(t,x) \rvert
    \lesssim
    \epsilon^{2}
    \langle t\rangle^{-1}
    \langle t-r\rangle_-^{-\delta}.
  \end{align}
\end{prop}

\begin{proof}
  Let $\Phi$ solve
  \begin{align}
    \square \Phi
    = G
    := \Gin(t,r)\1_{\{t\geq r\}} + \Gout(t,r)\1_{\{r>t\}},
  \end{align}
  with zero initial data.
  Since $ \lvert F \rvert \lesssim G $, by Kirchhoff's formula
  \begin{align}
    |\waveinhom(t,x)|\lesssim \Phi(t,r).
  \end{align}
  Moreover, since the source $ G $ is radial, so is $\Phi$. 
  The radial representation formula gives
  \begin{align}
    \Phi(t,r)
    &=
      \frac{1}{2r} \int_{ D(t,r)\cap \left\{ s\geq \rho \right\}  } \rho\Gin(s,\rho)  \dif
      \rho\dif s
      + \frac{1}{2r} \int_{ D(t,r)\cap \left\{ s < \rho \right\} } \rho\Gout(s,\rho)\dif
      \rho\dif s\\
    &=: (\mathrm{I}) + (\mathrm{II})
  \end{align}
  where $ D(t,r) $ is the truncated null rectangle emanating from $ (t,r) $ and
  represented in Figure \ref{fig:nullrec}.
  It is sufficient to show that
  \begin{align}
    (\mathrm{I}) \lesssim \epsilon ^2 \langle t \rangle  ^{ -1 } 
    \quad\text{and}\quad 
    (\mathrm{II}) \lesssim \epsilon ^2 \langle t \rangle  ^{ -1 }
    \langle t-r \rangle_-^{ -\delta }.
  \end{align} 
   Without loss of generality, we assume below that $ t\geq 2 $ and $ r > 0 $.
  
   \begin{figure}[h]
     \centering
     \begingroup
\definecolor{UCSDNavy}{HTML}{182B49}
\definecolor{UCSDBlue}{HTML}{00629B}
\definecolor{UCSDGold}{HTML}{C69214}
\colorlet{NullConeAxis}{UCSDNavy}
\colorlet{NullConeBoundary}{UCSDNavy}
\providecommand{\NullConeDomainColor}{UCSDBlue!55!white}
\providecommand{\NullConeCancellationColor}{UCSDGold!55!white}

\begin{tikzpicture}[
  x=.8cm,y=.8cm,
  font=\fontsize{10.95}{13.2}\selectfont,
  line cap=round,line join=round,
  axis/.style={draw=NullConeAxis,line width=.65pt,-{Stealth[length=5pt,width=4pt]}},
  boundary/.style={draw=NullConeBoundary,line width=.85pt},
  guide/.style={draw=NullConeAxis!65,line width=.45pt,dash pattern=on 2.5pt off 2.5pt},
  tick/.style={draw=NullConeAxis,line width=.6pt},
  every node/.style={text=NullConeBoundary,inner sep=2pt}
]
  \pgfmathsetmacro{\rr}{1.55}
  \pgfmathsetmacro{\tt}{4.10}
  \pgfmathsetmacro{\difference}{\tt-\rr}
  \pgfmathsetmacro{\total}{\tt+\rr}

  \coordinate (apex) at (\rr,\tt);
  \coordinate (crossing) at (0,\difference);
  \coordinate (left) at (-\difference,0);
  \coordinate (reflected) at (\difference,0);
  \coordinate (right) at (\total,0);

  \fill[fill=\NullConeCancellationColor] (left) -- (crossing) -- (reflected) -- cycle;
  \fill[fill=\NullConeDomainColor] (crossing) -- (apex) -- (right) -- (reflected) -- cycle;

  \draw[guide] (0,\tt) -- (apex);
  \draw[guide] (\rr,0) -- (apex);
  \draw[boundary] (left) -- (apex) -- (right);
  \draw[boundary] (crossing) -- (reflected);

  \draw[axis] (-3.05,0) -- (6.25,0) node[right=2pt] {$\rho$};
  \draw[axis] (0,-.35) -- (0,4.75) node[above=2pt] {$s$};

  \foreach \position/\label in {
    -\difference/{r-t},\rr/r,\difference/{t-r},\total/{r+t}}
    {\draw[tick] (\position,.065) -- (\position,-.065)
       node[below=4pt] {$\label$};}
  \draw[tick] (-.065,\tt) -- (.065,\tt);
  \node[left=5pt] at (0,\tt) {$t$};
  \node[below left=3pt] at (0,0) {$0$};

  \fill[NullConeBoundary] (apex) circle[radius=1.5pt];
  \node[above right=4pt] at (apex) {$(r,t)$};
  \fill[NullConeBoundary] (crossing) circle[radius=1.5pt];
  \node[above left=5pt] at (crossing) {$(0,t-r)$};
  \node[fill=none,inner sep=2pt]
    at (2.65,1.50) {$D(t,r)$};
\end{tikzpicture}
\endgroup

     \caption{Representation of the null rectangle $ D(t,r) $ in blue.
     The triangle in gold accounts for the cancellation coming from the even extension in the
   radial representation formula.}
     \label{fig:nullrec}
   \end{figure}
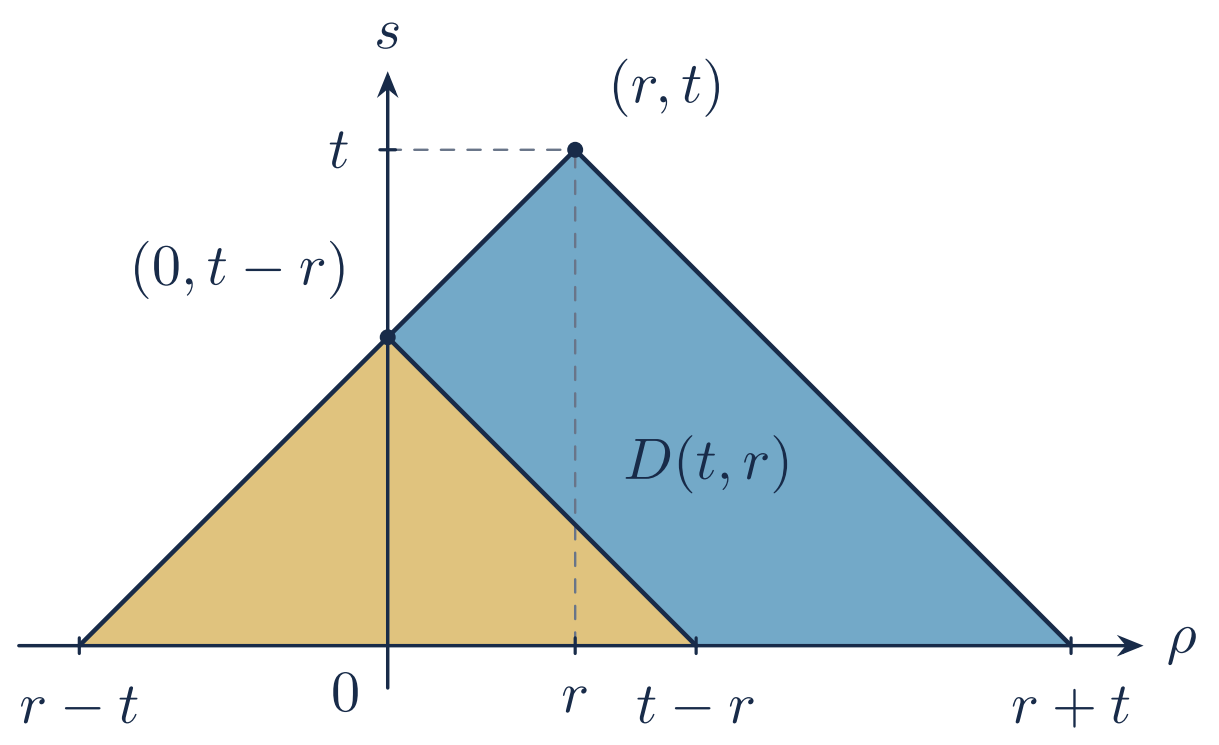

   We first estimate the contribution of $\Gout$ when $ r > t $.  Then $ D $
   corresponds to the full backward light cone
   emanating from $ (r,t) $.
   It is natural to use the null coordinates
  \begin{align}
    u=\rho-s,
    \qquad
    v=\rho+s,
  \end{align}
  where
  \begin{align}
    \rho=\frac{u+v}{2},
    \qquad
    \dif \rho\dif s \sim \dif u\dif v,
  \end{align}
  and $ D $ is parametrized by
  \begin{align}
    r-t\leq u\leq r+t,
    \qquad
    u\leq v\leq r+t.
  \end{align}
  Therefore, since $ (u+v)\leq 2v $ in $ D $,
  \begin{align}
    (\mathrm{II})
    &\lesssim
      \frac{\epsilon^{2}}{r} \int_{r-t}^{r+t} \int_{u}^{r+t} (u+v)\langle v\rangle^{-3+c\epsilon}
      \langle u\rangle^{-\delta} \dif v\dif u\\
    &\lesssim
      \frac{\epsilon^{2}}{r} \int_{r-t}^{r+t} \langle u\rangle^{-\delta} \int_{u}^{r+t}
      \langle v\rangle^{-2+c\epsilon} \dif v\dif u\\
    &\lesssim
      \frac{\epsilon^{2}}{r} \int_{r-t}^{r+t} \langle u\rangle^{-1-\delta} \dif u\\
    &\lesssim
    \epsilon^{2} t^{-1} \langle r-t\rangle^{-\delta}.
  \end{align}
  If $r<t$ and the backward cone meets the exterior-source region, the same
  null-coordinate estimate on $D(t,r)\cap\{\rho>s\}$ gives
  $(\mathrm{II})\lesssim\epsilon^2t^{-1}$.  

  It remains to estimate the contribution of $\Gin$, which only appears when $r<t$. 
  We use slightly different null coordinates
  \begin{align}
    u=s-\rho,
    \qquad
    v=s+\rho,
  \end{align}
  and parametrize $ D $ by
  \begin{align}
    0\leq u\leq t-r,
    \qquad
    t-r\leq v\leq t+r.
  \end{align}
  Since $ 2\rho = (v-u)\leq v $, we get
  \begin{align}
    (\mathrm{I}) 
    &\lesssim
    \frac{\epsilon^{2}}{r}
    \int_{0}^{t-r}
    \int_{t-r}^{t+r}
      (v-u)\langle v\rangle^{-3+c\epsilon}
      \langle u\rangle^{-c\epsilon}
    \dif v\dif u\\
    &\lesssim
    \frac{\epsilon^{2}}{r}
    \left(
      \int_{0}^{t-r}
        \langle u\rangle^{-c\epsilon}
      \dif u
    \right)
    \left(
      \int_{t-r}^{t+r}
        \langle v\rangle^{-2+c\epsilon}
      \dif v
    \right).
  \end{align}
  We now consider two cases depending on whether the $ r^{ -1 } $ contribution
  improves or worsens the decay in $ t $.
  If $t/2\leq r<t$, then
  \begin{align}
    (\mathrm I)
    \lesssim
    \frac{\epsilon^{2}}{t}
    \langle t-r\rangle^{1-c\epsilon}
    \langle t-r\rangle^{-1+c\epsilon}
    \lesssim
    \epsilon^{2}t^{-1}.
  \end{align}
  If instead $0<r\leq t/2$, then
  \begin{align}
    \int_{t-r}^{t+r}
      \langle v\rangle^{-2+c\epsilon}
    \,\dif v
    \lesssim
    r\langle t-r\rangle^{-2+c\epsilon},
  \end{align}
  and consequently
  \begin{align}
    (\mathrm I)
    \lesssim
      \frac{\epsilon^{2}}{r}\langle t-r\rangle^{1-c\epsilon}
      \cdot r\langle t-r\rangle^{-2+c\epsilon}
    \lesssim
      \epsilon^{2}t^{-1},
  \end{align}
  completing the proof.
\end{proof}

Combining Propositions \ref{prop:wave-hom-bound} and \ref{prop:wave-inhom-bound},
we obtain the global wave bound
\begin{align}
  \lvert\phi(t,x)\rvert
  &\lesssim
    \frac{\epsilon}{\langle t \rangle  },
\end{align}
inside the light cone, which closes the bootstrap assumption for the wave component.

\section{Modified scattering}\label{sec:scattering}
In this section, we use the asymptotic equation for the profiles
$\rho^{\pm}$, together with the wave bounds from section~\ref{sec:wavebounds},
to identify the leading asymptotic behavior of $\psi$ and $\phi$ inside
the light cone.  
The construction proceeds in three steps, each corresponding to one of the
subsections below.  
\begin{itemize}
  \item First, we obtain a preliminary asymptotic profile for $\psi$, with an
    arbitrary phase normalization. 

  \item Second, we use this profile for $ \psi $ to identify the radiation profile of $\phi$.

  \item Third, we use the asymptotic profile for $ \phi $ to identify the
    appropriate phase normalization for the Dirac profiles.
\end{itemize}

\begin{convention}
  In this section we work inside the wave packet region
  \begin{align}
    \Din
    := \left\{ t-r\geq 1 \right\} \cap \left\{ t\geq t_0(v) \right\},
    \qquad t_0(v) = (1-v^2)^{ -10 },
  \end{align} 
  where the wave packet estimates provide uniform control of the
  Dirac profiles. 
  Every fixed timelike ray eventually enters this region,
  while the wave generated by the source outside it contributes only a
  lower-order term to the timelike asymptotics. 
  Indeed, repeating the argument in Proposition \ref{prop:wave-inhom-bound} for the
  equation $ \square \phi^{ \text{ext}} = G^{ \text{ext} } $ with source
  $ G^{ \text{ext} }= 
  \1_{(\Din)^{ c }  }G $ gives the
  improved bound
  \begin{align}
    \lvert \phi^{ \text{ext} } (t,x)\rvert \lesssim \epsilon t^{ -1 }\langle
    t-r \rangle  ^{ -\delta }
  \end{align} 
  with a uniform constant, see \cite[Remark 7.4]{HIS}. 
  The source $ G $ is defined in \eqref{eq:wave-inhom}.
  
  All implicit constants in the
  estimates below are uniform for $ (t,vt)\in \Din $ , with
  velocity weights displayed explicitly.
\end{convention}

\subsection{The first Dirac profile}\label{sec:first-dirac}

For the convenience of the reader, we recall the asymptotic equation obtained
in Lemma~\ref{lemma:asymptoticeq}. 
For each $v\in B(0,1)$ and $t\geq t_{0}(v)$, we have
\begin{align}\label{eq:asymprecall}
  i\partial_{t}\rho^{\pm}(t,v)
  =
    \pm \sqrt{1-v^{2}}\phi(t,vt)\rho^{\pm}(t,v)
    + O_{H}\left( \epsilon t^{-\frac32+c\epsilon}(1-v^{2})^{-\frac14} \right).
\end{align}
The key point in \eqref{eq:asymprecall} is that the error term is integrable in time. 
Thus, the long-time behavior of $\rho^{\pm}$ is governed by the
effective profile equation
\begin{align}\label{eq:principal-rho-ode}
  i\partial_{t}\beta^{\pm}(t,v)
  =
  \pm \sqrt{1-v^{2}}\phi(t,vt)\beta^{\pm}(t,v).
\end{align}
Moreover, since the coefficient $\sqrt{1-v^{2}}\phi(t,vt)$ is real-valued,
\eqref{eq:principal-rho-ode} only affects the phase of $\beta^{\pm}$ and
does not change its norm.  

Therefore, solutions to \eqref{eq:principal-rho-ode} have the form
\begin{align}
  \beta^{\pm}(t,v)
  = e^{\mp i\theta(t,v)}A^{\pm}(v),
\end{align}
where $ A^{\pm} $ are time-independent profiles and $ \theta $ is explicitly given by
\begin{align}\label{eq:first-dirac-phase}
  \theta(t,v)
  :=
  \int_{t_{0}(v)}^{t}
    \sqrt{1-v^{2}}\phi(s,vs)\dif s.
\end{align}
Accordingly, we introduce the renormalized profiles
\begin{align}\label{eq:renormalized-rho}
  \widetilde\rho^{\pm}(t,v)
  :=
  e^{\pm i\theta(t,v)}\rho^{\pm}(t,v),
\end{align}
which solve the purely perturbative equation
\begin{align}\label{eq:perturbative-equation}
  i\partial_{t}\widetilde\rho^{\pm}(t,v)
  = O_{H}\left( \epsilon t^{-\frac32+c\epsilon}(1-v^{2})^{-\frac14} \right).
\end{align}

\begin{prop}\label{prop:first-dirac-profile}
  For each $v\in B(0,1)$, the limits
  \begin{align}\label{eq:rho-infty-def}
    \rho_{\infty}^{\pm}(v)
    :=
    \lim_{t\to\infty}
    \widetilde\rho^{\pm}(t,v)
  \end{align}
  exist.  
  Moreover, for $t\geq t_{0}(v)$,
  \begin{nalign}\label{eq:rho-infty-bounds}
    |\rho_{\infty}^{\pm}(v)|_{H}
    &\lesssim
      \epsilon(1-v^{2})^{1-C\epsilon},\\
    |\widetilde\rho^{\pm}(t,v)
      -\rho_{\infty}^{\pm}(v)|_{H}
    &\lesssim
      \epsilon t^{-\frac12+c\epsilon}(1-v^{2})^{-\frac14}.
  \end{nalign}
  Finally, the total charge is preserved by the asymptotic profiles
  \begin{align}\label{eq:profile-energy-conservation}
    \sum_{\pm}
    \lVert \rho_{\infty}^{\pm}  \rVert_{L^{2}(\H)}^{2}
    = \lVert \psi_{0} \rVert_{L^{2}_{x}}^{2}.
  \end{align}
\end{prop}

\begin{proof}
  For $t_{2}\geq t_{1}\geq t_{0}(v)$, we integrate \eqref{eq:perturbative-equation} to
  obtain
  \begin{align}
    |\widetilde\rho^{\pm}(t_{2},v)
      -\widetilde\rho^{\pm}(t_{1},v)|_{H}
    &\lesssim
      \epsilon \int_{t_{1}}^{t_{2}} s^{-\frac32+c\epsilon}(1-v^{2})^{-\frac14} \dif s\\
    &\lesssim
      \epsilon t_{1}^{-\frac12+c\epsilon}(1-v^{2})^{-\frac14}.
  \end{align}
  This proves the existence of the limits in \eqref{eq:rho-infty-def} and, when combined with
  the a priori bound \eqref{eq:finalrhobound}, furnishes \eqref{eq:rho-infty-bounds}.

  We now prove the charge conservation \eqref{eq:profile-energy-conservation}.
  The idea here is to formulate the charge conservation for the Dirac component $
  \psi $ in the hyperbolic metric, use the asymptotic decomposition
  \eqref{eq:dirac-profile-decomposition} to transfer information from $ \psi $ to the profiles $
  \rho^{ \pm }(t,v)  $, and take the limit as $ t\to\infty $ to write the
  conservation law in terms of the asymptotic profiles
  $ \rho_{\infty}^{\pm}(v) $.

  Fix $ T,\Lambda > 0 $ and let
  \begin{align}
    \Dcal_{ \Lambda,T } := \left\{ 0\leq t\leq T, t^2-r^2\leq \Lambda^2 \right\}
  \end{align} 
  denote the region delimited by $ t=0 $, $ t=T $ and the hyperboloid $ H_{
  \Lambda }  $.  
  Since $\phi$ is real, the current in \eqref{eq:dirac-current} is
  divergence free.  
  Integrating this conservation law on $ \Dcal_{ \Lambda,T }
  $ and applying the divergence theorem gives
  \begin{align}\label{eq:truncated-charge}
    \int_{\mathbb R^3} |\psi_0|^2\dif x
    &=
    \Lambda^3 \int_{H_\Lambda\cap\{t\leq T\}} |\psi|_H^2\dif V_H
    + \int_{\{t=T,\,T^2-r^2\leq \Lambda^2\}} |\psi|^2\dif x .
  \end{align}
  Equivalently, for $T\geq \Lambda$, the last integral is over the exterior
  part of the time slice
  \begin{align}
    \left\{ t=T,\, r\geq r_\Lambda(T) \right\},
    \qquad
    r_\Lambda(T):=(T^2-\Lambda^2)^{\frac12}.
  \end{align}

  We first show that the constant-time-slice contribution satisfies good bounds in
  $ T $ and $ \Lambda $.
  We split
  \begin{align}
    \{r\geq r_\Lambda(T)\}
    =
    \{r_\Lambda(T)\leq r\leq T\} \cup \{r\geq T\}
  \end{align}
  into contributions near the cone and in the exterior of the cone.
  For the near-cone shell, we use the pointwise bounds for $\psi$ inside the cone.
  Since
  \begin{align}
    T-r_\Lambda(T)
    =
    \frac{\Lambda^2}{T+r_\Lambda(T)}
    \lesssim
    \frac{\Lambda^2}{T},
  \end{align}
  we obtain
  \begin{align}
    \int_{\{r_\Lambda(T)\leq r\leq T\}} |\psi(T,x)|^2\dif x
    &\lesssim \epsilon^2 T^{-3+2C\epsilon}\vol\left( \left\{ r_{ \Lambda } (T)\leq
      r\leq T \right\} \right)\\
    &\lesssim \epsilon^2 T^{-3+2C\epsilon}T^{ 2 } \left(T-r_\Lambda(T)\right)\\
    &\lesssim \epsilon^2\Lambda^2 T^{-2+2C\epsilon}. 
  \end{align}
  For the exterior region, we use the exterior pointwise bound for $\psi$:
  \begin{align}
    \int_{\{r\geq T\}} |\psi(T,x)|^2\dif x
    \lesssim \epsilon^2 \int_T^\infty \langle T+r\rangle^{-3} \langle
      r-T\rangle^{-2\delta} r^2\dif r
    \lesssim \epsilon^2 T^{-2\delta}.
  \end{align}

  \noindent
  Now we go back to the hyperboloidal component of \eqref{eq:truncated-charge}.
  Write
  \begin{align}
    \psi(t,vt) = \sum_{\pm} \Lambda^{-\frac32}e^{\pm i\Lambda}\rho^\pm(t,v) +
    \err_\Lambda(t,v),
  \end{align}
  where the error term $ \err_{ \Lambda }  $ satisfies the pointwise bounds 
  \begin{align}\label{eq:truncated-dirac-profile-error}
    \lvert \err_{ \Lambda }  \rvert  _{ H } 
    &\lesssim \epsilon \Lambda^{ -\frac{7}{4}  } T^{ C\epsilon } 
  \end{align} 
  from Lemma \ref{asympdecomp}.
  Expanding the hyperboloidal flux gives
  \begin{align}
    \Lambda^3 \int_{H_\Lambda\cap\{t\leq T\}} |\psi|_H^2 \dif V_H
    &=
    \Lambda^3 \int_{H_\Lambda\cap\{t\leq T\}}
      \left\lvert \sum_{ \pm } \Lambda^{ -\frac{3}{2}  } e^{ \pm i\Lambda } \rho^{ \pm }
      (t,v) \right\rvert_{ H } ^2 \dif V_H\\
    &\qquad
    + 2\Lambda^3 \Re \int_{H_\Lambda\cap\{t\leq T\}}
      \left\langle \sum_{ \pm } \Lambda^{ -\frac{3}{2}  } e^{ \pm i\Lambda } \rho^{ \pm },
      \err_\Lambda\right\rangle_H \dif V_H\\
    &\qquad
    + \Lambda^3 \int_{H_\Lambda\cap\{t\leq T\}} |\err_\Lambda|_H^2\dif V_H.
  \end{align}
  The quadratic error term is controlled directly from the pointwise bound:
  \begin{align}
    \Lambda^3 \int_{H_\Lambda\cap\{t\leq T\}} |\err_\Lambda|_H^2\dif V_H
    &\lesssim \epsilon^2 \Lambda^3 \Lambda^{-\frac72} T^{2C\epsilon}
    \vol_{\dif V_H}(H_{ \Lambda } \cap \left\{ t\leq T \right\} )\\
    &\lesssim
    \epsilon^2 \Lambda^{ -\frac{5}{2}  } T^{ 2 + 2C\epsilon },
  \end{align}
  where
  \begin{align}
    \vol_{\dif V_H}(H_{ \Lambda }\cap\{t\leq T\})
    \lesssim \Lambda^{-2}T^2.
  \end{align}
  Repeating the same argument for the cross-term, we get
  \begin{align}
     \left\lvert 
       \Lambda^3 \Re \int_{H_\Lambda\cap\{t\leq T\}}
       \left\langle \sum_{ \pm } \Lambda^{ -\frac{3}{2}  } e^{ \pm i\Lambda } \rho^{ \pm },
       \err_\Lambda\right\rangle_H \dif V_H
     \right\rvert  
     \lesssim \epsilon^2 \Lambda^{ -\frac{1}{2}  } T^{ \frac{1}{4} +2C\epsilon },
  \end{align}
  which is still acceptable.
  Moreover, using the orthogonality between $ \rho^{ + }\in V^{ + }  $ and $ \rho^{ - } \in
  V^{ -}  $ in the $ H $-inner product, we can rewrite the remaining term as
  \begin{align}
    \Lambda^3 \int_{H_\Lambda\cap\{t\leq T\}}
      \left\lvert \sum_{ \pm } \Lambda^{ -\frac{3}{2}  } e^{ \pm i\Lambda } \rho^{ \pm }
      (t,v) \right\rvert^2 \dif V_H
    = \int_{ H_{ \Lambda } \cap \left\{ t\leq T \right\}  } \left(
      \lvert \rho^{ + } (t,v) \rvert_{H}^2
      +\lvert \rho^{ - } (t,v) \rvert_{H}^2\right)\dif V_{ H }.
  \end{align} 
  The time variable on $ H_{ \Lambda }  $ is given by
  \begin{align}
    t_{ \Lambda } (v) = \Lambda(1-v^2)^{ -\frac{1}{2}  },
  \end{align} 
  and so we can reparametrize the truncated hyperboloid $ H_{ \Lambda } \cap \{t\leq T_{
  \Lambda } \} $ by the hyperbolic ball
  \begin{align}
    B_{ \Lambda } := \left\{ v\in B(0,1) : t_{ \Lambda } (v)\leq T_{ \Lambda }  \right\} 
    = \left\{ v : v^2\leq 1- \Lambda^2 T_{\Lambda}^{ -2 } \right\}
    \subset \H.
  \end{align} 
  
  To tie everything together, we take  $ T_{ \Lambda }  := \Lambda^{ 1+\eta }  $ with $ \eta
  > 0 $ sufficiently small. 
  Then, as $ \Lambda\to\infty $, the error terms above vanish and
  the  balls $ B_{ \Lambda }  $ saturate $ \H $.
  Moreover, $ t_{ \Lambda } (v)\to \infty $ as $ \Lambda\to\infty $, and so
  \begin{align}
    \lim_{ \Lambda\to\infty } \lvert \rho^{ \pm } (t_{ \Lambda }(v),v)  \rvert _{ H } 
    = \lim_{ \Lambda\to\infty }
      \lvert \widetilde\rho^{\pm}(t_{\Lambda}(v),v)\rvert_{H}
    = |\rho_{\infty}^{\pm}(v)|_{H}.
  \end{align} 
  Moreover, by \eqref{eq:finalrhobound}, we get
  \begin{align}
    |\rho^+(t_\Lambda(v),v)|_H^2+|\rho^-(t_\Lambda(v),v)|_H^2
    \lesssim\epsilon^2(1-v^2)^{2-2C\epsilon}.
  \end{align}
  Since $\dif V_H=(1-v^2)^{-2}\dif v$, the right-hand side is integrable for
  $2C\epsilon<1$.  
  Applying the dominated convergence theorem, we conclude that
  \begin{align}
    \int_{ \R^{ 3 }  }  \lvert \psi_0 \rvert ^2 \dif x
    &= \lim_{ \Lambda\to\infty } \int_{B_\Lambda}
    \left( \lvert \rho^{ + } (t_{ \Lambda } (v),v)\rvert _{ H }^{ 2 }
      + \lvert \rho^{ - } (t_{ \Lambda } (v),v)\rvert _{ H }^{ 2 }  \right) \dif V_{ H } (v)\\
    &= \int_{ H } \left( |\rho_{\infty}^{+}|_{H}^{2}
      + \lvert \rho_{\infty}^{-}\rvert_{H}^{2} \right) \dif V_{ H },
  \end{align} 
  completing the proof.
\end{proof}

\subsection{The wave profile}

We now use the preliminary asymptotics for $\psi$ to identify the leading
profile of the wave component.  
Consider the inhomogeneous equation
\begin{align}
  \square\waveinhom
  = \langle \psi,\gamma^{0}\psi\rangle .
\end{align}
The preliminary Dirac asymptotic expansion gives
\begin{align}
  \psi(t,vt)
  = \Lambda^{-\frac32}
    \sum_{\pm} e^{\pm i\Lambda}e^{\mp i\theta(t,v)}
      \rho_{\infty}^{\pm}(v)
    + \err_{\psi}(t,v),
\end{align}
which we now use to rewrite the source term $ \langle \psi, \gamma^{ 0 }\psi \rangle $ in
terms of the preliminary asymptotic profiles $ \rho_{\infty}^{\pm} $.
Here, the error term $ \err_{ \psi } $ has better decay in $ t $ and shall not affect the
leading-order decay of $ \phi $.

Since the profiles satisfy
$ (P^{+}-P^{-})\rho_{\infty}^{\pm}
  = \mp \rho_{\infty}^{\pm}$,
by \eqref{eq:H-Euclidean-duality} we get
\begin{align}
  \langle \psi, \gamma^{ 0 }\psi \rangle
  &= \sum_{\pm_1,\pm_{2}}\Lambda^{-3}
    \left\langle e^{\pm_{1}i\Lambda \mp_1 i\theta}
      \rho_{\infty}^{\pm_1},
    \gamma^{0}e^{\pm_{2}i\Lambda \mp_{ 2 } i\theta}
      \rho_{\infty}^{\pm_2}\right\rangle
    + \err_{ \phi }\\
  &= \Lambda^{-3}\left(\lvert\rho_{\infty}^{-}\rvert_{H}^{2}
    -\lvert\rho_{\infty}^{+}\rvert_{H}^{2}\right)
    + \err_{ \phi },
\end{align} 
where $ \err_{ \phi } $ is the quadratic contribution generated by $ \err_{ \psi } $.

We set
\begin{align}
  b_{\infty}(v)
  := \lvert\rho_{\infty}^{-}(v)\rvert_{H}^{2}
    - \lvert\rho_{\infty}^{+}(v)\rvert_{H}^{2}.
\end{align}

It remains to solve the wave equation
\begin{align}\label{eq:phi-main-source}
  \square\phi
  = \Lambda^{-3}b_{\infty}(v) + \err_{ \phi }.
\end{align}
As we will show below, $ \err_{ \phi } $ does not affect the leading asymptotics of $ \phi
$, and so we focus on the main equation
\begin{align}
  \square\phimain
  &= \Lambda^{ -3 }b_{ \infty }(v).
\end{align} 
The source is homogeneous of degree $-3$.
Following the ansatz in
\cite{Tataru_hyperbolic},
we look for solutions which are homogeneous of degree $-1$. 
Set
\begin{align}
  \phimain(t,x)
  = \Lambda^{-1}\phi_{\infty}(v),
  \qquad v=\frac{x}{t}.
\end{align}
Applying \eqref{eq:wave-operator-hyperbolic} to $ \Lambda^{-1}a(v) $, we obtain
\begin{align}\label{eq:hyperbolic-wave-identity}
  \square\left(\Lambda^{-1}a(v)\right)
  =
  \Lambda^{-3}
  \left(-1-\Delta_{\H}\right)a(v).
\end{align}
Thus, $\phi_{\infty}$ is determined by the elliptic equation
\begin{align}
  \left(-1-\Delta_{\H}\right)\phi_{\infty}
  = b_{\infty}.
\end{align}
Below, we establish the relevant bounds for the formal computations above.

\begin{prop}\label{prop:wave-profile}
  Set
  \begin{align}
    \phi_{\infty}(v)
    := \left(-1-\Delta_{\H}\right)^{-1}b_{\infty}(v),
    \quad
    \lvert v \rvert  < 1.
  \end{align}
  Then, inside the light cone,
  \begin{align}\label{eq:asymptotic-phi}
    \phi(t,vt)
    = \Lambda^{-1}\phi_{\infty}(v) + O\left(\epsilon\Lambda^{-1}t^{-\delta}\right).
  \end{align}
  Moreover,
  \begin{align}\label{eq:Holder-phi}
    \lVert \phi_{ \infty }\rVert_{ C_{ 0 }^{ \frac{1}{2} }(\overline B) }
    &\lesssim \epsilon ^2,
  \end{align} 
  where $ \bar B $ denotes the unit ball $ \left\{ v\in \R^{ 3 } : \lvert v \rvert
  \leq 1\right\}  $ equipped with the Euclidean metric.
\end{prop}

\begin{remark}
  The elliptic regularity bounds for \((-1-\Delta_{\mathbb H})\) are more
  naturally stated in Hölder spaces of the type \(C^\alpha(\mathbb H)\).
  However, since \(|v|=1\) lies at infinite hyperbolic distance, regularity in
  this geometry alone does not imply continuity up to the Euclidean boundary.
  Hence, we establish our main regularity estimates in \textit{weighted} Hölder
  spaces, e.g. \eqref{eq:H-reg-bounds}, which also give the decay
  \(|\phi_\infty(v)|\lesssim\epsilon^2(1-|v|^2)^{1/2}\). Combining these bounds
  via Lemma \ref{lemma:holder-tradeoff} yields \(\phi_\infty\in
  C^{1/2}(\overline B)\), with \(\phi_\infty=0\) on \(\partial B\).
\end{remark}

\begin{remark}[Linear scattering for the scalar wave]\label{rmk:wave-scattering}
  As mentioned in Remark \ref{rmk:wave-scattering-intro}, the asymptotic formula
  \eqref{eq:asymptotic-phi} \textit{is compatible with linear scattering}, at
  least  if measured in a weak enough topology. 
  Indeed, a simple bilinear argument using the (dual) Sobolev embedding $ L_{ x
  }^{ \frac{3}{2-\eta}  } \hookrightarrow \dot H_{ x }^{ -\frac{1}{2} + \eta }
  $ shows that
  \begin{align}
    \lVert \square\phi \rVert_{\dot H_{ x }^{ -\frac{1}{2} +\eta }  } 
    &\lesssim \lVert \langle \psi, \gamma^{ 0 }\psi \rangle\rVert_{ L_{ x }^{
      \frac{3}{2-\eta}  }  }
      \lesssim \lVert \psi \rVert_{ L_{ x }^{ \frac{6}{3-2\eta}  }  } \lVert \psi
      \rVert_{ L_{ x }^{ 6 }  } 
    \lesssim \langle t \rangle  ^{ -1-\eta },
  \end{align} 
  where the decay rate $ \langle t \rangle  ^{ -1-\eta } $ follows 
  by interpolating the $ L_{ x}^{ 6 } $ decay $ O(t^{ -1 }) $ with $ L_{ x }^{
  2 }  $ charge conservation.
  When $ \eta > 0 $, $ \langle t \rangle  ^{ -1-\eta } $ is integrable and one
  obtains scattering in $ \dot H_{ x }^{ \frac{1}{2} + \eta } \times \dot H_{ x
  }^{ -\frac{1}{2} + \eta }  $ from Duhamel's formula. 
  In particular, while scattering in the (free) wave energy norm $ \dot
  H_{ x }^{ 1 } \times L_{ x }^{ 2 }  $ follows, scattering in the \eqref{dkg}
  scaling-invariant norm $ \mathcal H^{ 0 } $ does not.\footnote{Recall that the $
    \mathcal H^{ 0 } $ norm of the scalar term corresponds to
  $ \eta = 0 $.}
  
\end{remark}

We also define the time-dependent principal term
\begin{align}
  b(t,v)
  &:=
  \lvert \rho^{ - }(t,v) \rvert_{ H }^{ 2 } - \lvert \rho^{ + }(t,v) \rvert _{ H }^{ 2 }.
\end{align}  
Below, we first establish bounds on $ b(t,v) $ for finite $ t $ which we
then transfer to $ b_{ \infty } $.
Moreover, the Hölder regularity bounds \eqref{eq:Holder-phi} shall follow from
vector field bounds on $ b_{ \infty } $ combined with a Morrey-type elliptic
estimate.

\begin{lemma}\label{lemma:b-properties}
  We have the following bounds:
  \begin{align}\label{eq:b-pointwise-convergence}
    \lvert b(t,v) - b_{ \infty }(v) \rvert  
    &\lesssim \epsilon^2 (1-v^2)^{ \frac{3}{4} -C\epsilon } t^{ -\frac{1}{2} +
    C\epsilon }.
  \end{align} 
  \begin{align}\label{eq:interpolated-bound}
    \lVert (1-v^2)^{ -1}\Omega^{ \leq 2 } b_{ \infty } \rVert_{ L^{ 2 }(\H)
    } \lesssim \epsilon^2.
  \end{align} 
\end{lemma}

\begin{proof}[Proof of Lemma \ref{lemma:b-properties}]
  Since the vector fields $ \Omega $ are tangent to the hyperboloids $ H_{ \Lambda }
  $, we use hyperbolic coordinates $ (\Lambda,v) $.
  We recall the bounds on $ \rho^{ \pm } $ in \eqref{eq:Hnormrho},
  \eqref{eq:apriorirhoL2}, and \eqref{eq:rho-infty-bounds} written in coordinates
  $ (\Lambda,v)$,
  \begin{nalign}\label{eq:convergence-summary}
    \lvert \rho^{ \pm }(\Lambda,v) \rvert _{ H }
    &\lesssim
      \epsilon (1-v^2)^{ 1-C\epsilon }, \\
    \lVert \widehat\Omega^{ \leq 3 }\rho^{ \pm }(\Lambda) \rVert_{ L^{ 2 }(\H_{ <R
    }) }
    &\lesssim
      \epsilon \Lambda^{ C\epsilon }, \\
    \lvert \widetilde\rho^{\pm}(\Lambda,v)
      - \rho_{\infty}^{\pm}(v) \rvert _{ H }
    &\lesssim
      \epsilon \Lambda^{ -\frac{1}{2} +C\epsilon }
    (1-v^2)^{-C\epsilon },
  \end{nalign}
  where
  \begin{align}
    \H_{<R}
    &= \left\{v\in B(0,1):(1-v^2)^{-\frac12}\lesssim R\right\},\\
    \H_R
    &= \left\{v\in B(0,1):(1-v^2)^{-\frac12}\sim R\right\}.
  \end{align} 
  Here we also assume  $ \Lambda \gtrsim R^{ 19 } $ so that $ \H_{ <R }  $ is
  fully contained in the pullback of $ \Din\cap H_{ \Lambda } $ to $ \H $, 
  see Appendix~\ref{sec:hyperboloids} for more details.
  The pointwise bound \eqref{eq:b-pointwise-convergence} then follows immediately
  from \eqref{eq:convergence-summary}.
  
    Integrating \eqref{eq:b-pointwise-convergence} over the region where
    $t\geq t_0(v)$ and using the a priori bound \eqref{eq:rho-apriori-growth} on
    its complement, we obtain
    \begin{align}\label{eq:L2-b-difference}
      \lVert (1-v^2)^{ -1 }(b(\Lambda)-b_{ \infty }) \rVert_{ L^{ 2 }(\H) } 
      &\lesssim \epsilon ^2 \Lambda^{ -\frac{1}{20} }.
    \end{align} 

    Now, for the vector field bound \eqref{eq:interpolated-bound}, we proceed by
    interpolation.
    We start with a local $ L^{ 2 } $ bound on the hyperbolic shell $ \H_{ R }  $
    \begin{nalign}\label{eq:L2-b}
      \left\lVert (1-v^2)^{-1}b(\Lambda) \right\rVert_{L^2(\H_{R})}
      &\lesssim
        \epsilon^2 \left(
          \int_{\H_{R}}
          (1-v^2)^{2-C\epsilon}
          \dif V_H
        \right)^{\frac12} \\
      &\lesssim \epsilon^2 R^{-1+C\epsilon}.
    \end{nalign}
    
    For the higher-order bounds, we first note that 
    since $\lvert \Omega_v^{\leq3}(1-v^2)^{-1}\rvert\lesssim
    (1-v^2)^{-1}$, the vector fields may be passed through the weight at the
    expense of lower-order terms of the same form.
    First, we seek to control $ \lVert (1-v^2)^{ -1 }\langle
    \widehat\Omega\rho, \widehat\Omega^{ \leq 2 }\rho \rangle \rVert_{ L^{ 2 }(\H)
    }  $.
    By Hölder and Gagliardo--Nirenberg on $\H_R$, we get
    \begin{align}
      &\lVert (1-v^2)^{-1}
        \langle\widehat\Omega\rho(\Lambda),
        \widehat\Omega^{\leq2}\rho(\Lambda)\rangle_H
        \rVert_{L^2(\H_R)} \\
      &\qquad\lesssim
        \lVert(1-v^2)^{-1}\rVert_{L^\infty(\H_R)}
        \lVert\widehat\Omega\rho(\Lambda)\rVert_{L^4(\H_R)}
        \lVert\widehat\Omega^{\leq2}\rho(\Lambda)\rVert_{L^4(\H_R)} \\
      &\qquad\lesssim
        \lVert (1-v^2)^{-1}\rho(\Lambda)\rVert_{L^\infty(\H_R)}
        \lVert\widehat\Omega^{\leq3}\rho(\Lambda)\rVert_{L^2(\H_R)}.
    \end{align}
    Thus, by \eqref{eq:convergence-summary}, we have the vector field estimate
  \begin{nalign}\label{eq:vf-b}
    \lVert (1-v^2)^{-1}\Omega^{\leq 3} b(\Lambda) \rVert_{L^2(\H_R)}
    &\lesssim
      \lVert (1-v^2)^{-1} \Omega^{\leq 3} \langle
      \rho(\Lambda),\rho(\Lambda)\rangle_H \rVert_{L^2(\H_R)} \\
    &\lesssim \lVert (1-v^2)^{-1} \langle \rho(\Lambda), \widehat\Omega^{\leq
      3}\rho(\Lambda) \rangle_H \rVert_{L^2(\H_R)} \\
    &\hspace{6em}
      + \lVert (1-v^2)^{-1} \langle \widehat\Omega\rho(\Lambda),
      \widehat\Omega^{\leq 2}\rho(\Lambda) \rangle_H \rVert_{L^2(\H_R)} \\
    &\lesssim
      \lVert (1-v^2)^{-1}\rho(\Lambda) \rVert_{L^\infty(\H_R)} \lVert
      \widehat\Omega^{\leq 3}\rho(\Lambda) \rVert_{L^2(\H_R)} \\
    &\lesssim
      \epsilon^2\Lambda^{C\epsilon}.
  \end{nalign}
  Interpolating \eqref{eq:L2-b-difference} with \eqref{eq:vf-b} for the (Cauchy)
  sequence $(b(2\Lambda)-b(\Lambda))$ gives, for $\Lambda\gtrsim R^{19}$,
  \begin{align}
    &\lVert (1-v^2)^{-1}\Omega^{\leq2}
      (b(2\Lambda)-b(\Lambda))\rVert_{L^2(\H_R)} 
    \lesssim
      (\epsilon^2\Lambda^{-\frac1{20}})^{\frac13}
      (\epsilon^2\Lambda^{C\epsilon})^{\frac23}
      \lesssim \epsilon^2\Lambda^{-\frac1{60}+C\epsilon}.
  \end{align}
  Summing in dyadic $\Lambda$ shows that, for each fixed $R$,
  \begin{align}
    \lVert (1-v^2)^{-1}\Omega^{\leq2}
      (b(\Lambda)-b_\infty)\rVert_{L^2(\H_R)}
    &\lesssim \epsilon^2\Lambda^{-\frac1{60}+C\epsilon}
    \longrightarrow0.
  \end{align}
  The growth in the vector field bound prevents global $ H^{ 2 } $ convergence;
  however, we still retain boundedness.
  For each dyadic $ R $, we set $ \Lambda_{ R }  := R^{ 19 } $ and interpolate
  \eqref{eq:L2-b} and \eqref{eq:vf-b} to obtain
  \begin{align}
    \lVert (1-v^2)^{ -1 }\Omega^{ \leq 2 }b(\Lambda_{ R } ) \rVert_{ L^{ 2 }(\H_{
      R } ) }
    &\lesssim
      (\epsilon^2R^{-1+C\epsilon})^{\frac13}
      (\epsilon^2\Lambda_R^{C\epsilon})^{\frac23}
      \lesssim \epsilon ^2 R^{ -\frac{1}{4} }.
  \end{align} 
  The preceding Cauchy estimate also gives
  \begin{align}
    \lVert (1-v^2)^{-1}\Omega^{\leq2}
      (b(\Lambda_R)-b_\infty)\rVert_{L^2(\H_R)}
    &\lesssim \epsilon^2R^{-\frac{19}{60}+C\epsilon}
    \lesssim \epsilon^2R^{-\frac14},
  \end{align}
  provided $\epsilon$ is sufficiently small.
  Finally,
  \begin{align}
    \lVert (1-v^2)^{ -1 }\Omega^{ \leq 2 }b_{ \infty } \rVert_{ L^{ 2 }(\H) } 
    &\lesssim \sum_{ R \text{ dyadic} } \left(\lVert (1-v^2)^{ -1 }\Omega^{ \leq 2
      }b(\Lambda_{ R } ) \rVert_{ L^{ 2 }(\H_{ R } ) }\right. \\
    &\hspace{7em}\left.+ \lVert (1-v^2)^{ -1 }\Omega^{ \leq 2 }(b(\Lambda_{ R })-b_{
      \infty } ) \rVert_{ L^{ 2 }(\H_{ R } ) }\right)\\
    &\lesssim \sum_{ R \text{ dyadic} } \epsilon^2 R^{ -\frac{1}{4}  }\\
    &\lesssim \epsilon ^2,
  \end{align} 
  completing the proof.

\end{proof}

\begin{proof}[Proof of Proposition \ref{prop:wave-profile}]
  Using the notation from \eqref{eq:phi-main-source}, 
  we start by establishing the bounds on $ \phimain $.
  Since $ H^{ 2 }(\H)$ embeds into $C^{ \frac{1}{2}  }(\H) $ continuously,
  we immediately get
  \begin{align}
    \lVert (1-v^2)^{ -1 }b_{ \infty } \rVert_{ C^{ \frac{1}{2}  }(\H) } 
    &\lesssim \epsilon ^2.
  \end{align} 
  Now, by \eqref{eq:interpolated-bound}, the source term 
  \begin{align}
    G_{ \text{main} } (t,v) = \1_{ \{t>r\} }\Lambda^{ -3 }b_{ \infty }(v)
  \end{align}
  satisfies bounds similar to those in \eqref{eq:sourcebounds}, 
  and, repeating the same argument, we obtain
  \begin{align}
    \lvert \phimain (t,vt) \rvert 
    = t^{ -1 }(1-v^2)^{ -\frac{1}{2}  }\lvert \phi_{ \infty }(v) \rvert 
    &\lesssim \epsilon ^2 t^{ -1 }.
  \end{align} 
  Therefore, by local elliptic regularity,
  \begin{align}\label{eq:H-reg-bounds}
    \lVert (1-v^2)^{ -\frac{1}{2}  }\phi_{ \infty } \rVert_{ C^{ 2,\frac{1}{2}
      }(\H) } 
    &\lesssim \lVert (1-v^2)^{ -\frac{1}{2} } \phi_{ \infty } \rVert_{ L^{
        \infty }(\H) } + \lVert (1-v^2)^{ -\frac{1}{2} } b_{ \infty } \rVert_{
        C^{ \frac{1}{2}  }(\H) } 
    \lesssim \epsilon^2.
  \end{align} 
  The Hölder bounds in the Euclidean metric follow directly from Lemma
  \ref{lemma:holder-tradeoff}.

  As for the error term, on $ \Din $ we use Lemma~\ref{asympdecomp} and
  \eqref{eq:rho-infty-bounds} to directly estimate
  \begin{align}
    \lvert \err_{ \phi }(t,v) \rvert  
    \lesssim \lvert \psi \rvert \lvert \err_{ \psi } \rvert 
    \lesssim \epsilon^2 t^{ -3-\delta } 
  \end{align} 
  for a small enough $ \delta > 0 $.
  Repeating the proof of Proposition \ref{prop:wave-inhom-bound} again,
  we can shrink $ \delta $ and obtain
   \begin{align}
    \lvert \square^{ -1 }\err_{ \phi } \rvert 
    &\lesssim \epsilon t^{ -1-\delta }.
  \end{align}
  Finally, on $ \left\{ t > r \right\} \setminus \Din $, 
  the improved source bounds \eqref{eq:sourcebounds} suffice for
  establishing \eqref{eq:asymptotic-phi}.
\end{proof}

\begin{remark}\label{rmk:refine-sobolev}
  The choice of two vector fields in \eqref{eq:interpolated-bound} is not sharp,
  but it suffices for our purposes.
  Since $ \lVert (1-v^2)^{ -1 }\Omega^{ \leq 3 } b_{ \infty } \rVert_{ L^{ 2
  }(\H_{ R } ) } $ grows in $ R $, 
  we could use the same interpolation as above to show uniform boundedness
  \begin{align}
    \lVert (1-v^2)^{ -1 }\Omega^{ \leq 3-\delta } b_{ \infty }\rVert_{ L^{ 2 }(\H) } 
    &\lesssim_{ \delta } \epsilon^2,
  \end{align} 
  for any $ \delta > 0 $, but this fails at $ \delta = 0 $.
\end{remark}

\begin{remark}
    While $ b_{ \infty } $ has strong decay $ O\left( \epsilon^2(1-v^2)^{ 2-C\epsilon }
    \right)  $ near $ \lvert v \rvert  = 1 $, the solution $ \phi_{ \infty } $
    only decays with $ O( (1-v^2)^{ \frac{1}{2} } )  $ instead.
    This is sharp, as can be seen by 
    directly solving the elliptic equation
    \begin{align}\label{eq:ellipticradial}
      (-1-\Delta_{ \H })u = f
    \end{align} 
    for a (radial) function $ f $  depending only on $ \rho = d_H(0,v) $, see
    Appendix \ref{sec:hyperboloids}.
    For $ f $ sufficiently regular and with enough decay near $ \lvert v \rvert  =
    1$, we write
    \begin{align}
      \Delta_{ \H } u = \partial_{ \rho }^2 u + 2\coth\rho \partial_{ \rho }u,   
    \end{align} 
    and so \eqref{eq:ellipticradial} simplifies to
    \begin{align}
      -\partial_{ \rho }^2u -2\coth\rho \partial_{ \rho }u - u = f.
    \end{align}
    If we define $ U(\rho) := \sinh\rho\cdot u $, \eqref{eq:ellipticradial} can be
    simplified further to 
    \begin{align}
      -\partial_{ \rho }^2 U = \sinh\rho f.
    \end{align} 
    The unique solution to this ODE with $U(0)=0$ and bounded as
    \(\rho\to\infty\) is
    \begin{align}
      U(\rho) 
      &= \int_{ 0 }^{ \infty } \min\left\{ \rho,\rho' \right\} \sinh\rho'
        f(\rho') \dif \rho'\\
      &= \int_{ 0 }^{ \rho } \rho'\sinh\rho' f(\rho') \dif \rho'
      + \rho \int_{ \rho }^{ \infty }\sinh\rho' f(\rho') \dif \rho'.
    \end{align} 
    Denoting $ M := \int_{ 0 }^{ \infty } \rho'\sinh\rho' f(\rho')\dif \rho' $,
    which we assume to be finite, we can rewrite
    \begin{align}\label{eq:ellipticsol}
      u(\rho) = \frac{ M }{ \sinh\rho } 
      - \frac{ 1 }{ \sinh\rho } \int_{ \rho }^{ \infty}
      (\rho'-\rho)\sinh\rho'f(\rho')\dif \rho'.
    \end{align} 
    Now, while the second term in \eqref{eq:ellipticsol} can be made smaller by
    imposing stronger boundary decay on $ f $, the leading term satisfies
    \begin{align} 
      \frac{M}{\sinh\rho} =  M\frac{(1-v^2)^{ \frac{1}{2}}}{\lvert v \rvert} ,
    \end{align} 
    see \eqref{eq:d_H}, and cannot be improved further.
    This can also be seen at the level of the kernel $ K $ of $
    (-1-\Delta_{ \H })^{ -1 } $, which is given explicitly by
    \begin{align}
      (-1-\Delta_{ \H })^{ -1 }b(v) 
      = \int_{ \H } K(d_H(v,v')) b(v') \dif V_H(v'),
      \qquad
      K(\rho) = \frac{1}{4\pi \sinh\rho},
    \end{align} 
    see \cite{Tataru_hyperbolic}.
\end{remark}

\subsection{The Dirac profile}

In subsection \ref{sec:first-dirac}, we solved the asymptotic ODE for the
profiles $ \rho^{ \pm } $ using a preliminary choice of phase.  Having now
identified the wave profile $ \phi_{ \infty } $, we return to the Dirac
profiles, fix this choice canonically, and prove regularity bounds for $
\rho^{ \pm }$.

In \eqref{eq:renormalized-rho}, we used the preliminary normalization
\begin{align}
  \widetilde\rho^{\pm}(t,v)
  = e^{\pm i\theta(t,v)}\rho^{\pm}(t,v),
  \qquad
  \theta(t,v) = \int_{ t_0(v) }^{ t } \sqrt{ 1-v^2 } \phi(s,vs)\dif s.
\end{align} 
Using \eqref{eq:asymptotic-phi}, we can write 
\begin{align}
  \widetilde\rho^{\pm}(t,v)
  \approx 
  e^{ \pm i \phi_{ \infty }(v)\log t }\rho^{ \pm }(t,v)e^{ i\vartheta }
\end{align} 
for some small error $ \vartheta $. 
We use this heuristic to choose the canonical normalization below, retaining the
notation $\widetilde\rho^{\pm}$ and $\rho_{\infty}^{\pm}$, and effectively initializing
$ \vartheta(t=\infty) = 0 $.
The reason we use $ \log (t^2-x^2) $ instead of $ \log t $ is to preserve the
Lorentz invariance, which will be convenient when establishing the regularity
bounds below.

\begin{prop}
  Let
  \begin{align}
    \widetilde\rho^{ \pm }(t,v)
    &:= e^{\pm \frac{i}{2}\log(t^2-\lvert vt\rvert^2)\phi_\infty(v)}
      \rho^{\pm}(t,v).
  \end{align} 
  For each fixed $ v\in B(0,1) $, the limit 
  \begin{align}\label{eq:limit}
    \rho_{ \infty }^{ \pm }(v)
    &:= \lim_{ t\to\infty } \widetilde\rho^{ \pm }(t,v)
  \end{align} 
  exists.
  Moreover, the redefined profiles satisfy the bounds
  \begin{align}
    \lvert \rho^{\pm}(t,v)
      - e^{\mp \frac{i}{2}\phi_{\infty}(v)
        \log(t^2-\lvert vt\rvert^2)}\rho_{\infty}^{\pm}(v) \rvert_{H}
    &\lesssim \epsilon \langle t-r \rangle^{ -\delta },\label{eq:difference-rho}\\
      \lVert (1-v^2)^{-\eta}\widehat\Omega^{\leq\frac52}
        \rho_{\infty}^{\pm}
      \rVert_{ L^{ 2}(\H) } 
    &\lesssim \epsilon\label{eq:energy-rho},
  \end{align} 
  for some sufficiently small $ \delta,\eta>0$.
  The profiles also satisfy the Hölder bounds
  \begin{align*}
    \lVert (1-v^2)^{-\frac14-\eta}\rho_\infty^\pm
      \rVert_{C^{1-\eta}(\H)}
    + \lVert \rho_\infty^\pm
      \rVert_{C^{\frac12}(\overline{B(0,1)})}
    &\lesssim \epsilon,
    \qquad
    \rho_\infty^\pm|_{\partial B(0,1)}=0.
  \end{align*}
\end{prop}

\begin{proof}
  Using the asymptotic equation \eqref{eq:asymp} for $ \rho^{ \pm } $ and the
  asymptotic expansion for $ \phi $ \eqref{eq:asymptotic-phi}, we compute
  \begin{align}
     i\partial_t \widetilde\rho^{ \pm }(t,v)
    &=
      e^{\pm \frac{i}{2}\log(t^2-(vt)^2)\phi_\infty(v)}
      \left(
        i\partial_t\rho^\pm(t,v)
        \mp \frac{\phi_\infty(v)}{t}\rho^\pm(t,v)
      \right)\\
    &=
      \pm \sqrt{1-|v|^2}
        \left(
          \phi(t,vt)-\Lambda^{-1}\phi_\infty(v)
        \right)
        \widetilde\rho^\pm(t,v)
      + O_H\left(
        \epsilon t^{-1+C\epsilon}\Lambda^{-\frac12}
      \right)\\
      &= O_{ H } \left( \epsilon t^{ -1-\delta }(1-v^2)^{ -\frac{1}{4}  }
    \right) .
  \end{align}
  Since the error is integrable along each ray $ x = vt $, the limit
  \eqref{eq:limit} exists, and integrating the error yields the difference bound
  \eqref{eq:difference-rho}.

  For the remaining regularity estimates, we fix a branch and suppress the label $\pm$.
  As for the vector field bounds \eqref{eq:energy-rho}, we proceed as in Lemma
  \ref{lemma:b-properties}.
  This requires some care, since $ b_{ \infty } $ had much stronger decay at the
  boundary.
  
  It is sufficient to establish the global bound
  \begin{align}
    \lVert (1-v^2)^{-\eta}
      \widehat\Omega^{\leq\frac52}\rho_{ \infty }
    \rVert_{L^2(\H)}
    &\lesssim \epsilon
  \end{align}
  for some $ \eta > 0 $ sufficiently small. Indeed,
  by Morrey's inequality (Lemma \ref{lemma:morrey}), we get
  $(1-v^2)^{-\frac{1}{4} -\eta}\rho_\infty\in C^{1-\eta}(\H)$,
  and Lemma \ref{lemma:holder-tradeoff} then gives
  $\rho_\infty\in C^{\frac12}(\overline{B(0,1)})$ with  
  $\rho_\infty|_{\partial B}=0$.
 
  As before, we start by obtaining estimates on hyperbolic shells $ \H_{ R }  $ 
  for finite $ \Lambda $-time, use interpolation to control higher regularity
  norms, and take appropriate limits as $ \Lambda,R\rightarrow\infty $ in
  order to pass these bounds to $ \rho_{ \infty } $. 
  The first observation is that 
  by refining Lemma \ref{lemma:b-properties} slightly, see Remark
  \ref{rmk:refine-sobolev},  we have
  \begin{align}
    \lVert \Omega^{\leq 3-\eta}\phi_\infty\rVert_{L^2(\H_R)}
    \lesssim \epsilon^2.
  \end{align}
  Since $H^{3-\eta}(\H_{ R } )$ is an algebra, \eqref{eq:apriorirhoL2} gives
  \begin{align}\label{eq:weighted-sobolev-tilde-rho}
    \lVert \widehat\Omega^{ \leq 3-\eta }\widetilde\rho(\Lambda) \rVert _{ L^{
      2 }(\H_{ R }) }
    \lesssim \epsilon\Lambda^{ C\epsilon },
  \end{align}
  where, as in \eqref{eq:convergence-summary}, we pick $ \Lambda \gtrsim R^{ 19 } $
  so that $ \H_{ R } $ is contained in the pullback of $ \Din\cap H_{ \Lambda }$
  to $ \H $.
  On the other hand, the pointwise bound in \eqref{eq:convergence-summary} and the
  shell volume give
  \begin{align}
    \lVert (1-v^2)^{ -\eta }\widetilde\rho(\Lambda) \rVert _{ L^{ 2 }(\H_{ R } ) } 
    &\lesssim R^{ 2\eta } \lVert \widetilde\rho(\Lambda) \rVert_{ L^{ \infty }(\H_{
    R}) } \vol_{ \dif V_{ H }  }(\H_{ R } )^{ \frac{1}{2}  }
    \lesssim \epsilon R^{ 2\eta-1+C\epsilon }.
  \end{align}
  Since multiplication by $(1-v^2)^{ -\eta }$ costs $R^{ 2\eta }$
  on $\H_R$, including its Lorentz derivatives,
  \eqref{eq:weighted-sobolev-tilde-rho} gives
  \begin{align}
    \lVert (1-v^2)^{-\eta}
      \widehat\Omega^{ \leq 3-\eta } \widetilde\rho(\Lambda)
      \rVert_{L^{ 2 }(\H_R)}
    &\lesssim \epsilon R^{2\eta}\Lambda^{C\epsilon}.
  \end{align}
  Pick $ \eta > 0 $ sufficiently small and interpolate
  \begin{align}
    \lVert (1-v^2)^{-\eta}
     \widehat\Omega^{ \leq \frac{5}{2}  }
     \widetilde\rho(\Lambda) \rVert_{L^{ 2 }(\H_R)}
    &\lesssim
    \left( 
      \epsilon R^{ 2\eta - 1+ C\epsilon } 
    \right) ^{ \frac{1}{8}  }
    \left( 
      \epsilon R^{ 2\eta }\Lambda^{C\epsilon}
    \right) ^{ \frac{7}{8}  }\\
    &\lesssim
      \epsilon R^{2\eta - \frac{1}{8} + C\epsilon }\Lambda^{C\epsilon}.
  \end{align}
  Choosing $\Lambda_R=R^M$ for sufficiently large $ M > 0 $, we use $ R
  $-decay to absorb the  $ \Lambda $-growth and obtain 
  \begin{align}\label{eq:rho-shell-finite-time}
    \lVert (1-v^2)^{-\eta}
      \widehat\Omega^{ \leq \frac{5}{2}  } \widetilde\rho(\Lambda_R) \rVert_{L^{ 2 }(\H_R)}
    &\lesssim \epsilon R^{-\frac{\eta}{4} }.
  \end{align}
  It remains to pass this estimate to the limit.  
  For dyadic $\Lambda\geq \Lambda_R$, we have $(t-r)\sim\Lambda/R$ on $\H_R$.  Thus
  \eqref{eq:difference-rho}, \eqref{eq:weighted-sobolev-tilde-rho}, and the
  shell volume give
  \begin{align}
    \lVert (1-v^2)^{-\eta}
      (\widetilde\rho(2\Lambda)-\widetilde\rho(\Lambda))
      \rVert_{L^2(\H_R)}
    &\lesssim \epsilon R^{2\eta+1+\delta}\Lambda^{-\delta}, \\
    \lVert (1-v^2)^{-\eta}\widehat\Omega^{\leq 3-\eta}
      (\widetilde\rho(2\Lambda)-\widetilde\rho(\Lambda))
      \rVert_{L^2(\H_R)}
    &\lesssim \epsilon R^{2\eta}\Lambda^{C\epsilon}.
  \end{align}
  Interpolation yields
  \begin{align}
    \lVert (1-v^2)^{-\eta}\widehat\Omega^{\leq\frac52}
      (\widetilde\rho(2\Lambda)-\widetilde\rho(\Lambda))
      \rVert_{L^2(\H_R)}
    &\lesssim \epsilon R^{2\eta+\frac{1+\delta}{8}}
      \Lambda^{-\frac\delta8+C\epsilon}.
  \end{align}
  Hence, for $\epsilon$ sufficiently small, 
  \begin{align}\label{eq:rho-shell-tail}
    \lVert (1-v^2)^{-\eta}\widehat\Omega^{\leq\frac52}
      (\rho_\infty-\widetilde\rho(\Lambda_R))
      \rVert_{L^2(\H_R)}
    &\lesssim \epsilon R^{2\eta+\frac{1+\delta}{8}}
      \Lambda_R^{-\frac\delta{16}}
    \lesssim \epsilon R^{-\frac\eta4}.
  \end{align}
  Combining \eqref{eq:rho-shell-finite-time} and \eqref{eq:rho-shell-tail} and summing over the dyadic shells completes
  the proof.

\end{proof}

\subsection{The Majorana class and linear scattering} 

We recall that $\psi_0$ satisfies the \textit{Majorana condition} if there
is a fixed $z\in\C$, $\lvert z\rvert=1$, such that
\begin{align}
  \psi_0+z\gamma^2\psi_0^\ast &= 0.
\end{align}
A direct computation shows that this condition is preserved by the Dirac
evolution with a real-valued scalar potential, see \cite{CandyHerr2018Majorana}. 
More precisely, if $\psi$ is a global solution to
\begin{align}
  (-i\gamma^\mu\partial_\mu+1)\psi &= V\psi,
\end{align}
where $V$ is real-valued, and $\psi(0)$ satisfies the Majorana condition for
some fixed $z$, then
\begin{align}\label{eq:majorana-evolution}
  \psi(t)+z\gamma^2\psi^*(t) &= 0,
\end{align}
for all $t>0$. Moreover, it can be directly checked that Majorana solutions satisfy
\begin{align}\label{eq:majorana-source-cancellation}
  \langle\psi,\gamma^0\psi\rangle &= 0.
\end{align}

\begin{proof}[Proof of Proposition~\ref{prop:majorana-scattering}]
  Let $(\psi,\phi)$ be the global solution to \eqref{dkg} with initial data in
  the weighted class of Theorem~\ref{main-thm-gwp}, and suppose that $\psi(0)$
  satisfies the Majorana condition for a fixed $z\in\C$, $\lvert z\rvert=1$.
  Given a spinor $ \rho $, we define
  \begin{align}
    M_z\rho
    &:= z\gamma^2\rho^\ast,
  \end{align}
  and assume $ \psi_0 + M_{ z } \psi_0 = 0 $.
  Recall that we must show
  \begin{align}\label{eq:algebraic}
    \rho_{ \infty }^{ + }(v) + M_{ z } \rho_{ \infty }^{ - }(v) = 0. 
  \end{align} 
  First, we note that $ M_{ z }  $ swaps subspaces $ V^{ \pm } $, that is,
  $ M_{ z } P_{ v }^{ \pm }   = P_{ v }^{ \mp } M_{ z } $.
  This follows directly from computing
  \begin{align}
    M_z\gamma^H\rho
    &= z\gamma^2\overline{\gamma^H}\rho^\ast 
    = -\gamma^H M_z\rho,
  \end{align}
  where we recall that $ \gamma^{ 2 }\overline{\gamma^{ \mu }} = -\gamma^{ \mu 
  }\gamma^{ 2 } $.
  Recall the asymptotic decomposition for $ \psi $ on the truncated hyperboloid
  $ \H_{ <T/\Lambda } $ from \eqref{eq:truncated-dirac-profile-error}
  \begin{align}
    \psi(t,vt)
    &=
    \Lambda^{-\frac32}
    \left(
      e^{i\Lambda-i\phi_\infty(v)\log\Lambda}\rho^+(t,v)
      +
      e^{-i\Lambda+i\phi_\infty(v)\log\Lambda}\rho^-(t,v)
    \right)
    +\err_{ \Lambda },
  \end{align}
  for the approximate profiles $ \rho^{ \pm }(t,v) \in V_{ v }^{ \pm }  $.
  Applying $ M_{ z }  $ to this expression and imposing the Majorana
  cancellation \eqref{eq:majorana-evolution} yields
  \begin{align}
    0
    &= \Lambda^{\frac32}
    \left(\psi+M_z\psi\right) \\
    &=
    e^{i\Lambda-i\phi_\infty(v)\log\Lambda}
    \left(
      \rho^+(t,v)+M_z\rho^-(t,v)
    \right) 
    + e^{-i\Lambda+i\phi_\infty(v)\log\Lambda}
    \left(
      \rho^-(t,v)+M_z\rho^+(t,v)
    \right)\\
    &\hspace{31em}+ \Lambda^{\frac32}
    \err_{ \Lambda }.
  \end{align}
  Projecting onto the two subspaces, we get
  \begin{align}
    \rho^\pm(t,v)+M_z\rho^\mp(t,v)
    &=
    O_H\left(
      \epsilon\Lambda^{-\frac14}T^{C\epsilon}
    \right). 
  \end{align}
  Taking $ \Lambda = T^{ \frac{19}{20} } $ as before and passing to the limit 
  completes the proof.
\end{proof}

As discussed in the introduction, the relation \eqref{eq:algebraic} together
with the asymptotic completeness construction in
section~\ref{sec:completeness} shows that the Majorana class does not exhaust
the class of linearly scattering solutions.

\begin{remark}
  In fact, if we do not require the explicit relation \eqref{eq:algebraic},
  scattering can be established directly from
  \eqref{eq:majorana-source-cancellation} and the asymptotic formula
  \eqref{eq:asymptotic-phi}.
  Indeed, by  \eqref{eq:asymptotic-phi},
  $\phi(t,vt)=O(t^{-1})$ along each fixed timelike ray $\lvert v\rvert<1$.  
  Since $\langle\psi,\gamma^0\psi\rangle=0$, the scalar component $\phi$ is a
  free wave, and it has improved decay along timelike rays 
  $ O(t^{ -1 }\langle t-r \rangle ^{ -\delta }) $ 
  \eqref{eq:free-wave-pointwise}. 
  The two decay rates are only compatible if $ \phi_{ \infty }\equiv
  0 $,
  and we directly conclude that $ \psi $ scatters linearly.
\end{remark}

\section{Asymptotic completeness}\label{sec:completeness}
Our goal in this section is to  establish asymptotic
completeness for a class of high-regularity initial profiles.
We prove the following version of Theorem \ref{main-thm-completeness}, now stated with the full
regularity assumptions.

\begin{prop}[Asymptotic completeness]\label{prop:main-back}
  Fix $0<\nu_0<\frac12$, and let $\epsilon>0$ be sufficiently
  small. Let
  $(\rho_\infty^+,\rho_\infty^-,\Fren)$ be prescribed
  asymptotic data, with $\rho_\infty^\pm(v)\in V_v^\pm$
  and $\Fren:\R_q\times\S_\omega^2\to\R$, satisfying
  \begin{align}
    \sum_\pm\sum_{|\beta|\leq20}
      \bigl\lVert w^{-50+2|\beta|}
        \partial_v^\beta\rho_\infty^\pm
      \bigr\rVert_{L^2(B(0,1),\dif v)}
    &\ll\epsilon,
    \quad w(v):=(1-|v|^2)^{1/2},
  \end{align}
  and
  \begin{align}
    &\sum_{3(k+\ell)+m\leq27}
      \bigl\lVert(q\partial_q)^k\nabla_\omega^\ell
        \partial_q^m\Fren
      \bigr\rVert_{\dot H_q^{1/2}L_\omega^2}
    + \sum_{3(k+\ell)+m\leq27}
      \bigl\lVert(q\partial_q)^k\nabla_\omega^\ell
        \partial_q^m\Fren
      \bigr\rVert_{\dot H_q^{1/2-\nu_0}L_\omega^2}
    \ll\epsilon.
  \end{align}
  Then, there exists a global solution $(\psi,\phi)$
  to \eqref{dkg} whose initial data satisfy the
  hypotheses of Theorem \ref{main-thm-gwp}.
  Inside the forward light cone,  $\psi$ and $\phi$
  have the asymptotic expansions \eqref{eq:psi-asympt-mod} and
  \eqref{eq:phi-asympt-mod}, respectively.

  Moreover, the wave component satisfies
  \begin{align}
    \lim_{r\to\infty}
      r\bigl(\phi-\phit\bigr)(r+q,r\omega)
      =\Fren(q,\omega)
  \end{align}
  for every fixed $q\in\R$ and $\omega\in\S^2$, where
  $\phit$ is the solution with zero Cauchy data at $t=1$
  of
  \begin{align}
    \square\phit(t,x)
      =\1_{ \left\{  \lvert x \rvert  < t \right\}  }(t,x)
      (t^2-|x|^2)^{-3/2}
      \left(
        |\rho_\infty^-(x/t)|_H^2
        -|\rho_\infty^+(x/t)|_H^2
      \right).
  \end{align}
\end{prop}

\begin{remark}
  The number of derivatives was chosen so that our solutions recover the
  forward weighted class, but the construction itself could be done with
  substantially fewer assumptions.
\end{remark}

\begin{remark}
  The different regularity assumptions on the Dirac and wave profiles reflect
  how we estimate the corresponding backgrounds: we use the explicit Dirac
  profile and its nonresonant correction, whereas the free wave is controlled
  by a commuted pointwise estimate. We give a refined derivative and weight count in
  Remark \ref{rmk:backward-regularity-count} after introducing some
  notation.
\end{remark}

The idea is to build the prescribed asymptotics into an approximate solution
$(\psiap,\phiap)$, then solve backward for a small correction. 
Writing
$\Psi=\psi - \psiap$ and $\Phi = \phi - \phiap$, we formally obtain
\begin{nalign}\label{eq:untrunc}
  &(-i\gamma^\mu\partial_\mu+1)\Psi
    =\Phi\Psi+\Phi\psiap+\phiap\Psi-R_\psi,\\
  &\square\Phi
    =\langle\Psi,\gamma^0\Psi\rangle
      +2\Re\langle\Psi,\gamma^0\psiap\rangle-R_\phi,\\
  &(\Psi,\Phi,\partial_t\Phi)\big|_{t=\infty}= 0,
\end{nalign}
where $R_\psi$ and $R_\phi$ are the errors obtained by inserting the
approximate solution into \eqref{dkg}. 
Our goal is then to prove that 
$ (\Psi(t),\Phi[t]) \longrightarrow 0 $ as $ t\to \infty $, in an appropriate
norm, and that both $ (\Psi(T_0),\Phi[T_0]) $ and $ (\psiap(T_0),\phiap[T_0]) $
belong to the weighted class in Theorem \ref{main-thm-gwp}, for a fixed $ T_0
\geq 2 $.

\vspace{3pt}
\noindent
We make this idea precise in three steps:
\begin{enumerate}[(i)]
  \item Solve the difference equations with zero final data on the time 
    slice $ \left\{ t = T \right\} $.
    A contraction argument gives bounds on
    $[T_0,T]$ that are uniform in $T$.

  \item Compare solutions $ (\Psi_{ T },\Phi_{ T } ) $ with different terminal
    times along the slice $ \left\{ t = T_0 \right\} $ using weighted norms.

  \item Pass to the limit $T\to\infty$ in $ C_{ t }^{ 0 } \mathcal H^{ 0 } $ to
    construct a solution to the untruncated equation \eqref{eq:untrunc}.
    This gives the required data bounds on the fixed slice $ \{t = T_0 \}$.  
    Local well-posedness then transfers these bounds to $t=0$.
\end{enumerate}
Thus, the main obstruction is choosing $(\psiap,\phiap)$ so that the residual
terms $ R_{ \psi } $ and $ R_{ \phi } $ decay sufficiently fast as $ t\to\infty
$.

\smallskip
\noindent
\textit{The Dirac background.}
We start with the expansion obtained in the forward theory. In the
hyperboloidal coordinates $(\Lambda,v)$, it suggests the approximation
\begin{align}\label{eq:psiasymp}
  \psia^\pm(\Lambda,v)
    =\Lambda^{-3/2}e^{\pm i\Lambda\mp i\phi_\infty(v)\log\Lambda}
      \rho_\infty^\pm(v),
\end{align}
where $\phi_\infty$ is determined by the prescribed Dirac profiles as in the
forward problem. 

The difficulty is that having the right leading behavior does not give a
small enough residual. Differentiating the amplitude and phase leaves terms
whose $L_x^2$ norms decay only as $t^{-1}$, with logarithmic losses. 
As in the wave-packet construction \eqref{eq:wp-correction}, we can avoid losses by adding a
nonresonant correction to the definition of the profiles. 
We replace $\rho_\infty^\pm$ by
\begingroup
\begin{align}\label{eq:corrected-profiles}
  \varrho^\pm=\rho_\infty^\pm+\Lambda^{-1}r^\pm,
  \qquad r^\pm\in V_v^\mp.
\end{align}
\endgroup
Within the high-regularity profile class from Proposition \ref{prop:main-back}, one should  
think of $ r^{ \pm } $ as being $ O(\log t) $; the factor $ \Lambda^{ -1 } $
keeps the correction lower-order.
$ r^{ \pm } $ is explicit and given in the proof of Lemma
\ref{lemma:background}.
Its definition involves derivatives of $\rho_\infty^\pm$ and $\phi_\infty$.
This derivative cost, together with the later commutations, is, in part, why we work with
high-regularity profiles.

For the backward evolution, we must also address the choice of coordinate
systems.
The $(\Lambda,v)$ coordinates only cover the interior of the cone, whereas the
terminal data must be prescribed on a constant-time slice. 
We therefore keep the velocity variable but return to the physical time $t$:
\begin{align}\label{eq:coorddef}
  v=\frac{x}{t},\qquad
  w:=(1-|v|^2)^{1/2},\qquad
  \Lambda=tw \quad (|v|<1).
\end{align}
We now regard $\psiap^\pm$ as functions of $(t,v)$ inside the cone, and extend
them by zero outside.

Writing $\Lambda=tw$ also exposes a second issue: both the amplitude and
phase correction in \eqref{eq:psiasymp} contain negative powers of $w$, which
become singular as $|v|\to1$.
This is harmless in the forward problem, since the Klainerman-Sobolev bounds
\ref{prop:klainerman-sobolev} are enough to close the bootstrap argument in the
exterior of the inner-cone $ \left\{ t - r = 1\right\}  $, but becomes an issue
once we prescribe $ \psiap $ on a constant-time slice and need to control it
across $ \lvert v \rvert   = 1 $.
We require enough boundary decay of the spinor profiles and their
derivatives to compensate for these factors and justify the zero extension.
This $ w $-decay also provides a crucial gain in the Dirac--wave interaction
near the cone, where the wave itself may decay only as $t^{-1}$ in time, see
\eqref{eq:weighted-free-wave-product}.

\smallskip
\noindent
\textit{The wave background.}
The Dirac profiles determine the leading scalar source
\begin{align}
  b_\infty(v)=|\rho_\infty^-(v)|_H^2-|\rho_\infty^+(v)|_H^2.
\end{align}
As in the forward theory,
$\phi_\infty=(-1-\Delta_\H)^{-1}b_\infty$, so the corresponding timelike
model is $\phimain=(tw)^{-1}\phi_\infty$. This still leaves part of the
wave undetermined: free waves carry nonzero radiation at null infinity
whereas their leading timelike profile vanishes.

To retain both pieces of information, we first construct the $ \psi $-induced
timelike wave $ \phit $ by solving
\begin{align}
  \square\phit=\1_{\{|v|<1\}}(tw)^{-3}b_\infty(v),
  \qquad \phit[1]=0.
\end{align}
This wave agrees with $\phimain$ when $t-r\geq2$ and supplies the transition
near the cone, see Lemma \ref{lemma:wave-asymptotics} for a proof. 
We then let $\phif$ be the free wave with radiation profile
$\Fren$ and define
\begin{align}
  \phiap=\phit+\phif.
\end{align}
Here $ \Fren $ accounts for the $ \phi $-radiation escaping to future
null infinity. 

We organize this section as follows.
In subsection \ref{sec:wave-asymptotics}, we collect the main properties of $ (\psiap,\phiap) $.
In subsection \ref{sec:backward-contraction}, we define the time-truncated difference equations and
solve them using a contraction argument.
In subsection \ref{sec:backward-differences}, 
we show that this family of time-truncated solutions is Cauchy, which completes
the proof of Proposition \ref{prop:main-back}. 
Finally, in subsection \ref{sec:backward-residuals} we provide proofs
of the residual bounds necessary for the contraction.

\subsection{Wave asymptotics}\label{sec:wave-asymptotics}
In order to simplify the exposition leading up to the main contraction
argument, we summarize the properties of $ (\psiap, \phiap) $ here,
but relegate the longer proofs to subsection {\ref{sec:backward-residuals}}.
The regularity and $ w  $-decay  assumptions on $ (\rho_{ \infty }^{ + },\rho_{
\infty }^{ - },\Fren) $ from Proposition {\ref{prop:main-back}} were chosen precisely so that
the bounds below hold.
We do not expect these assumptions to be sharp.

The problem of scattering from infinity for free waves has long been studied,
and Lemma \ref{lemma:wave-asymptotics} summarizes the key facts needed for our
analysis. We refer the reader to \cite[section~4]{LindbladSchlue2023} for a more thorough
treatment of this problem using the vector field method, and to \cite[section
6.2]{Hormander} and \cite{Friedlander} for more classical approaches.
\begin{lemma}\label{lemma:wave-asymptotics}
  Fix $0<\nu_0<\frac12$. Let $ \epsilon > 0 $ be sufficiently small and let $ \Fren: \R_{ q } \times \S_{
  \omega }^{ 2 }\to \R $ be a function satisfying
  \begin{align}\label{eq:regasymp}
    \sum_{3(k+\ell) + m\leq27}
    \lVert(q\partial_q)^k\nabla_\omega^\ell \partial_q^m\Fren\rVert_{\dot H_q^{\frac12}L_\omega^2}
    +\sum_{3(k+\ell) + m\leq27}
    \lVert(q\partial_q)^k\nabla_\omega^\ell \partial_q^m\Fren\rVert_{\dot H_q^{\frac12-\nu_0}L_\omega^2}
    &\ll \epsilon.
  \end{align} 
  Then, there exists a solution $ \phif $ to the free wave equation $ \square \phif = 0$
  with asymptotics
  \begin{align}
    \Fren (q,\omega)
    &= \lim_{ r\to\infty } r\phif(r+q,r\omega),
  \end{align} 
  and satisfying the bounds
  \begin{align}
    \sum_{ 3\lvert J \rvert  + K \leq 9 } 
      \lVert x^{ J } \partial_{ x }^{ J+K }\phi_0 \rVert _{ \dot H_{ x }^{
      \frac{1}{2}  }  } 
      + \lVert x^{ J } \partial_{ x }^{ J+K } \phi_1 \rVert _{ \dot H_{ x }^{
      -\frac{1}{2}  }  }
    &\leq\epsilon\label{eq:free-wave-weighted-data},\\
    \lVert \phi_{ 0 } \rVert_{ \dot H_{ x }^{ \frac{1}{2} -\nu_0 }  } 
    + \lVert \phi_1 \rVert_{ \dot H_{ x }^{ -\frac{1}{2} -\nu_0 }  }
    &\leq \epsilon\label{eq:free-wave-low-frequency-data},
  \end{align}
  where $ (\phi_0,\phi_1) := \phif[0] $.
  The parameters $\delta,\nu>0$ above are chosen so that $\delta+\nu<\nu_0/6$
  for compatibility with Proposition \ref{prop:wave-hom-bound}.
\end{lemma}

\begin{proof}
  The inverse scattering transform for the free wave equation is given
  explicitly by 
  \begin{align}\label{eq:inverse-free}
    \phif(t,x)
    &= -\frac{1}{2\pi} \int_{ \S^{ 2 } } \partial_{ q } \Fren (t -
    x\cdot\omega,\omega)\dif\omega,
  \end{align}
  where $ \Fren = \Fren (q,\omega) $. 
  Moreover, given $ s\in \R $, we effectively get an isometry between the
  initial data and the asymptotic radiation:
  \begin{align}\label{eq:radnorm}
    \lVert \phi_0 \rVert_{ \dot H_{ x }^{ s  }  }^2
    + \lVert \phi_1 \rVert_{ \dot H_{ x }^{ s-1  }  }^2
    &= 2\lVert \Fren \rVert _{ \dot H_{ q }^{ s  } L_{ \omega }^{ 2 } }^2,
  \end{align} 
  see \cite{Friedlander}.
  Since $ \square(\Omega\phif) = 0 $, we can study the radiation profile of $
  \Omega\phif $ in a similar fashion.
  Symbolically, 
  \begin{align}
    \lim_{ r\to\infty } r\Omega\phif
    &= a(\omega)q\partial_{ q } \Fren + b(\omega)\nabla_{ \omega }^{ \leq 1 }\Fren
  \end{align} 
  for some functions $ a,b\in C^{ \infty }(\S_{ \omega }^{ 2 }) $.
  Hence, applying \eqref{eq:radnorm} for $ s = \frac{1}{2} $ yields
  \begin{align}
    \lVert x_j\nabla_x f\rVert_{\dot H_x^{\frac{1}{2} }} + \lVert x_j\nabla_x g\rVert_{\dot
    H_x^{-\frac{1}{2} }}
    \lesssim \lVert q\partial_q \Fren\rVert_{\dot H_q^{\frac{1}{2} }L_\omega^2}
    + \lVert\nabla_\omega^{ \leq 1 } \Fren\rVert_{\dot H_q^{\frac{1}{2} }L_\omega^2}.
  \end{align}
  Repeating the argument above for higher-order  derivatives and vector fields
  establishes \eqref{eq:free-wave-weighted-data}.
  The low-frequency bound \eqref{eq:free-wave-low-frequency-data} follows by
  the same argument using $ s = \frac{1}{2} -\nu_0 $ instead.

\end{proof}

Under \eqref{eq:regasymp}, the wave initial data belong to the weighted class
in Theorem \ref{main-thm-gwp} and the linear theory yields the following global bounds.
\begin{coro}
  The free wave constructed in Lemma \ref{lemma:wave-asymptotics} satisfies,
  for $0\leq k\leq18$ and $t\geq1$,
  \begin{align}
    \lVert \Gamma^{ \leq k }\phif \rVert_{ L_{ x }^{ 3 }  } 
    + t\lVert \Gamma^{ \leq k }\phif \rVert_{ L_{ x }^{ \infty }  } 
    \lesssim \epsilon.
  \end{align} 
\end{coro}

The following lemma contains the main set of bounds on the background terms $
\psiap $ and  $ \phiap $ necessary for the nonlinear analysis.

\begin{lemma}\label{lemma:background}
  The background terms satisfy the bounds
  \begingroup
\begin{align}
    &\lVert \Gamma^{ k } \psiap\rVert_{ L^{ 2 }_{ x }  } 
      + t\lVert \Gamma^{ k } \psiap  \rVert_{ L_{ x }^{ 6 }  } 
      \lesssim \epsilon (1 + \log t)^{ k } 
     \label{eq:psiap-bound}
      \\
    &\lVert \Gamma^{ k } \phiap \rVert_{ L_{ x }^{ 3 } }
      + t \lVert \Gamma^{ \leq k }\phiap \rVert_{ L_{ x }^{ \infty }  } 
      \lesssim \epsilon  (1+\log t)^{ k }.
     \label{eq:phiap-bound}
  \end{align}
\endgroup

  Moreover, the residual terms satisfy
  \begingroup
\begin{nalign}\label{eq:residual-bounds}
    \lVert\Gamma^k R_{\psi}(t)\rVert_{L_x^2}
    +\lVert\Gamma^k R_{\phi}(t)\rVert_{\dot H_x^{-1/2}}
    &\lesssim \epsilon^2 t^{-1-\delta-\nu}(1+\log t)^k,
    \quad k\leq18,\\
    \lVert\Gamma^k R_{\phi}(t)\rVert_{\dot H_x^{-1/2-\nu}}
    &\lesssim \epsilon^2 t^{-1-\delta}(1+\log t)^k,
    \quad k\leq9,
  \end{nalign}
\endgroup
  for $t\geq T_0$, with $\delta,\nu>0$ sufficiently small and
  $T_0$ sufficiently large.
\end{lemma}

\begin{remark}\label{rmk:backward-regularity-count}
  We explain the regularity count in Proposition \ref{prop:main-back}.
  Recall that $\Gamma^{\leq k}$ uses the weighted order
  $\lvert I\rvert+3\lvert J\rvert\leq k$ with Lorentz fields having weight
  three.
  The forward energy class has pointwise weighted order $9$ and so, 
  by the Klainerman--Sobolev embedding \ref{prop:klainerman-sobolev}, we need
  to provide $ L^{ 2 } $ control to weighted order $ 18 $.  

  In the Dirac residual, commutation through weighted order $18$ costs at most
  $18$ ordinary velocity derivatives, with one more for the nonresonant
  correction and one for the Dirac operator. This gives $18+2=20$ velocity
  derivatives of $\rho_\infty^\pm$ for \eqref{eq:residual-bounds}.
  For the free wave, we need $L_x^\infty$ control through weighted order $18$.
  Applying Proposition \ref{prop:wave-hom-bound} to these commuted fields costs
  nine additional weighted energy degrees, giving $18+9=27$ in
  \eqref{eq:regasymp}.  

  The weight count comes from having to control terms of the form
  \begin{align}
    \lVert \partial_{ v } (w^{ -\alpha }\rho_{ \infty }) \rVert_{ L_{ v }^{ 2 }
    } \lesssim \lVert \partial_{ v } (w^{ -\alpha }) \rho_{ \infty }^{ \pm }
    \rVert_{ L_{ v }^{ 2 }  } 
    + \lVert w^{ -\alpha }\partial_{ v }\rho_{ \infty }^{ \pm }  \rVert _{ L_{
    v}^{ 2 } }.
  \end{align}
  But $ \lvert\partial_{ v } (w^{ -\alpha })\rvert \lesssim w^{ -\alpha-2 } $;
  this factor of $ w^{ -2 } $ is behind the $ w^{ 2\lvert \beta \rvert } $
  term in Proposition \ref{prop:main-back}.
  This means that the worst decaying terms occur by repeatedly differentiating
  $ w^{ -\alpha } $, and terms with higher-order derivatives on $ \rho_{ \infty
  } $ need less boundary decay.
  Allowing a loss of $w^{-2}$ for each of the $20$ derivatives, together with
  a safety margin, gives the non-optimized weight $w^{-50+2\lvert\beta\rvert}$
  in the profile assumption.
\end{remark}

\subsection{Contraction argument and proof of Proposition
\ref{prop:T-bounds}}\label{sec:backward-contraction}

Setting 
\begin{align}
  \Psi = \psi - \psiap,
  &\qquad
  \Phi = \phi - \phiap
\end{align} 
we get the following equations for the difference
\begin{align}\label{eq:difference}
  (-i\gamma^{ \mu }\partial_{ \mu } + 1)\Psi
  &= \Phi\Psi + \Phi \psiap + \phiap \Psi - R_{ \psi }\\
  \square \Phi
  &= \langle \Psi, \gamma^{ 0 }\Psi \rangle + 2\Re\, \langle \Psi, \gamma^{ 0
    }\psiap \rangle - R_{ \phi }
\end{align} 
with residuals
\begin{align}
  R_{ \psi } &:= (-i\gamma^{ \mu }\partial_{ \mu } + 1)\psiap -  \phiap \psiap \\
  R_{ \phi } &:= \square \phiap - \langle \psiap, \gamma^{ 0 }\psiap \rangle  
\end{align} 
independent of $ (\Psi,\Phi) $.

In reality, instead of solving \eqref{eq:difference} with data prescribed at $
t = \infty$ directly, we first solve a family of intermediate equations with
final-state data prescribed at $ t = T$ 
\begin{nalign}\label{eq:differenceT}
  (-i\gamma^{ \mu }\partial_{ \mu } + 1)\Psi_{ T } 
  &= 
    \vartheta(t / T)
    \big(
      \Phi_{ T } \Psi_{ T }  + \Phi_{ T }  \psiap + \phiap \Psi_{ T }
      - R_{ \psi }
    \big)\\
  \square \Phi_{ T } 
  &= 
    \vartheta(t / T)
    \big(
      \langle \Psi_{ T } , \gamma^{ 0 }\Psi_{ T }  \rangle + 2\Re\, \langle
      \Psi_{ T } , \gamma^{ 0 }\psiap \rangle - R_{ \phi }
    \big) \\
  &\Psi_{ T } \,\big|_{ t = T }  = 0,
  \quad (\Phi_{ T } , \partial_{ t } \Phi_{ T } )\,\big|_{ t = T } = (0,0).
\end{nalign}
Later on, comparing the solutions $ (\Psi_{ T } ,\Phi_{ T } )_{ T \geq T_0 }  $
along  a fixed initial hypersurface $ \left\{ t = T_0 \right\}  $ will then
produce a solution to \eqref{eq:difference}.
The addition of the cutoff $ \vartheta $ is convenient for technical reasons
but not strictly necessary. In particular, the cutoff lets us set $\Psi_{ T }
(t) = \Phi_{ T } (t)=  0,$ for all $ t \geq T $.

Throughout this section, $\delta>0$ denotes a sufficiently small decay exponent
whose value may decrease from one estimate to the next. Only finitely many such
reductions are needed, and the final estimates hold with a common positive
exponent. We now state the main proposition of this section:
\begin{prop}\label{prop:T-bounds}
  Let $ T \geq T_0 > 0 $ be fixed.
  The equation \eqref{eq:differenceT} admits solutions
  $(\Psi_{ T }, \Phi_{ T } ) $ satisfying 
  \begingroup
\begin{nalign}\label{eq:T-solution-bounds}
    \lVert \Gamma^k(\Psi_T,\Phi_T)(t)\rVert_{\Hcal^0}
    &\lesssim \epsilon^2 t^{-\delta}(1+\log t)^k,
    && k\leq18,\\
    \lVert\Gamma^k\Psi_T(t)\rVert_{L_x^6}
    &\lesssim \epsilon^2 t^{-1-\delta}(1+\log t)^{k+9},
    && k\leq9,\\
    \lVert\Gamma^k\Phi_T(t)\rVert_{L_x^\infty}
    &\lesssim \epsilon^2 t^{-1-\delta}(1+\log t)^{k+9},
    && k\leq9,
  \end{nalign}
\endgroup
  for some fixed small $ \delta > 0 $, and uniformly in  $ t\in
  [T_0,T]$.
  Moreover, $ \Phi_{ T }  $ satisfies the additional low-frequency bounds
  \begingroup
\begin{align}\label{eq:T-low-frequency-bounds}
      \lVert\Gamma^k\Phi_T(t)\rVert_{\dot H_x^{\frac12-\nu}}
      +\lVert\partial_t\Gamma^k\Phi_T(t)\rVert_{\dot H_x^{-\frac12-\nu}}
    &\lesssim \epsilon^2 t^{-\delta}(1+\log t)^k,
     \quad k\leq 9,
  \end{align}
\endgroup
  for some small $ \nu > 0 $.
\end{prop}

\begin{proof}
We prove Proposition \ref{prop:T-bounds} via a contraction argument.
Consider the solution norms
\begin{align}
    \lVert (\Psi,\Phi)\rVert_{S_T}
  &:= \max_{0\leq k\leq18}\sup_{t\in[T_0,T]}
    \frac{t^{\delta+\nu}}{(1+\log t)^k}
    \lVert\Gamma^k(\Psi,\Phi)(t)\rVert_{\Hcal^0}\\
  &\quad+\max_{0\leq k\leq9}\sup_{t\in[T_0,T]}
    \frac{t^\delta}{(1+\log t)^k}
    \left(
      \lVert\Gamma^k\Phi(t)\rVert_{\dot H_x^{\frac12-\nu}}
      +\lVert\partial_t\Gamma^k\Phi(t)\rVert_{\dot H_x^{-\frac12-\nu}}
    \right)\\
  &\quad+\max_{0\leq k\leq9}\sup_{t\in[T_0,T]}
    \frac{t^{1+\delta}}{(1+\log t)^{k+9}}
    \lVert\Gamma^k\Psi(t)\rVert_{L_x^6}.
\end{align}
Additionally, by a Sobolev embedding and 
\cite[Remark 7.2]{HIS},
we have control over the norms
\begin{align}
    \lVert\Gamma^k\Phi(t)\rVert_{L_x^3}
  &\lesssim \lVert\Gamma^k\Phi(t)\rVert_{\dot H_x^{\frac12}}
    \lesssim t^{-\delta-\nu}(1+\log t)^k
    \lVert\Phi\rVert_{S_T},\quad k\leq18,\\
    \lVert\Gamma^k\Phi(t)\rVert_{L_x^\infty}
  &\lesssim t^{-1-\delta}(1+\log t)^{k+9}
    \lVert\Phi\rVert_{S_T},\quad k\leq9.
\end{align}
By the linear estimates {\eqref{eq:linear-energy-estimates}} and
the Klainerman-Sobolev embedding (Proposition {
\ref{prop:klainerman-sobolev}})
\begin{nalign}\label{eq:linear-sufficient}
  \lVert(\Psi,\Phi)\rVert_{S_T}
  &\lesssim
    \max_{k\le18}\sup_{t\in[T_0,T]}
    \frac{t^{1+\delta+\nu}}{(1+\log t)^k}
    \left(
    \lVert\Gamma^kF(t)\rVert_{L_x^2}
    +\lVert\Gamma^kG(t)\rVert_{\dot H_x^{-1/2}}
    \right)\\
  &\qquad+
    \max_{k\le9}\sup_{t\in[T_0,T]}
    \frac{t^{1+\delta}}{(1+\log t)^k}
    \lVert\Gamma^kG(t)\rVert_{\dot H_x^{-1/2-\nu}}\\
  &\qquad+
    \max_{k\le9}\sup_{t\in[T_0,T]}
    \frac{ t^{\delta} }{ (1 + \log t)^{k+9} }
    t^{1/2}\lVert F\rVert_{L^2_{t,x}([t/2, 2t]\times \R^3)}
\end{nalign}
where $ F = (-i\gamma^{ \mu }\partial_{ \mu }+1) \Psi  $, $ G = \square\Phi $,
and with implicit constant independent of $ T $.
Thus, in order to establish a contraction, it is sufficient to show that the
right-hand side of \eqref{eq:linear-sufficient} is controlled by 
$ C\epsilon ^2 + C\epsilon\lVert (\Psi,\Phi) \rVert_{ S_{ T }  }  + C\lVert
(\Psi,\Phi) \rVert _{ S_{ T }  } ^2 $.
Control over the $ L^{ 2 }_{ t,x }  $ norm is immediate from the other terms in
\eqref{eq:linear-sufficient}, since $ \vartheta(t /T)  $ is supported on a time
interval of length $ O(T) $, so we only need to provide estimates for the
energy and low-frequency $ \dot H^{ \frac{1}{2} -\nu } $ norms.

We start by estimating the Dirac nonlinearities.
Let $k\leq18$. 
By Hölder's inequality
\begin{align}
  \lVert\Gamma^k(\phiap\Psi)\rVert_{L_x^2}
  &\lesssim\sum_{j+\ell\leq k}
    \lVert\Gamma^j\phiap\rVert_{L_x^\infty}
    \lVert\Gamma^\ell\Psi\rVert_{L_x^2}
    \lesssim\epsilon t^{-1-\delta-\nu}(1+\log t)^k
    \lVert\Psi\rVert_{S_T},\\
  \lVert\Gamma^k(\Phi\psiap)\rVert_{L_x^2}
  &\lesssim\sum_{j+\ell\leq k}
    \lVert\Gamma^j\Phi\rVert_{L_x^3}
    \lVert\Gamma^\ell\psiap\rVert_{L_x^6}
    \lesssim\epsilon t^{-1-\delta-\nu}(1+\log t)^k
    \lVert\Phi\rVert_{S_T}.
\end{align}
For $ \Gamma^{ k } (\Phi\Psi) $, since there can only be at most one derivative
of order $ > 9 $ falling on either $ \Phi $ or $ \Psi $, we estimate
\begin{align}
  \lVert\Gamma^k(\Phi\Psi)\rVert_{L_x^2}
  &\lesssim\sum_{\substack{j+\ell\leq k\\ \ell > 9}}
  \lVert\Gamma^j\Phi\rVert_{L_x^\infty}
  \lVert\Gamma^\ell\Psi\rVert_{L_x^2}
  +\sum_{\substack{j+\ell\leq k\\ \ell \leq 9}}
  \lVert\Gamma^j\Phi\rVert_{L_x^3}
  \lVert\Gamma^\ell\Psi\rVert_{L_x^6}\\
  &\lesssim t^{-1-2\delta-\nu}(1+\log t)^{k+9}
  \lVert\Phi\rVert_{S_T}\lVert\Psi\rVert_{S_T}.
\end{align} 
For the wave nonlinearities, we estimate
\begin{align}
    \lVert\Gamma^k\langle\psiap,\gamma^0\Psi\rangle
  \rVert_{\dot H_x^{-\frac12}}
  &\lesssim\sum_{j+\ell\leq k}
    \lVert\Gamma^j\psiap\rVert_{L_x^6}
    \lVert\Gamma^\ell\Psi\rVert_{L_x^2}
    \lesssim\epsilon t^{-1-\delta-\nu}(1+\log t)^k
    \lVert\Psi\rVert_{S_T},\\
    \lVert\Gamma^k\langle\Psi,\gamma^0\Psi\rangle
    \rVert_{\dot H_x^{-\frac12}}
  &\lesssim\sum_{\substack{j+\ell\leq k\\j\leq\ell}}
    \lVert\Gamma^j\Psi\rVert_{L_x^6}
    \lVert\Gamma^\ell\Psi\rVert_{L_x^2}
    \lesssim t^{-1-2\delta-\nu}(1+\log t)^{k+9}
    \lVert\Psi\rVert_{S_T}^2.
\end{align}
Thus, for $k\leq9$,
\begin{align}
  \lVert\Gamma^k\langle\psiap,\gamma^0\Psi\rangle
  \rVert_{\dot H_x^{-\frac12-\nu}}
  &\lesssim\sum_{j+\ell\leq k}
    \lVert\Gamma^j\psiap\rVert_{L_x^{6/(1+2\nu)}}
    \lVert\Gamma^\ell\Psi\rVert_{L_x^2}
    \lesssim\epsilon t^{-1-\delta}(1+\log t)^k
    \lVert\Psi\rVert_{S_T},\\
    \lVert\Gamma^k\langle\Psi,\gamma^0\Psi\rangle
    \rVert_{\dot H_x^{-\frac12-\nu}}
  &\lesssim\sum_{j+\ell\leq k}
    \lVert\Gamma^j\Psi\rVert_{L_x^{6/(1+2\nu)}}
    \lVert\Gamma^\ell\Psi\rVert_{L_x^2}
    \lesssim t^{-1-2\delta-\nu^2}(1+\log t)^{k+9}
    \lVert\Psi\rVert_{S_T}^2,
\end{align}
where we have used the following interpolated bounds
\begin{align}
  \lVert\Gamma^j\psiap\rVert_{L_x^{6/(1+2\nu)}}
  &\lesssim\lVert\Gamma^j\psiap\rVert_{L_x^2}^{\nu}
  \lVert\Gamma^j\psiap\rVert_{L_x^6}^{1-\nu}
  \lesssim\epsilon t^{-1+\nu}(1+\log t)^j,\quad j\leq18,\\
  \lVert\Gamma^j\Psi\rVert_{L_x^{6/(1+2\nu)}}
  &\lesssim\lVert\Gamma^j\Psi\rVert_{L_x^2}^{\nu}
  \lVert\Gamma^j\Psi\rVert_{L_x^6}^{1-\nu}
  \lesssim t^{-1-\delta+\nu-\nu^2}(1+\log t)^{j+9(1-\nu)}
  \lVert\Psi\rVert_{S_T},\quad j\leq9.
\end{align}
A similar argument furnishes the wave estimates in $L_x^{3/2}$ and
$L_x^{3/(2+\nu)}$.
We should point out that the estimates above were proved for the nonlinearities
without the extra time cutoff.
The addition of $ \vartheta (t /T) $, however, is harmless; all nonlinear interactions
containing either $ \psiap $ or $ \phiap $ are supported inside the truncated
cone $ \left\{ \lvert t-r \rvert  \leq T, t \geq T_0 \right\}  $ and, since
\begin{align}
  [\partial_{ t } , \vartheta (t /T)] = O(T^{ -1 })
  \quad\text{and}\quad 
  [\Omega, \vartheta (t /T)] = O(1)
\end{align} 
in this region, the same bounds as above follow. 
The estimates for the $ \vartheta(t /T) \Phi\Psi $ and $ \vartheta (t /T)
\langle \Psi, \gamma^{ 0 }\Psi \rangle $ nonlinearities similarly follow from
the fast decay enforced by the $ S_{ T }  $-norm.

\end{proof}

\subsection{The difference equations}\label{sec:backward-differences}

In this section, we  show that the family $ (\Psi_{ T } ,\Phi_{ T } )_{ T\geq
T_0 } $ is Cauchy in $ T $ and, thus, must converge to a solution to
equation {\eqref{eq:untrunc}} with zero data at infinity. 
This is enough to complete the proof of Proposition \ref{prop:main-back}.

Given  $ T_2 > T_1 \geq T_0 $, define the differences
\begin{align}
  \delta\Psi := \Psi_{ T_2 } - \Psi_{ T_1 } 
  \quad\text{and}\quad 
  \delta\Phi := \Phi_{ T_2 } - \Phi_{ T_1 } .
\end{align} 
We show that 
\begin{align}
  \sup_{ t\in [T_0,\infty) } 
  \lVert \Gamma^{ \leq 9 }(\delta\Psi(t),\delta\Phi[t]) \rVert_{ \Hcal^{ 0 } }
  \lesssim \epsilon T_{ 1 }^{ -\delta }. 
\end{align} 

First, we note that $ \delta\Psi (t) = \delta\Phi(t) = 0$ for $ t\geq T_2 $.
Moreover, the a priori bounds {\eqref{eq:T-solution-bounds}} immediately yield
\begin{align}\label{eq:good-T}
  \lVert \Gamma^{ \leq 9 }(\delta\Psi,\delta\Phi) (t) \rVert_{ \Hcal^{ 0 } }
  \lesssim \lVert \Gamma^{ \leq 9 }(\Psi_{ T_2 },\Phi_{ T_2 } )  (t)\rVert_{
    \Hcal^{ 0} } 
    + \lVert \Gamma^{ \leq 9 }(\Psi_{ T_1 },\Phi_{ T_1 } ) (t)\rVert_{ \Hcal^{ 0}}
  \lesssim \epsilon T_{ 1 }^{ -\delta } 
\end{align} 
uniformly for $ t\in [\frac{T_1}{2} , T_2] $. 

In the region $ [T_0,\frac{T_1}{2} ] $, we instead consider the difference
equations 
\begin{align}
  (-i\gamma^{ \mu }\partial_{ \mu } + 1)\delta\Psi
  &= \vartheta(t/T_2)\mathcal N_{ 1 }^{ D } +(\vartheta(t/T_2)-\vartheta(t/T_1))
    \mathcal N_{ 2 }^{ D }\\
      \square \delta\Phi
  &= \vartheta(t/T_2)\mathcal N_{ 1 }^{ W }
    + (\vartheta(t/T_2)-\vartheta(t/T_1)) \mathcal N_{ 2 }^{ W },
\end{align} 
where
\begin{align}
  \mathcal N_{ 1 }^{ D }
  &:= (\Phi_{ T_2 }+\phiap)\delta\Psi
      +(\Psi_{ T_1 }+\psiap)\delta\Phi\\
  \mathcal N_{ 2 }^{ D } 
  &:= \Phi_{ T_1 }\Psi_{ T_1 }+\Phi_{ T_1 }\psiap
      +\phiap\Psi_{ T_1 }-R_{ \psi },\\
  \mathcal N_{ 1 }^{ W } 
  &:= 
    \Re\langle \delta\Psi , \gamma^{ 0 }\Psi_{ T_1 }  \rangle
    + \Re\langle \delta\Psi, \gamma^{ 0 }\Psi_{ T_2 }  \rangle
    + 2\Re \langle \delta\Psi, \gamma^{ 0 }\psiap \rangle,\\
  \mathcal N_{ 2 }^{ W }
  &:=  \langle \Psi_{ T_1 },\gamma^{ 0
    }\Psi_{ T_1 }\rangle +2\Re\,\langle \Psi_{ T_1 },\gamma^{ 0 }\psiap\rangle
      -R_{ \phi }.
\end{align} 
In fact, since  $ \vartheta(t /T_2) - \vartheta(t /T_1) = 0 $ for $ t\in [T_0,
T_1 /2]$, it suffices to consider the simpler equations
\begin{align}\label{eq:differencesT}
    (-i\gamma^{ \mu }\partial_{ \mu } + 1)\delta\Psi
  &= \mathcal N_{ 1 }^{ D } 
  \\
   \square \delta\Phi
  &= \mathcal N_{ 1 }^{ W }
\end{align} 
which we estimate directly.
Placing scalars in $L^\infty$, known spinors in $L^6$, and differences in $L^2$
and $\dot H^{1/2}$, we get
\begin{align}
  &\lVert \Gamma^{k}\big((\Phi_{T_2}+\phiap)\delta\Psi\big) \rVert_{L^2}
  \lesssim \sum_{j=0}^k
    \lVert \Gamma^{k-j}(\Phi_{T_2}+\phiap)\rVert_{L^\infty}
    \lVert \Gamma^{j}\delta\Psi\rVert_{L^2},\\
  &\lVert \Gamma^{k}\big((\Psi_{T_1}+\psiap)\delta\Phi\big) \rVert_{L^2}
  \lesssim \sum_{j=0}^k
    \lVert \Gamma^{k-j}(\Psi_{T_1}+\psiap)\rVert_{L^6}
    \lVert \Gamma^{j}\delta\Phi\rVert_{\dot H^{1/2}},
\end{align}
and
\begin{align}
  &\lVert \Gamma^{k}\langle\delta\Psi,\gamma^0\Psi_{T_1}\rangle
    \rVert_{\dot H^{-1/2}}
    +\lVert \Gamma^{k}\langle\delta\Psi,\gamma^0\psiap\rangle
    \rVert_{\dot H^{-1/2}}\\
  &\hspace{14em}\lesssim \sum_{j=0}^k
    \left(
      \lVert \Gamma^{k-j}\Psi_{T_1}\rVert_{L^6}
      +\lVert \Gamma^{k-j}\psiap\rVert_{L^6}
    \right)
    \lVert \Gamma^{j}\delta\Psi\rVert_{L^2},
\end{align} 
for $ k\leq 9 $.
The remaining term $ \langle \delta\Psi, \gamma^{ 0 }\Psi_{ T_2 }   \rangle $
can be controlled analogously.
In particular, by {Lemma \ref{lemma:background} and Proposition \ref{prop:T-bounds}},  
\begin{align}\label{eq:N1-bounds}
  \lVert \Gamma^{ \leq k }\mathcal N_{ 1 }^{ D } \rVert_{  L_{ x }^{ 2 }  }
  + \lVert \Gamma^{ \leq k }\mathcal N_{ 1 }^{ W } \rVert_{  \dot H_{ x }^{
  -\frac{1}{2}  }  }
  &\leq C\epsilon t^{ -1 }
  \sum_{ j=0 }^{ k } (1+\log t)^{ k-j } 
  \lVert \Gamma^{  j }(\delta \Psi,\delta\Phi)(t) \rVert_{ \Hcal^{ 0 }  } .
\end{align} 
Now we want to apply Grönwall's inequality backward in time to transfer the
good bounds from $ t = \frac{T_1}{2}  $ to $ t = T_0$.
To simplify the notation, denote
\begin{align}
   E_{ k }(t)  := \lVert \Gamma^{ k } (\delta\Psi,\delta\Phi) (t)\rVert_{ \Hcal^{
   0} }.
\end{align} 
At the ground level $ k = 0 $, the borderline non-integrable decay rate 
$ t^{ -1 }  $
is harmless, since it gets absorbed by the factor $ t^{ -1-\delta } $ from
Proposition \ref{prop:T-bounds}.
Indeed, by the standard energy estimate {\eqref{eq:linear-energy-estimates}}, Grönwall's inequality,
Proposition \ref{prop:T-bounds}, and \eqref{eq:good-T}
\begin{align}
  E_{ 0 }(t)
  \leq
    E_{ 0 }
    \left( 
      T_1 /2 
    \right) 
    \left( 
      \frac{T_1}{t}  
    \right)^{ C\epsilon }
  \lesssim \epsilon T_{ 1 }^{ -\delta }
\end{align} 
uniformly in $ \left[ T_0,\frac{T_1}{2}  \right]  $.
For the higher-order energies, repeating the previous argument yields
\begin{align}
  E_{ k } (t)
  &\leq E_{ k } (T_1 /2) 
  \left( 
    \frac{T_1}{t}  
  \right) ^{ C\epsilon } 
  + C\epsilon \int_{ t }^{ T_1 /2 } \sum_{ j=0 }^{ k-1 } 
  (1+\log s)^{ k-j } s^{ -1 }\left( 
    \frac{s}{t}  
  \right)^{ C\epsilon }   E_{ j } (s) \dif s.
\end{align} 
Assuming inductively that $ \sup_{ s\in [T_0,T_1 /2] } E_{ j } (s) \lesssim
\epsilon T_{ 1 }^{ -\delta }$ for all $ j < k $, and 
integrating the expression above, we conclude that
\begin{align}
  E_{ k } (t)
  &\lesssim \epsilon T_{ 1 }^{ -\delta },
\end{align} 
uniformly for $ t\in [T_0, \frac{T_1}{2} ] $, as desired.

While the argument above shows the existence of $ (\Psi,\Phi) := \lim_{
T\to\infty } (\Psi_{ T } ,\Phi_{ T } ) $ with uniform trace bounds on  $
\lVert \Gamma^{ \leq 9 }(\Psi(T_0), \Phi[T_0]) \rVert_{ \mathcal H^{ 0 } } $,
it does not, by itself, imply the low-frequency bound assumed in the forward
Theorem {\ref{main-thm-gwp}}. 
Instead of repeating the difference estimates for the low-frequency norms, it
is sufficient to note that the proof of Proposition {\ref{prop:T-bounds}} establishes a priori
bounds for the untruncated equations, and since the limit above is taken in a
strong $ C_{ t }^{ 0 } \mathcal H^{ 0 } $ sense, $ \Phi $ directly
satisfies {\eqref{eq:T-low-frequency-bounds}}.

\subsection{Residual bounds and proof of Lemma \ref{lemma:background}}\label{sec:backward-residuals}

Here we prove the background and residual bounds needed for the contraction
argument used in Proposition \ref{prop:T-bounds}.

\begin{proof}[Proof of Lemma \ref{lemma:background}]
  We use coordinates $ (t,x) = (t,vt) $ and the adapted vector fields
  \begin{align} 
    \breve\partial_t & := \partial_t + \frac{1}{t} (x\cdot\nabla_{ x } ),
    \qquad
    \breve\Omega_{0j} :=\partial_{v_j}-v_j(v\cdot\nabla_v),
  \end{align} 
  and start by establishing the background bounds {\eqref{eq:psiap-bound}} and {\eqref{eq:phiap-bound}}.
  
  We directly compute
  \begin{align}
    \Omega_{0j}\psiap^\pm
    &= \left( 
      tv_{ j } \breve \partial_{ t } + \breve \Omega_{ 0j } 
    \right) \psiap^{ \pm }\\
    &=(tw)^{-3/2}e^{\pm itw\mp i\Theta}
    \left[
      tv_j\breve{\partial}_t\varrho^\pm
      +\breve\Omega_{0j}\varrho^\pm
      \mp i
      \left(
        tv_j\breve{\partial}_t\Theta
        +\breve\Omega_{0j}\Theta
    \right)
    \varrho^\pm
    \right].
  \end{align} 
  Similarly,
  \begin{align}
    \partial_{ x_{ j } } \psiap^{ \pm }
    &= (tw)^{-3/2}e^{\pm itw\mp i\Theta}
      \left[
        \frac1t\partial_{v_j}\varrho^\pm
        +\left(
          \frac{3v_j}{2tw^2}
          \mp i\left(\frac{v_j}{w}+\frac1t\partial_{v_j}\Theta\right)
        \right)\varrho^\pm
      \right],\\
    \partial_{ t } \psiap^{ \pm }
    &= (tw)^{-3/2}e^{\pm itw\mp i\Theta}
      \left[
        \breve{\partial}_t\varrho^\pm
        -\frac1t(v\cdot\nabla_v)\varrho^\pm
      \right.\\
    &\qquad\left.
        +\left(
          -\frac{3}{2tw^2}
          \pm\frac{i}{w}
          \mp i\left(
            \breve{\partial}_t\Theta-\frac1t(v\cdot\nabla_v)\Theta
          \right)
        \right)\varrho^\pm
      \right].
  \end{align}

  Using that $ \lVert \varrho^{ \pm } \rVert_{ L_{ x }^{ p }  } = t^{ \frac{3}{p}
  }\lVert \varrho^{ \pm } \rVert_{ L_{ v }^{ p }  }  $ and iterating the
  computation above, we get
  \begin{align}
    \lVert \widehat\Omega^{ \leq k }\psiap^{ \pm } \rVert_{ L_{ x }^{ 2 } }
    &\lesssim (1+\log t)^k
      \sum_{j+|J|\leq k}
      \left\lVert w^{-\frac{3}{2}}(t\breve\partial_t)^j
      \breve\Omega^J\varrho^\pm \right\rVert_{L_v^2},\\
    \lVert \widehat\Omega^{ \leq k }\psiap^{ \pm } \rVert_{ L_{ x }^{ 6 } }
    &\lesssim t^{-1}(1+\log t)^k
      \sum_{j+|J|\leq k}
      \left\lVert w^{-\frac{3}{2}}(t\breve\partial_t)^j
      \breve\Omega^J\varrho^\pm \right\rVert_{L_v^6}.
  \end{align} 

  Notice that the bounds for the translation fields $ \partial_{ \mu } $ are
  worse than the corresponding bounds for $ \Omega $; this follows from the
  Lorentz invariance of $ tw = \Lambda $:
  \begin{align}
    \partial_{ x_j }(tw)^{ -a } \sim x_j(tw)^{ -a - 2}
    \quad\text{whereas}\quad 
    \Omega (tw)^{ -a } =   0,
    \quad  a > 0.
  \end{align} 

  The bounds on the free wave component $ \phif $ follow from Lemma {\ref{lemma:wave-asymptotics}}, and
  so it is enough to establish
  {\eqref{eq:phiap-bound}} for $ \phit $.
  Denote
  \begin{align}
    \Sm := \1_{ \left\{ \lvert v \rvert  < 1 \right\}  }(tw)^{ -3 } (\lvert
    \rho_{ \infty }^{ - }\rvert _{ H }^{ 2 } - \lvert \rho_{ \infty }^{ + }
    \rvert _{ H }^{ 2 }  )
  \end{align}
  so that $ \square\phit = \Sm $.
  By a direct computation,
  \begin{align}
    \Omega_{ 0j } \Sm
    &= \1_{ \left\{ \lvert v \rvert < 1 \right\}  } (tw)^{ -3 }
      \big(
        tv_{ j } \breve \partial_{ t } b_{\infty} + \breve \Omega_{ 0j }b_{\infty}
      \big),\\
      \partial_{x_j}S_{\mathrm{main}}
    &=\mathbf1_{\{|v|<1\}} t^{ -3 } w^{ -4 }
      \left[
        \frac{w^2}{tw}\partial_{v_j}b_{\infty}
        +\frac{3v_j}{tw}b_{\infty}
      \right].
  \end{align} 
  Iterating the computation above yields
  \begin{align}\label{eq:pointwise-Sm}
    \lvert \Gamma^{ \leq k } \Sm \rvert 
    \lesssim t^{ -3 }
    w^{ -2k-3 }
    \sum_{ j + |J| \leq k } 
    \lvert (t\breve \partial_{ t } )^{ j } \partial_{ v }^{ J } b_{\infty} \rvert.
  \end{align} 
  We recall Kirchhoff's formula in coordinates $ (t,v) $ with zero initial data,
  \begin{align}\label{eq:kirchhoff}
    \phit (t,v)
    &= \frac{1}{4\pi} 
      \int_{ 1 }^{ t } (t-s) \int_{ \S^{ 2 } } \Sm \left(
        s, \frac{ tv + (t-s)\omega }{ s } 
      \right) 
      \dif\omega \dif s.
  \end{align} 
  We note that 
  \begin{align}\label{eq:timelike-wave-agreement}
    \phit(t,v) = (tw)^{ -1 } \phi_{ \infty }(v)
  \end{align}
  on $ \left\{ t-r\geq 2 \right\}  $, 
  as their Kirchhoff representations differ only by integration over
  \(0<s<1\), whose contribution vanishes in this region.

  The iterates $ \Gamma^{ j } \phit $ solve the wave equation
  \begin{align}
    \square\Gamma^{ j } \phit 
    &= \Gamma^{ j } \Sm
  \end{align} 
  with constrained $ \Gamma^{ j } \phit [0]$: higher-order derivatives must be
  compatible with the lower-order equations and the ground-state initial data $
  (\phi_0,\phi_1) $. 
  Since such states contribute homogeneous terms of size
  $O(\epsilon^2t^{-1})\1_{\{|v|<1\}}$, by
  \eqref{eq:kirchhoff} and \eqref{eq:pointwise-Sm},
  \begin{align}
    \lvert \Gamma^{ \leq k }\phit \rvert 
    &\lesssim \frac{\epsilon^2}{t}\1_{\{|v|<1\}}
      +\int_{ 1 }^{ t } (t-s) \int_{ \S^{ 2 } }
      \left\lvert 
        (\Gamma^{ \leq k }\Sm)
        \left( 
          s, \frac{ tv + (t-s)\omega }{ s } 
        \right)  \right\rvert 
      \dif\omega \dif s\\
    &\lesssim \frac{1}{t} \1_{ \left\{ \lvert v \rvert  < 1 \right\}  }
      \left(\epsilon^2+
      \sum_{ j + \lvert J \rvert  \leq k } 
      \left\lVert w^{ -2k-3 } (s\breve \partial_{ s } )^{ j } \partial_{ v }^{ J }
      b_{\infty} \right\rVert_{ L^{ \infty }_{ s }
        \left( 
          [1,t]
      \right)   
      L_{ v }^{ \infty } 
      }
      \right).
  \end{align} 
  
  In summary,  $ \lVert \Gamma^{ \leq k }\phit \rVert_{ L^{ \infty }_{ x }  }
  \lesssim \epsilon^2 t^{ -1 } $,
  and integrating this on $ \supp \phit \subset\left\{ t \geq r \right\}
  $, we conclude that
  \begin{align}
    \lVert \Gamma^{ \leq k }\phit \rVert_{ L_{ x }^{ 3 }  } 
    \lesssim \epsilon^2.
  \end{align} 

  Now we turn our attention to the residuals $ \Gamma^{ k } R_{ \psi } $ and $
  \Gamma^{ k } R_{ \phi } $.

  \vspace{5pt}
  \noindent
  \textit{The $ k = 0 $ case.}
  We write the Dirac operator $ -i\gamma^{ \mu 
  }\partial_\mu $ in terms of $ (\breve \partial_{ t }, \partial_{ v_j } ) $ 
  \begin{align}
    -i\gamma^{ \mu }\partial_{ \mu } 
    &= -i\gamma^{ 0 }\breve \partial_{ t } - \frac{i}{t} (\gamma^{ j } -v_{ j }
    \gamma^{ 0 })\partial_{ v_{ j }  }.
  \end{align} 
  Using the following identities
  \begin{align}
    &\gamma^\mu \partial_\mu(tw)
    = \frac{\gamma^0 - v_j\gamma^j}{w}
    = -\gamma^H,\\
    &\gamma^\mu \partial_\mu(tw)^{-3/2}
    = \frac{3}{2tw}\gamma^H(tw)^{-3/2},\\
    &\breve{\partial}_t\Theta
    = \frac{\phi_\infty}{t},\\
    &\partial_{v_j}\Theta
    = \log(tw)\partial_{v_j}\phi_\infty
    - \frac{v_j}{w^2}\phi_\infty,\\
    &\gamma^\mu\partial_\mu\Theta
    = -\frac{\phi_\infty}{tw}\gamma^H
    + \frac{\log(tw)}{t}
    (\gamma^j-v_j\gamma^0)\partial_{v_j}\phi_\infty,
  \end{align}
  we directly compute
  \begin{nalign}\label{eq:dirac-psiap}
    (-i\gamma^\mu\partial_\mu + 1)\psiap^{\pm} - \phiap \psiap^{\pm}
    &= \frac{e^{\pm itw\mp i\Theta}}{(tw)^{3/2}}
    (1\mp\gamma^H)\varrho^\pm
    -\frac{i e^{\pm itw\mp i\Theta}}{(tw)^{3/2}}
    \gamma^0\breve{\partial}_t\varrho^\pm\\
    &\quad+\frac{e^{\pm itw\mp i\Theta}}{t(tw)^{3/2}}
    \left[
    -i(\gamma^j-v_j\gamma^0)\partial_{v_j}\varrho^\pm
    -\frac{3i}{2w}\gamma^H\varrho^\pm
    \pm\frac{\phi_\infty}{w}\gamma^H\varrho^\pm
    \right]\\
    &\quad\mp\frac{\log(tw)e^{\pm itw\mp i\Theta}}{t(tw)^{3/2}}
    (\gamma^j-v_j\gamma^0)(\partial_{v_j}\phi_\infty)\varrho^\pm
    - \phiap \psiap^{\pm}\\
  &=: R_{ 1 }^{ \pm  } + R_{ 2 }^{ \pm }+ R_{ 3 }^{ \pm }+ R_{ 4 }^{ \pm }
  - \phiap\psiap^{\pm}.
  \end{nalign}
  Up to this point, we have not used the explicit definition of the regularized
  profiles $ \varrho^{ \pm } $ in {\eqref{eq:corrected-profiles}}.
  In order to properly motivate this definition, we first consider the more
  natural first choice of profiles $ \varrho^{ \pm }(t,v) = \rho_{ \infty }^{
  \pm }(v) $.
  In this case,  $ R_1^{ \pm } = 0 $, since $ \rho_{ \infty }^{ \pm }\in V^{
  \pm } $; however, 
  \begin{align}
    \lVert R_{ 3 }^{ \pm } \rVert _{ L_{ x }^{ 2 }  } 
    &\lesssim  t^{ -1 }\left( 
      \lVert w^{ -\frac{5}{2}  }(1+\lvert\phi_\infty\rvert) \rho_{ \infty }^{ \pm } \rVert_{ L_{ v }^{ 2 }  }
      + \lVert w^{ -\frac{3}{2}  }\partial_v\rho_{ \infty }^{ \pm }  \rVert_{ L_{ v }^{ 2 }  }
    \right), 
  \end{align} 
  which does not satisfy {\eqref{eq:residual-bounds}} and cannot be improved by assuming additional
  $ w $-decay on $ \rho_{ \infty }^{ \pm } $. 
  A similar computation reveals the same issue with  $ R_{ 4 }^{ \pm }
  $, with the additional loss
  \begin{align}
    \lVert R_{ 4 }^{ \pm } \rVert _{ L_{ x }^{ 2 }  } 
    &\lesssim t^{ -1 }(1+\log t) \lVert w^{ -\frac{5}{2}  } \partial_v\phi_{
    \infty }  \rho_{ \infty }^{ \pm } \rVert _{ L_{ v }^{ 2 }  }.
  \end{align} 

  As in the construction of the corrected wave packets for the Dirac
  equation {\eqref{eq:wp-correction}}, the more slowly decaying terms can be addressed by adding an
  appropriate nonresonant correction in the complementary subspace to the
  definition of $ \varrho^{ \pm } $; we set
  \begin{align}
    \varrho^{ \pm }(t,v) 
    &:= \rho_{ \infty }^{ \pm }(v) + (tw)^{ -1 } r^{ \pm }(t,v),
    \quad r^{ \pm }(t,v) \in V_{ v }^{ \mp }
  \end{align} 
  chosen so that
  \begin{align}
    R_{ 1 }^{ \pm } + R_{ 3 }^{ \pm } + R_{ 4 }^{ \pm } - \phiap\psiap^{\pm}
    &= O_{ L_{ x }^{ 2 }  } \left(
      t^{ -1-\delta }(1+\log t)
    \right).
  \end{align} 
  A direct computation shows that choosing
  \begin{align}
      r^{ \pm }(t,v)
    &:= \frac{iw}{2}
        (\gamma^j-v_j\gamma^0)\partial_{v_j}\rho_\infty^\pm
        \pm \frac{3i}{4}
        \rho_\infty^\pm
      \pm \frac{w}{2}\log(tw)
      (\gamma^j-v_j\gamma^0)
      (\partial_{v_j}\phi_\infty)\rho_\infty^\pm
  \end{align}
  yields
  \begin{nalign}\label{eq:error-residual}
    &R_{ 1 }^{ \pm } + R_{ 3 }^{ \pm } + R_{ 4 }^{ \pm } - \phiap\psiap^{\pm}\\
    &\hspace{6em}= \frac{ e^{ \pm itw\mp i\Theta } }{ (tw)^{ 5/2  } }\phi_{ \infty
    }\rho_{ \infty }^{ \pm } - \phiap\psiap^{\pm}
    +O_{L_v^2}\left(t^{-\frac72}(1+\log t)^2\right)\\
    &\hspace{6em}= -\left(
      \phiap - (tw)^{ -1 } \phi_{ \infty }
    \right) \psiap^{ \pm }
    - \frac{ e^{ \pm itw\mp i\Theta } }{ (tw)^{ \frac{7}{2}  } } \phi_{ \infty
    }r^{ \pm }
    +O_{L_v^2}\left(t^{-\frac72}(1+\log t)^2\right).
  \end{nalign} 
  
  Here we have used that $ \phiap = (\phiap - (tw)^{ -1 }\phi_{ \infty }) +
  (tw)^{ -1 }\phi_{ \infty }  $.
  For all practical purposes, one should think of the correction $ r^{ \pm } $
  as simply satisfying
  \begin{align}
    \lVert r^{ \pm }(t,v) \rVert_{ L_{ v }^{ 2 } }
    \lesssim_{ \eta } 
    \lVert w\partial_{ v } \rho_{ \infty }^{ \pm } \rVert _{ L_{ v }^{ 2 } } 
    + \lVert \rho_{ \infty }^{ \pm } \rVert_{ L_{ v }^{ 2 }  } 
    + \epsilon ^2 (1+\log t)\lVert w^{ -\eta }\rho_{ \infty }^{ \pm }
    \rVert_{ L_{ v }^{ 2 }  },
  \end{align} 
  
  for any $ \eta > 0  $.
  Thus, integrating \eqref{eq:error-residual} and using the difference bounds
  from Lemma {\ref{lemma:approx-background}} below,
  \begin{align}
    \lVert R_{ \psi } \rVert _{ L_{ x }^{ 2 } }
    &\lesssim \sum_\pm
      \left\lVert\left(\phiap-(tw)^{-1}\phi_\infty\right)
      \psiap^\pm\right\rVert_{L_x^2}\\
    &\qquad+t^{-2}(1+\log t)^2\sum_\pm\left(
      \left\lVert w^{-\frac72}\partial_v\rho_\infty^\pm
      \right\rVert_{L_v^2}
      +\left\lVert w^{-\frac32}\partial_v^2\rho_\infty^\pm
      \right\rVert_{L_v^2}\right),
  \end{align} 
  where the bound on $ R_{ 2 }^{ \pm } $ follows from noting that $ \breve
  \partial_{ t } \varrho^{ \pm } = O \left( t^{ -2 }(1+\log t) \right) $.

  Notice that, since $ \rho_{ \infty }^{ \pm }\in V^{ \pm } $ and $ r^{ \pm
  }\in V^{ \mp } $, 
  \begin{align}
    \langle\rho_\infty^+,\gamma^0\rho_\infty^-\rangle
    =\langle\rho_\infty^\pm,\gamma^0r^\pm\rangle
    =\langle r^+,\gamma^0r^-\rangle
    &=0.
  \end{align}
  Expanding $ R_{ \phi } $ in terms of $ \rho_{ \infty }^{ \pm } $ and $ r^{
  \pm } $ inside the light cone, we get
  \begin{align}
    R_\phi
    &=-\frac{2}{(tw)^4}\Re\left[
      e^{2i(tw-\Theta)}
      \left(
        \langle\rho_\infty^+,\gamma^0r^-\rangle
        +\langle r^+,\gamma^0\rho_\infty^-\rangle
      \right)\right]
    +\frac{1}{(tw)^5}
      \left(\lvert r^-\rvert_H^2-\lvert r^+\rvert_H^2\right).
  \end{align}
  Using the Sobolev embedding $ \dot H_{ x } ^{ \frac{1}{2}  }\hookrightarrow
  L_{ x } ^{ 3 }$, we estimate
  \begin{align}
    \lVert R_{ \phi } \rVert _{ L_{ x }^{ \frac{3}{2}  }  }
    &\lesssim t^{ -4 } \sum_{ \pm }
    \left( 
      \lVert w^{ -\frac{3}{2}  } \rho_{ \infty }^{ \pm } \rVert_{ L_{ x }^{ 3 }  } 
      \lVert w^{ -\frac{5}{2}  } r^{ \mp } \rVert _{ L_{ x }^{ 3 }  } 
    \right) 
    + t^{ -5 }
    \left( 
      \lVert w^{ -\frac{5}{2}  }r^{ + } \rVert_{ L_{ x }^{ 3 }  }^{ 2 } 
      +\lVert w^{ -\frac{5}{2}  }r^{ - } \rVert_{ L_{ x }^{ 3 }  }^{ 2 } 
    \right)\\
    &\lesssim \epsilon^2 t^{ -1-\delta-\nu }.
  \end{align} 
  Since $ R_{ \phi }(t) $ is supported on a ball $ B_{ t }  $ of
  radius $ t $, the low-frequency bound follows from the Sobolev embedding
  $\dot H_{ x } ^{ \frac{1}{2} + \nu } \hookrightarrow L_{ x } ^{
  \frac{3}{1-\nu} }$ and Hölder:
  \begin{align}
    \lVert R_{ \phi }(t) \rVert_{ \dot H_{ x } ^{ -\frac{1}{2} -\nu } }
    \lesssim \lVert R_{ \phi }(t) \rVert _{ L_{ x }^{ \frac{3}{2+\nu}  }  } 
    &\lesssim t^{ \nu } \lVert R_{ \phi }(t) \rVert_{ L_{ x }^{ \frac{3}{2} }  }
    \lesssim \epsilon^2 t^{ -1-\delta }.
  \end{align} 

  \vspace{5pt}
  \noindent
  \textit{The general case.} Since the precise value of $ k $ does not affect
  the structure of the argument below, we assume $ 1\leq k\leq 18 $ whenever we
  discuss energy-type bounds, and $ 1\leq k\leq 9 $ for the remaining ones.
  We directly compute the effect of the $ \Gamma $ vector fields on the
  remainder \eqref{eq:error-residual} obtained above.
  Using $\Omega_{0j}(tw)=0$, we compute 
  \begin{align}
    \Omega_{0j}\left(
      \frac{e^{\pm itw\mp i\Theta}}{t(tw)^{5/2}}
    \right)
    &=\left(-v_j\mp i\log(tw)\breve\Omega_{0j}\phi_\infty\right)
      \frac{e^{\pm itw\mp i\Theta}}{t(tw)^{5/2}},\\
    \partial_t(tw)&=w^{-1},
    \qquad\partial_{x_j}(tw)=-v_jw^{-1},\\
    t\breve\partial_t r^\pm
    &=\pm\frac w2(\gamma^j-v_j\gamma^0)
      (\partial_{v_j}\phi_\infty)\rho_\infty^\pm.
  \end{align}
  Thus, commutation between $ \Gamma $ and $ (-i\gamma^{ \mu }\partial_{ \mu }) $
  preserves time decay modulo harmless additional velocity derivatives,
  boundary weights, and logarithms.
  Consequently,
  \begin{align}
    \lVert\Gamma^kR_\psi\rVert_{L_x^2}
    &\lesssim_k\sum_\pm
      \left\lVert\Gamma^k\left\{
      \left(\phiap-(tw)^{-1}\phi_\infty\right)\psiap^\pm
    \right\}\right\rVert_{L_x^2}\\
    &\quad+t^{-2}(1+\log t)^{k+2}\sum_\pm\left(
      \left\lVert w^{-\frac72}\partial_v^{k+1}\rho_\infty^\pm
      \right\rVert_{L_v^2}
      +\left\lVert w^{-\frac32}\partial_v^{k+2}\rho_\infty^\pm
      \right\rVert_{L_v^2}\right).
  \end{align}
  For the wave residual, we again compute 
  \begin{align}
    \Omega_{0j}\left(\frac{e^{2i(tw-\Theta)}}{(tw)^4}\right)
    &=-2i\log(tw)(\breve\Omega_{0j}\phi_\infty)
      \frac{e^{2i(tw-\Theta)}}{(tw)^4},\\
    \partial_\mu\left(\frac{e^{2i(tw-\Theta)}}{(tw)^4}\right)
    &=\left[-4\frac{\partial_\mu(tw)}{tw}
      +2i\partial_\mu(tw)-2i\partial_\mu\Theta\right]
      \frac{e^{2i(tw-\Theta)}}{(tw)^4}.
  \end{align}
  The same weighted Sobolev estimates give
  \begin{align}
    \lVert\Gamma^kR_\phi\rVert_{L_x^{3/2}}
    &\lesssim_k t^{-2}(1+\log t)^{k+2}
      \sum_\pm\left(
      \left\lVert w^{-\frac72}\partial_v^{k+1}\rho_\infty^\pm
      \right\rVert_{L_v^2}^2
      +\left\lVert w^{-\frac32}\partial_v^{k+2}\rho_\infty^\pm
      \right\rVert_{L_v^2}^2\right).
  \end{align}
  The low-frequency bounds follow as in the case $k=0$.
\end{proof}

\begin{lemma}\label{lemma:approx-background}
  We have the following approximation bounds
  \begingroup
  \mathtoolsset{showonlyrefs=false}
  \begin{align}\label{eq:approx-background}
    \sum_{ \pm }\lVert \Gamma^{ k } \left\{ (\phit - (tw)^{ -1 }\phi_{ \infty
      }) \psiap^{ \pm } \right\}  \rVert_{ L_{ x }^{ 2 }  }  
    + \lVert \Gamma^{ k } (\phif\psiap^{ \pm }) \rVert_{ L_{ x }^{ 2 }  }
    &\lesssim \epsilon t^{ -1-\delta-\nu } (1+\log t)^{ k }
  \end{align}
  \endgroup
  for $ 0\leq k\leq 18 $.
\end{lemma}

\begin{proof}
  For convenience, we prove the estimate for $ \psiap = \psiap^{ + }  $; the
  estimate for $ \psiap^{ - } $ follows analogously.
  As in the proof of Lemma {\ref{lemma:background}} above, the estimate for general $ k $ is
  obtained via a direct computation, and so we simplify the notation by only
  considering the case $ k = 0$.

  We split
  \begin{align}
    \phiap - (tw)^{ -1 }\phi_{ \infty } = (\phit - (tw)^{ -1 }\phi_{
    \infty } ) +  \phif
  \end{align}
  and estimate the contribution of each one of these
  terms separately.
  For the inhomogeneous term, we recall from {\eqref{eq:timelike-wave-agreement}} that $ \phit =  (tw)^{
  -1 }\phi_{ \infty } $ inside the cone and outside the conic shell $ Z = \left\{ 0\leq t-r \leq 2
  \right\} $.
  Moreover, on $ Z $, the $ L_{ x } ^{ 2 } $ bound on $ \psiap $ from
  {\eqref{eq:psiap-bound}} improves:
  \begin{align}
    \lVert w^{-\alpha}\psiap\rVert_{L_x^2}
    &\leq \lVert w^{-\alpha-\frac32}\rho_\infty^+\rVert_{L_v^2}
      +t^{-1}\lVert w^{-\alpha-\frac52}r^+(t)\rVert_{L_v^2}
    \lesssim \epsilon\bigl(1+t^{-1}(1+\log t)\bigr)
      \lesssim \epsilon,
  \end{align}
  for $ 0\leq \alpha < 2 $.
  Thus, fixing a parameter $ \alpha > 0$ to be specified later, we estimate
  \begin{align}
    \lVert (\phit - (tw)^{ -1 }\phi_{ \infty })\psiap \rVert_{ L_{ x }^{ 2 }  }
    &\lesssim \lVert w^{ \alpha }(\phit - (tw)^{ -1 }\phi_{ \infty }) \rVert_{
      L_{ x }^{ \infty } (Z)}
      \lVert w^{ -\alpha }\psiap \rVert _{ L_{ x }^{ 2 } (Z) }\\
    &\lesssim \epsilon^2 t^{ -1-\frac{\alpha}{2} } \lVert w^{ -\alpha }\psiap
    \rVert _{ L_{ x }^{ 2 } (Z) },
  \end{align} 
  since $ w^2\leq 4t^{ -1 } $ on $ Z $.
  We also used the bounds $t\lVert\phit\rVert_{L_x^\infty} +\lVert
  w^{-1}\phi_\infty\rVert_{L_v^\infty}\lesssim\epsilon^2$,
  which follow from a direct argument using the radial wave equation as in
  Proposition {\ref{prop:wave-inhom-bound}}.
  Then, choosing $ \alpha = 2(\delta + \nu)  $ furnishes
  {\eqref{eq:approx-background}}.

  As for the contribution from $ \phif\psiap $, we repeat the 
  localization argument using the weights $ w^{ \alpha } $:
  \begin{align}\label{eq:weighted-free-wave-product}
    \lVert\phif\psiap\rVert_{L_x^2}
    &\leq\lVert w^\alpha\phif\rVert_{L_x^\infty(r<t)}
      \lVert w^{-\alpha}\psiap\rVert_{L_x^2(r < t)}
      \lesssim\epsilon^2t^{-1-\delta-\nu}
  \end{align}
  where we have used the free wave bound from Lemma
  {\ref{lemma:wave-asymptotics}} and, again,
  that we may trade $ w^{ \alpha } $ for $ t $-decay near the cone.
  As before, choosing $ \alpha = 2(\delta + \nu) $ suffices.
  Note that the use of $ w $-weights is crucial, since the standard free-wave decay
  rate $ \lvert \phif \rvert \sim t^{ -1 } $ is insufficient for establishing
  {\eqref{eq:weighted-free-wave-product}}.

\end{proof}

\appendix
\addtocontents{toc}{\protect\setcounter{tocdepth}{1}}
\section{Analysis on hyperboloids}\label{sec:hyperboloids}
In this section we collect the hyperbolic identities used throughout the paper,
most notably in section \ref{sec:wp}.

In the interior of the forward cone, that is, $\{t>\lvert x\rvert\}$, we set
\begin{align}
  \Lambda := \left( t^2 - x^2 \right) ^{ \frac{1}{2}  },
  \quad
  r := \lvert x \rvert  ,
  \quad
  v := \frac{x}{t}   \, \left( \lvert v \rvert  < 1 \right),
\end{align}
and consider the hyperboloids
\begin{align}
  H_\Lambda&:=\left\{t^2-\lvert x\rvert^2=\Lambda^2,\ t>0\right\} 
\end{align}
with $ H := H_{ 1 } $.
These are equivalent to the hyperbolic polar coordinates
\begin{align}
  t
  &=\Lambda\cosh\rho
   =\frac{\Lambda}{\sqrt{1-\lvert v\rvert^2}},
  & r&=\Lambda\sinh\rho,
  & \lvert v\rvert&=\tanh\rho,
  & \omega&=\frac{x}{\lvert x\rvert}\in\mathbb S^2,
\end{align}
in $ H_{ \Lambda } $, where $\rho=\arctanh(|v|)$.  
Using these coordinates, we can rewrite the Minkowski metric $g$ inside the light
cone as
\begin{nalign}\label{eq:metricH}
  g
  &=-\dif t^2+\dif r^2+r^2\dif\omega^2\\
  &=-\dif\Lambda^2 +\Lambda^2\left(\dif\rho^2+\sinh^2\rho\dif\omega^2\right)\\
  &=-\dif\Lambda^2+\Lambda^2g_H.
\end{nalign}
That is, we foliate the interior of the light cone by hyperboloids $H_\Lambda$
with induced metric $g_{H_\Lambda}=\Lambda^2g_H$.
With the convention
$\square=\partial_t^2-\Delta_x=-g^{\mu\nu}\partial_\mu\partial_\nu$, the
wave operator takes the form
\begin{align}\label{eq:wave-operator-hyperbolic}
  \square
  &=\partial_\Lambda^2+3\Lambda^{-1}\partial_\Lambda
    -\Lambda^{-2}\Delta_\H.
\end{align}

It is also useful to look at $ H $ intrinsically rather than through its embedding into  $
\R^{ 4 }$.
We let 
\begin{align}
  \H := \left\{ v\in \R^{ 3 } : \lvert v \rvert  < 1 \right\} 
\end{align} 
denote the Klein--Beltrami hyperbolic disk model.
Its distance function is 
\begin{align}\label{eq:d_H}
  d_H(v,w)
  &:=\operatorname{arcosh}\left(
    \frac{1-v\cdot w}
    {\sqrt{(1-\lvert v\rvert^2)(1-\lvert w\rvert^2)}}\right).
\end{align}
In particular,
\begin{align}
  d_{ H } (v,0)
  &= \arctanh(|v|)
  = \frac{1}{2} \log \left( 
    \frac{ 1+\lvert v \rvert   }{ 1-\lvert v \rvert   }  
  \right).
\end{align} 

Intuitively, one should regard $\H$ as ``velocity space'' and $H$ as a ``reference
leaf'', akin to the hyperplane $\{t=0\}$ in the standard foliation by constant-time
slices.  
We often abuse notation and write $\H$ and $H$ interchangeably.  

The Lorentz fields $\Omega_{\mu\nu}$ are tangent to each $H_\Lambda$ and act as
the standard first-order derivatives on $H$.  
For a function $f(t,x)=f(v)$, their expressions in velocity coordinates are
\begin{align}
  \Omega_{0j}f
  &=x_j\partial_tf+t\partial_{x_j}f
   =\partial_{v_j}f-v_j(v\cdot\nabla_vf), \\
  \Omega_{ij}f
  &=v_j\Omega_{0i}f-v_i\Omega_{0j}f
   =v_j\partial_{v_i}f-v_i\partial_{v_j}f.
\end{align}
In particular, the boosts reduce to the coordinate derivatives at $v=0$ and so,
by Lorentz invariance, give normalized directional derivatives around any point
of $H$.

\subsection{Changes of variables and integration}\label{subsec:metric}

We use several coordinate systems throughout the text, and so we record
their respective changes of variables below.
On a fixed hyperboloid $ H_{ \Lambda } $, let $ \dif S_{ \Lambda } $ and $ \dif
V_{ H }  $ denote the restrictions of the Riemannian and
Lorentzian Euclidean measures on $ \R^{ 1+3 } $ to $ H_{ \Lambda } $, respectively.
Similarly, we can also parametrize $ H_{ \Lambda } $ by $ v\in \H $ via the
identification
\begin{align}
  \H \ni v \mapsto \frac{\Lambda}{\sqrt{ 1-v^2 } } (1,v) \in H_{ \Lambda }.
\end{align}
Then, one may directly check that
\begin{align}\label{eq:dVH-velocity}
  \dif V_{ H_{ \Lambda } } &= \Lambda^{ 3 }(1-v^2)^{ -2 }\dif v,
  \qquad
  \dif S_{ \Lambda } = \Lambda ^2 \sqrt{ t^2+x^2 } \dif V_{ H }. 
\end{align}
Here, $ \dif V_{ H }  $ should be thought of as the natural Lorentz-invariant
measure on $ \H $.

Below we record some volume computations for truncated hyperboloids.
For $R\geq2$, write
\begin{align}
  \H_{<R}
  &:=\left\{v\in\H:(1-\lvert v\rvert^2)^{-\frac12}<R\right\},
  & \H_R
  &:=\left\{v\in\H:(1-\lvert v\rvert^2)^{-\frac12}\sim R\right\}.
\end{align}
Since the time component of $\iota_\Lambda(v)$ is
$\Lambda(1-\lvert v\rvert^2)^{-\frac12}$, we have
\begin{align}
  \iota_\Lambda^{-1}\left(H_\Lambda\cap\{t<T\}\right)
  &=\H_{<T/\Lambda},
  & \iota_\Lambda^{-1}\left(H_\Lambda\cap\{t\sim T\}\right)
  &=\H_{T/\Lambda}.
\end{align}
Moreover, setting $r=\lvert v\rvert$ and
$s=(1-r^2)^{-\frac12}$ gives
\begin{align}
  \vol_{\dif V_H}(\H_{<R})
  &=4\pi\int_0^{\sqrt{1-R^{-2}}}\frac{r^2}{(1-r^2)^2}\dif r
    =4\pi\int_1^R\sqrt{s^2-1}\dif s
    \sim R^2.
\end{align}
Consequently, for $T\geq2\Lambda$,
\begin{align}\label{eq:truncated-hyperbolic-volume}
  \vol_{\dif V_H}\left(H_\Lambda\cap\{t<T\}\right)
  &=\vol_{\dif V_H}(\H_{<T/\Lambda})
    \sim T^2\Lambda^{-2},
\end{align}
where the volume of a subset of $H_\Lambda$ denotes that of its pullback to
$\H$.

The following local identities are used in the wave packet construction.
\begin{lemma}\label{lemma:metric-expansion}
  Let $v,w\in \H $, and let $ (t,x) \in \R_{ > 0 } \times \R^{ 3 } $ satisfy
  $t>\lvert x\rvert$.
  Then, for $\lvert w\rvert\ll1$,
  \begin{align}
    d_H^2(0,w)
    &=\lvert w\rvert^2+\frac23\lvert w\rvert^4
      +O(\lvert w\rvert^6).
  \end{align}
  Moreover, denoting  $ d_{ H }  = d_{ H } (x /t,v) $, we get
  \begin{align}
    g^{\mu\nu}\partial_\mu d_H^2\partial_\nu d_H^2
    &=4\Lambda^{-2}d_H^2, \notag\\
    -\square d_H^2
    &=2\Lambda^{-2}\left(1+2d_H\coth d_H\right).
    \label{eq:distance-contractions}
  \end{align}
  In particular, if $\Lambda d_H^2 = O(1)$, then
  \begin{align}
    -\square d_H^2
    &=6\Lambda^{-2}+O(\Lambda^{-3}).
  \end{align}
\end{lemma}

\begin{proof}
  The first identity follows by expanding
  $d_H(0,w)=\operatorname{artanh}\lvert w\rvert$.  By Lorentz invariance, for
  the remaining identities we may take $v=0$, so $d_H(x/t,0)=\rho$.  
  Then,
  \begin{align}
    g^{-1}(\nabla\rho^2,\nabla\rho^2)
    &=4\Lambda^{-2}\rho^2, \\
    -\square\rho^2
    &=\Lambda^{-2}\left(\partial_\rho^2
      +2\coth\rho\,\partial_\rho\right)\rho^2
     =2\Lambda^{-2}(1+2\rho\coth\rho),
  \end{align}
  where the last claim follows from
  $\rho\coth\rho=1+\frac13\rho^2+O(\rho^4)$.
\end{proof}

Let $\Lambda_1>\Lambda_0>0$ and let $H_{[\Lambda_0,\Lambda_1]}$ denote the
region between $H_{\Lambda_0}$ and $H_{\Lambda_1}$.  
For a compactly supported smooth function $f:\R^{1+3}\to\C$, let
\begin{align}
  F(\Lambda_1)
  &= \int_{H_{[\Lambda_0,\Lambda_1]}}f(t,x)\dif x\dif t.
\end{align}
Below we compute $ F'(\Lambda_0) $, which is needed in section \ref{sec:wp}.
Using the first line of \eqref{eq:dVH-velocity} and the Fundamental Theorem of
Calculus, 
\begin{align}
  F'(\Lambda_0)
  &=\lim_{\Lambda_1\to\Lambda_0}\frac{1}{\Lambda_1-\Lambda_0}
    \int_{\Lambda_0}^{\Lambda_1}\int_{\R^3}
    f\left(\sqrt{\Lambda^2+\lvert x\rvert^2},x\right)
    \frac{\Lambda}{\sqrt{\Lambda^2+\lvert x\rvert^2}}
    \dif x\dif\Lambda \notag\\
  &=\int_{H_{\Lambda_0}}f\,\dif V_{H_{\Lambda_0}}\\
  &=\Lambda_0^3\int_Hf(\Lambda_0t,\Lambda_0x)\dif V_H.
  \label{eq:derivH}
\end{align}

We also record the spinorial integration-by-parts identity used in wave packet
testing.  For compactly supported smooth spinors $\psi$ and $\psi'$ on $H$,
\begin{align}
  \int_H\langle\gamma^0\psi,\widehat\Omega_{\mu\nu}\psi'\rangle\dif V_H
  &=-\int_H\langle\gamma^0\widehat\Omega_{\mu\nu}\psi,\psi'\rangle\dif V_H.
  \label{eq:ibp}
\end{align}
Indeed, $\Omega_{\mu\nu}$ is a divergence-free Killing field on $H$.  If
$S_{\mu\nu}=\frac12\gamma_\mu\gamma_\nu$ denotes its spinorial correction,
then \eqref{eq:gamma-adjoint-relations} and \eqref{eq:gamma-relations} give
\begin{align}
  S_{\mu\nu}^\dagger\gamma^0+\gamma^0S_{\mu\nu}
  &=0.
\end{align}
Combining these two properties furnishes \eqref{eq:ibp}.

\subsection{Localized averages in hyperboloidal coordinates}
\label{app:hyperbolic-average}

We compare hyperbolic localized averages with flat ones in geodesic normal
coordinates.  Fix $v\in H$, choose an orthonormal identification
$T_vH\simeq\R^3$, and write $w=\exp_v(z)$.
The result below adapts Corollary \ref{coro:rescaledaverage} to the geometry of 
$ H $.

\begin{lemma}
  \label{lem:hyperbolic-flat-average}
  \label{coro:hyperbolic-flat-average-lambda}
  Let $\chi\in C_c^\infty(\R)$ and $f\in L^\infty(H)$.  For
  $0<\tau\leq1$, one has
  \begin{align}
    \tau^{-3}\int_H
      \chi\left(\frac{d_H(w,v)^2}{\tau^2}\right)f(w)\dif V_H(w)
    &=\tau^{-3}\int_{\R^3}
      \chi\left(\frac{\lvert z\rvert^2}{\tau^2}\right)
      f\bigl(\exp_v(z)\bigr)\dif z\\
    &\hspace{9em}+O\left(\tau^2\lVert f\rVert_{L^\infty(H)}\right).
  \end{align}
  In particular, for $\Lambda\geq1$,
  \begin{align}
    \Lambda^{\frac32}\int_H
      \chi\bigl(\Lambda d_H(w,v)^2\bigr)f(w)\dif V_H(w)
    &=\Lambda^{\frac32}\int_{\R^3}
      \chi\bigl(\Lambda\lvert z\rvert^2\bigr)
      f\bigl(\exp_v(z)\bigr)\dif z\\
    &\hspace{9em}+O\left(\Lambda^{-1}\lVert f\rVert_{L^\infty(H)}\right).
  \end{align}
  The implicit constants depend only on $\chi$.
\end{lemma}

\begin{proof}
  In geodesic polar coordinates centered at $v$, write
  $r=d_H(w,v)$.  Since $H$ has constant curvature $-1$, its metric and volume
  form are
  \begin{align}
    g_H
    &=\dif r^2+\sinh^2r\,\dif\omega^2,
    & \dif V_H
    &=\sinh^2r\,\dif r\dif\omega.
  \end{align}
  Under the orthonormal identification $T_vH\simeq\R^3$, we have
  $z=r\omega$ and $\dif z=r^2\dif r\dif\omega$.  Thus the defining property
  of normal coordinates and the comparison of the two polar volume forms give
  \begin{align}
    d_H(w,v)
    &=\lvert z\rvert,
    & \dif V_H(w)
    &=\left(\frac{\sinh\lvert z\rvert}{\lvert z\rvert}\right)^2\dif z
     =\left(1+\frac13\lvert z\rvert^2
       +O(\lvert z\rvert^4)\right)\dif z.
  \end{align}
  Moreover, the exponential map is a global diffeomorphism on $H$.  
  Changing variables by $w=\exp_v(z)$ therefore yields
  \begin{align}
    &\tau^{-3}\int_H
      \chi\left(\frac{d_H(w,v)^2}{\tau^2}\right)f(w)\dif V_H(w)\\
    &\qquad=\tau^{-3}\int_{\R^3}
      \chi\left(\frac{\lvert z\rvert^2}{\tau^2}\right)
      f\bigl(\exp_v(z)\bigr)
      \left(\frac{\sinh\lvert z\rvert}{\lvert z\rvert}\right)^2\dif z.
  \end{align}
  Since $\chi$ is compactly supported, the last integral is restricted to
  $\lvert z\rvert\lesssim\tau$.  Subtracting the flat average and using the
  quadratic error in the volume form, the absolute value of the difference is
  bounded by
  \begin{align}
    \tau^{-3}\lVert f\rVert_{L^\infty(H)}
    \int_{\lvert z\rvert\lesssim\tau}\lvert z\rvert^2\dif z
    &\lesssim\tau^2\lVert f\rVert_{L^\infty(H)}.
  \end{align}
  Taking $\tau=\Lambda^{-\frac12}$ gives the second identity.
\end{proof}

\section{Sobolev inequalities}\label{sec:sobolev}
\begin{prop}\label{thm:trace}
  Let $I\subset \R$ be an interval of length $S$, and let $f\in W^{1,1}(I)$ be a nonnegative
  function. Then, 
  \[
    f(s_0)\leq \frac{1}{S}\int_{I}f(s)\dif s + \int_{I}\lvert f'(s)\rvert \dif s,
  \] 
  for all $s_0\in I$.
\end{prop}

\begin{proof}
  The proof follows directly from the fundamental theorem of calculus
  \[
    \left\lvert f(s_0)-\frac{1}{S}\int_{I}f\right\rvert
    \leq \frac{1}{S}\int_{I}\lvert f(s_0) - f(s)\rvert \dif s
    \leq \frac{1}{S}\int_{I}\int_{\min\{s_0,s\}}^{\max\{s_0,s\}}\lvert f'(\tau)\rvert
    \dif\tau\dif s
    \leq \int_{I}\lvert f'(\tau)\rvert \dif\tau.
  \] 
\end{proof}

In particular, repeating the same argument with $ f^2 $ instead of $ f $ and
applying Hölder's inequality gives us
\begin{align}\label{eq:traceproduct}
  f^2(s_0) \lesssim \lVert f \rVert_{ L^{ 2 } (I) } \lVert f' \rVert  _{ L^{ 2 } (I) } 
  + S^{ -1 } \lVert f \rVert  _{ L^{ 2 } (I) } ^2.
\end{align} 

\begin{prop}
  Let $\chi\in C_{c}^{\infty}(\R^{3})$ be a radial function with  $\int \chi = 1$,
  and let  $g\in \dot H^{2}(\R^{3})\cap C(\R^{3})$.
  Then,
  \begin{equation}\label{eq:radialtest}
    \left\lvert g(0) - \int \chi(x)g(x)\dif x \right\rvert
    \lesssim \lVert g \rVert_{\dot H^{2}}.
  \end{equation} 
\end{prop}

\begin{proof}
  Without loss of generality, assume that $g$ is smooth and compactly supported.
  As such, one should view the left-hand side of \eqref{eq:radialtest} as the distribution
  $\left(\delta_{0}-\chi\right)$ acting on $g$. 
  Therefore,  
  \begin{align}
    \lVert \delta_{0}-\chi \rVert_{\dot H^{-2}}
    = \lVert  | \xi |^{-2} (1- \widehat \chi) \rVert_{L^{2}_{\xi}}.
  \end{align}
  By Taylor-expanding $\widehat \chi$, we get
  \[
    \widehat \chi (\xi) 
    =  \widehat\chi (0) + \nabla_{\xi}\widehat\chi(0) \cdot \xi
      + O\left(\lvert \xi \rvert^2\right)
    =  1 + O\left(\lvert \xi \rvert^2\right),
  \] 
  where since $\chi$ is radial, $\int x\chi(x)\dif x = 0$.
  Therefore, $\frac{(1-\widehat \chi)}{\lvert \xi \rvert^2}$ is bounded near $\xi = 0$,
  and since $\lvert \xi \rvert^{-2}$ is square-integrable for large $\xi$, we
  conclude that $\lVert \delta_{0}-\chi \rVert_{\dot H^{-2}}$ is finite.
\end{proof}

By rescaling the average in the previous result, we get the following corollary.

\begin{coro}\label{coro:rescaledaverage}
  Let $\chi\in C_{c}^{\infty}(\R^{3})$ be a radial function with  $\int \chi = 1$,
  and let  $g\in \dot H^{2}(\R^{3})\cap C(\R^{3})$.
  For $\tau>0$, define
  \[
    \chi_{\tau}(x) := \tau^{-3}\chi\left(\frac{x}{\tau}\right).
  \] 
  Then,
  \begin{align}
    \left\lvert g(0) - \int \chi_{\tau}(x)g(x)\dif x \right\rvert
    \lesssim \tau^{\frac{1}{2}}\lVert g \rVert_{\dot H^{2}}.
  \end{align} 
\end{coro}

\begin{lemma}[Embeddings of Hölder spaces]\label{lemma:holder-tradeoff}
  Let $ \alpha,\beta,\gamma\in (0,1) $ with
  $ \gamma>\max\{\beta,\frac{1}{2}\} $, and let $ f:\H\to \C^{ 4 } $ be a
  continuous function. Assume $ \lvert f(v) \rvert  \lesssim
  (1-v^2)^{ \alpha } $ and
  \begin{align}
    (1-v^2)^{ -\beta }f \in C^{ \gamma }(\H). 
  \end{align} 
  If $\alpha\geq\max\left\{\frac{1}{2},
  \frac{\gamma-\beta}{2\gamma-1}\right\},$ then $ f \in C^{
  \frac{1}{2}}(\overline{B(0,1)}) $ with
  $ f|_{ \partial B } = 0$.
\end{lemma}

\begin{proof}
  Since $ |f(v)|\lesssim (1-v^2)^{ \alpha } $, the zero extension of $ f $ is
  continuous, and for each $ \lvert \zeta \rvert = 1  $,
   \begin{align}
     \lvert f(v) - f(\zeta) \rvert
     =  \lvert f(v) \rvert
     \lesssim (1-v^2)^{ \alpha }
     \lesssim |\zeta - v|^{ \alpha }
     \lesssim |\zeta - v|^{ \frac{1}{2}  },
   \end{align}
  so $ f $ is $ \frac{1}{2}  $-Hölder continuous on $ \left\{ \lvert \zeta \rvert
  =1 \right\} $.
  It remains to compare interior points $ v,w\in B(0,1) $.

  If $ \lvert v-w \rvert  \gtrsim \min\left\{ 1-v^2,1-w^2 \right\}  $, the
  triangle inequality gives
  \begin{align}
    \lvert f(v) - f(w) \rvert  
    &\lesssim \lvert f(v) \rvert  + \lvert f(w) \rvert 
    \lesssim (1-v^2)^{ \alpha } + (1-w^2)^{ \alpha }
    \lesssim \lvert v-w \rvert^{ \alpha }
    \lesssim \lvert v-w \rvert^{ \frac{1}{2}  }.
  \end{align} 

  Now, if $ \lvert v-w \rvert  \ll \min \left\{ 1-v^2, 1-w^2 \right\}  $, we
  instead use the small-scale Euclidean bound
  \begin{align}
    \sinh^2 d_{ H }(v,w)
    &=\frac{\lvert v-w\rvert^2
      -\left(\lvert v\rvert^2\lvert w\rvert^2-(v\cdot w)^2\right)}
      {(1-\lvert v\rvert^2)(1-\lvert w\rvert^2)}
    \lesssim \frac{\lvert v-w\rvert^2}{(1-v^2)^2},
  \end{align}
  where we used $1-w^2\sim 1-v^2$ in this regime.  Thus
  $d_{ H }(v,w)\lesssim \lvert v-w\rvert/(1-v^2)$, and we expand
  \begin{align}
    \lvert f(v) - f(w)  \rvert 
    &\lesssim (1-v^2)^{ \beta }  
      \lvert (1-v^2)^{ -\beta }f(v) - (1-\lvert w \rvert  ^2)^{ -\beta }f(w) \rvert \\
    &\hspace{5em} + \lvert (1-v^2)^{ \beta } - (1-w^2)^{ \beta } \rvert 
      |(1-w^2)^{ -\beta }f(w)|\\
    &\lesssim (1-v^2)^{ \beta } d_{ H } (v,w)^{ \gamma } 
      + \lvert (1-v^2)^{ \beta } - (1-w^2)^{ \beta } \rvert
      (1-w^2)^{ \alpha-\beta } \\
    &\lesssim (1-v^2)^{ \beta - \gamma } \lvert v-w \rvert  ^{ \gamma }
      + (1-v^2)^{ \alpha - 1 } \lvert v-w \rvert.
  \end{align}
  Since $ \lvert v-w \rvert \ll 1-v^2 $ and $ \alpha\geq\frac{1}{2} $, the
  second term is bounded by $ \lvert v-w \rvert^{ \alpha }\lesssim
  \lvert v-w \rvert^{ \frac{1}{2} } $.  Combining this with the trivial bound
  $ \lvert f(v)-f(w)\rvert\lesssim (1-v^2)^{ \alpha } $, we obtain
  \begin{align}
    \lvert f(v)-f(w)\rvert
    &\lesssim \lvert v-w\rvert^{ \frac{1}{2} }
      +\min\left\{(1-v^2)^{ \alpha },
        (1-v^2)^{ \beta-\gamma }\lvert v-w\rvert^{ \gamma }\right\}\\
    &\lesssim \lvert v-w\rvert^{ \frac{1}{2} }
      +\lvert v-w\rvert^{ \frac{\alpha\gamma}{\alpha-\beta+\gamma} }
    \lesssim \lvert v-w\rvert^{ \frac{1}{2} },
  \end{align}
  where the last inequality follows from
  $ \alpha(2\gamma-1)\geq\gamma-\beta $, completing the proof.  
   
\end{proof}

\begin{lemma}[Weighted Morrey inequalities]\label{lemma:morrey}
  Let $ \alpha\in(0,1) $ and $ s>\frac{3}{2}+\alpha $.
  For a smooth scalar function $ g:\H\to\C $, we have
  \begin{align}\label{eq:g}
    \lVert g \rVert_{ C^{\alpha}(\H) }
    &\lesssim_{s,\alpha}
    \lVert (1-v^2)^{\frac{1}{2}}\Omega^{\leq s}g
      \rVert_{ L^2(\H) }.
  \end{align}
  For a smooth spinor $ f:\H\to\C^{4} $, we have
  \begin{align}\label{eq:f}
    \lVert f \rVert_{ C^{\alpha}(\H) }
    &\lesssim_{s,\alpha}
    \lVert (1-v^2)^{\frac{1}{4}}\widehat\Omega^{\leq s}f
      \rVert_{ L^2(\H) }.
  \end{align}
  Here the Hölder norms use the hyperbolic distance $ d_H $ and
  the Euclidean norm on the target, while the spinorial $ L^2 $
  norm uses $ \left\lvert \,\cdot\, \right\rvert_H $. 
\end{lemma}

\begin{proof}
  The difference in boundary decay between \eqref{eq:g} and \eqref{eq:f} 
  simply tracks the degeneracy of the $ H $-norm: for fixed $ v \in \R^{ 3 } $ with $
  \lvert v \rvert
  < 1$, we have $ \lvert z \rvert  \lesssim(1-v^2)^{ -\frac{1}{4}  }\lvert z \rvert
  _{ H } $,  $ z\in \C^{ 4 } $.
  As such, it is enough to consider the case of scalar functions.

  Let $ R \geq 1 $ be a dyadic number, and discretize the problem to $ \H_{ R }
  $ with $ \H_{ 1 } := \H_{ \lesssim 1 } $.
  By setting $ \varphi: \H \longrightarrow \R^{3} $,
  $ \varphi(v)=v/(R\sqrt{1-v^2}) $, we rescale
  \begin{align}
    g_{ R }(x) &= g\bigl(\varphi^{-1}(x)\bigr)
      =g\left(\frac{x}{\sqrt{R^{-2}+|x|^2}}\right)
  \end{align} 
  and use the Euclidean Morrey inequality, together with
  $ |\varphi(v)-\varphi(w)|\lesssim d_H(v,w) $ on $ \H_R $, to obtain
  \begin{align}
    \lVert g \rVert _{ C^{ \alpha }(\H_R) }
    &\lesssim \lVert g_{ R }  \rVert_{ C^{ \alpha }(\varphi(\H_R)) }
    \lesssim \lVert g_{ R }  \rVert _{ H^{ s }(\varphi(\H_R)) }.
  \end{align} 
  In these coordinates, $ t=\sqrt{R^{-2}+|x|^2}\sim1 $ and
  \begin{align}
    \dif V_H &= R^2 t^{ -1 }\dif x,
    \qquad 
    (\Omega_{0j}g)\circ\varphi^{-1} = t\partial_{x_j}g_R.
  \end{align} 
  The result follows from noting that 
  \begin{align}
    \lVert g_{ R }  \rVert _{ H^{ s }(\varphi(\H_R)) }
    &\lesssim_s R^{-1} \lVert \Omega^{ \leq s }g \rVert _{ L^{ 2 }(\H_{ R } ) },
  \end{align}
  and summing over the dyadic shells in $ \ell^{ 2 }_{ R }  $.

\end{proof}

\addtocontents{toc}{\protect\setcounter{tocdepth}{2}}

\end{document}